\documentclass[10pt,twoside,a4paper]{article} 
\usepackage{amsfonts}
\usepackage[utf8]{inputenc}
\usepackage[T1]{fontenc}
\usepackage{lmodern}
\usepackage[french]{babel}
\usepackage[all]{xy}
\usepackage{mathrsfs}
\usepackage{stixgras}
\RequirePackage{amsbsy}
\usepackage{xspace}

\usepackage[bookmarksopen=false,breaklinks=true,%
      backref=page,pagebackref=true,plainpages=false,%
      hyperindex=true,pdfstartview=FitH,%
      pdfpagelabels=true,colorlinks=true,linkcolor=blue,%
      citecolor=red,urlcolor=red,hypertexnames=false,colorlinks=false%
      ]%
   {hyperref}

\newcommand\hum[1]{} 

\newcounter{bidon}
\newcommand{\rdb}{\refstepcounter{bidon}}

\newcommand\Subsubsection[1]{
\rdb\addcontentsline{toc}{subsubsection}{#1}\subsubsection*{#1}}

\newtheorem{theorem}{Théorème}[section]   

\newtheorem{proposition}[theorem]{Proposition}
\newtheorem{lemma}[theorem]{Lemme}
\newtheorem{corollary}[theorem]{Corolaire} 
\newtheorem{propdef}[theorem]{Proposition et définition}
\newtheorem{remark}[theorem]{Remarque}

\newtheorem{fact}[theorem]{Fait}
\newtheorem{definition}[theorem]{Définition}
\newtheorem{defa}[theorem]{Définition alternative}

\newtheorem{notations}[theorem]{Notations}

\newcommand {\junk}[1]{}
\newcommand\scP{\mathscr{P}}

\newenvironment{Proof}[1]{
\trivlist \item[\hskip \labelsep{\it #1.}]\hskip 0pt\\}{\hfill 
\mbox{$\Box$}
\endtrivlist}

\newenvironment{proof}{
\trivlist \item[\hskip \labelsep{\it Démonstration.}]\hskip 0pt\\}{\hfill 
\mbox{$\Box$}
\endtrivlist}

\newcommand\oge{\leavevmode\raise.3ex\hbox{$\scriptscriptstyle\langle\!\langle\,$}}
\newcommand\feg{\leavevmode\raise.3ex\hbox{$\scriptscriptstyle\,\rangle\!\rangle$}}
\newcommand\gui[1]{\oge{#1}\feg}
\newcommand\som{\sum\nolimits}
\newcommand\bul{^\bullet }
\newcommand \dar { \,\downarrow\! }
\newcommand \uar { \,\uparrow\! }

\newcommand{\SCO}[6]{
\xymatrix @R = .4em @C =3em{
                                  &1 \ar@{-}[dl] \ar@{-}[dr] \\
#3 \ar@{-}[ddr]                   &   & #6 \ar@{-}[ddl] \\
                                  &\bullet\ar@{-}[d] \\
                                  &\bullet   \\
#2 \ar@{-}[ddr] \ar@{-}[uur]      &   & #5 \ar@{-}[ddl] \ar@{-}[uul] \\
                                  &\bullet \ar@{-}[d] \\
                                  &\bullet  \\
#1 \ar@{-}[uur]                   &   & #4 \ar@{-}[uul] \\
                                  & 0 \ar@{-}[ul] \ar@{-}[ur] \\
}
}

\newcommand{\SCOR}[8]{
\xymatrix @R =.4em @C =3em {
                                  &\bullet \ar@{-}[ddl] \ar@{-}[ddr] \\
                                  &#7\ar@{-}[u] \\
#3 \ar@{-}[ddr]                   &   & #6 \ar@{-}[ddl] \\
                                  &\bullet\ar@{-}[d] \\
                                  &\bullet   \\
#2 \ar@{-}[ddr] \ar@{-}[uur]      &   & #5 \ar@{-}[ddl] \ar@{-}[uul] \\
                                  &\bullet \ar@{-}[d] \\
                                  &\bullet  \\
#1 \ar@{-}[uur]                   &   & #4 \ar@{-}[uul] \\
                                  &#8\ar@{-}[d] \\
                                  &\bullet \ar@{-}[uul] \ar@{-}[uur] \\
}
}

\def\.@{\char'76}

\newcommand \wh {\widehat}

\def\={\hbox{~{\bf =}~}}
 
\def \aigu{\mathaccent19}    
\def \grave{\mathaccent18}    
\def \noi {\noindent}
\def \ss {\smallskip}
\def \sni {\ss\noi}
\def \ms {\medskip}
\def \mni {\ms\noi}

\def \alb {\allowbreak}

\def \im {\longrightarrow}

\def \aoe {{\land,\lor,\exists}}

\def \these {\quad\vdash\quad}
\def \thesu {\quad\vdash_\aoe\quad}

\def \cl {{\circ}}

\def\equidef{\buildrel{{\rm def}}\over{\quad\Longleftrightarrow\quad}} 

\def\gen#1{\left\langle{#1}\right\rangle} 
\def\avec#1{\,\&\left\{{#1}\right\}}

\newcommand \so[1] { \{#1\} }
\newcommand \sotq[2]{\so{#1\,\vert\,#2}}

\def \fait#1#2 {\vspace{.1cm} \begin{tabular}{p{2cm}p{2cm}l l l}  
$\cl\;#1\;\cl$ & $\these$ & $#2$ \end{tabular}\vspace{.1cm}}
\def \faut#1#2 {\vspace{.1cm} \begin{tabular}{p{2cm}p{2cm}l l l}
               $\cl\;#1\;\cl$ & $\thesu$ & $#2$
                  \end{tabular}}

\newcommand \ov[1] {\overline{#1} }

\newcommand \wi {\widetilde}
\def \sd#1#2 {\widetilde#1_#2 }

\newcommand \Z{\mathbb{Z}}
\newcommand \N{\mathbb{N}}

\newcommand \R{\mathbb{R}}

\newcommand \cC {{\cal C}}

\newcommand \cI {{\cal I}}

\newcommand \cN {{\cal N}}
\newcommand \cP {{\cal P}}
\newcommand \cQ {{\cal Q}}
\newcommand \cM {{\cal M}}
\newcommand \cT {{\cal T}}

\newcommand \cS {{\cal S}}

\newcommand \vu {\,\vee\,} 
\newcommand \vi {\,\wedge\,} 
\newcommand \Vu {\bigvee}
\newcommand \Vi {\bigwedge}
\newcommand \vda {\,\vdash\,}
\newcommand \vdi[1] {\,\vdash_{#1}\,}

\newcommand \Un {{\bf 1}}
\newcommand \Deux {{\bf 2}}
\newcommand \Trois {{\bf 3}}
\newcommand \Quatre {{\bf 4}}

\newcommand \Pfe {{{\rm P}_{{\rm fe}}}}

\newcommand \Hom {{\rm Hom}}

\renewcommand \dim {{\rm dim}}

\newcommand \Spec {{\rm Spec}}
\newcommand \Zar {\,{\rm Zar}}
\newcommand \Kr {\,{\rm Kr}}
\newcommand \Kru {\,{\rm Kru}}

\newcommand \num {{n$^{{\rm o}}$}}
 
\newcommand \hdr {hypo\-thèse de récur\-rence\xspace}
\newcommand \cad {c'est-à-dire\xspace}
\newcommand \ssi {si, et seu\-lement si,\xspace}

\newcommand \spdg {sans perte de géné\-ra\-lité\xspace}
\newcommand \Propeq {Les propri\-étés sui\-vantes sont 
équiva\-lentes.}
\newcommand \propeq {les propri\-étés sui\-vantes sont 
équiva\-lentes.}
\newcommand \disept {17$^{{\rm \grave{e}me}}$ problème de Hilbert\xspace}

\newcommand \alg {algèbre\xspace}
\newcommand \agB {algèbre de Boole\xspace}
\newcommand \agL {algèbre de Lindenbaum\xspace}

\newcommand \coli {combi\-naison liné\-aire\xspace}

\newcommand \com {comaximaux\xspace}

\newcommand \dfn {defi\-nition\xspace}
\newcommand \dfns {defi\-nitions\xspace}

\newcommand \elt {élément\xspace}
\newcommand \elts {éléments\xspace}

\newcommand \entrel {rela\-tion impli\-ca\-tive\xspace}
\newcommand \entrels {rela\-tions impli\-ca\-tives\xspace}

\newcommand \homo {homo\-mor\-phisme\xspace}
\newcommand \homos {homo\-mor\-phismes\xspace}

\newcommand \idemas {idéaux ma\-xi\-maux\xspace}

\newcommand \idep {idéal pre\-mier\xspace}
\newcommand \ideps {idéaux pre\-miers\xspace}

\newcommand \iso {isomor\-phisme\xspace}

\newcommand \itf {idéal \tf}
\newcommand \itfs {idéaux \tf}

\newcommand \mo {monoïde\xspace}
\newcommand \moco {\mos\com}

\newcommand \mos {monoïdes\xspace}

\newcommand \nst {Null\-stellen\-satz\xspace}
\newcommand \nsts {Null\-stellen\-sätze\xspace}

\newcommand \pf {de pré\-sen\-ta\-tion finie\xspace}

\newcommand \pols {poly\-nômes\xspace}

\newcommand \polcar {poly\-nôme carac\-té\-ris\-tique\xspace}

\newcommand \prt {pro\-pri\-\'et\'e\xspace}

\newcommand \proi {premier idéal\xspace}
\newcommand \prois {premiers idéaux\xspace}

\newcommand \proc {chaî\-ne idéale\xspace}
\newcommand \procs {chaî\-nes idéales\xspace}

\newcommand \proel {\proc élé\-men\-taire\xspace}
\newcommand \proels {\procs élé\-men\-taires\xspace}
\newcommand \proelo {\proel de lon\-gueur\xspace}

\newcommand \prolo {\proc de lon\-gueur\xspace}
\newcommand \prolos {\procs de lon\-gueur\xspace}

\newcommand \recu {récurrence\xspace}

\newcommand \sad {struc\-ture algé\-brique dyna\-mique\xspace}

\newcommand \tf {de type fini\xspace}
\newcommand \trdi {treil\-lis dis\-tri\-bu\-tif\xspace}
\newcommand \trdis {treil\-lis dis\-tri\-bu\-tifs\xspace}

\newcommand \cof {cons\-truc\-tif\xspace}
\newcommand \cofs {cons\-truc\-tifs\xspace}

\newcommand \cov {cons\-truc\-ti\-ve\xspace}
\newcommand \covs {cons\-truc\-ti\-ves\xspace}

\newcommand \coma {\maths\covs}
\newcommand \clama {\maths classiques\xspace}

\renewcommand \cot {cons\-truc\-ti\-ve\-ment\xspace}

\newcommand \LLPO{{\bf LLPO}}

\newcommand \maths {mathé\-ma\-tiques\xspace}

\newcommand \prco {démonstration \cov}
\newcommand \prcos {démonstrations \covs}

\newcommand \pte {principe du tiers exclu\xspace}

\newcommand \tcg {théo\-rème de compa\-cité\xspace} 
 
\newcommand \Tcgi {Le \tcg implique le résultat suivant.\xspace}

\RequirePackage{etoolbox}

\begin{document}
\title{Constructions cachées en algèbre abstraite \\ 
Dimension de Krull, Going Up, Going Down}
\author{
Thierry Coquand
(\thanks {~
Chalmers, University of Göteborg, Sweden,
email: coquand@cse.gu.se}~)
Henri Lombardi
(\thanks{~
Laboratoire de Mathématiques de Besançon, CNRS UMR 6623, UFR des Sciences et
Techniques,
Université Bourgogne Franche-Comté, 25030 BESANCON cedex, FRANCE,
email: henri.lombardi@univ-fcomte.fr}~),
}
\date{2018, \\version révisée d'un rapport technique de 2001
}
\maketitle

\begin{abstract}
Nous rappelons des versions \covs de la théorie de la dimension de Krull 
dans les anneaux commutatifs et dans les \trdis, dont les bases ont 
été posées par Joyal, Espan\~ol et les deux auteurs. Nous montrons 
sur les exemples de la dimension des algèbres \pf, du Going Up, du 
Going Down, 
\ldots\ que cela nous permet de donner une version \cov de 
théorèmes classiques importants, et par conséquent de récupérer un contenu 
calculatoire explicite lorsque ces théorèmes abstraits sont utilisés 
pour démontrer l'existence d'objets concrets. Nous pensons ainsi mettre 
en oeuvre une réalisation partielle du programme de Hilbert pour 
l'algèbre abstraite classique.  
\end{abstract}
\mni MSC 2000: 13C15, 03F65, 13A15, 13E05

\mni Mots clés :  Dimension de Krull, Going Up, Going Down,  
Mathématiques constructives.

\mni Key words: Krull dimension, Going Up, Going Down,  
Constructive Mathematics.

\vspace{2em} \centerline{\bf Avertissement}

\medskip  Cet article, basé sur un rapport technique de 2001, peut être vu comme complétant le chapitre XIII du livre \cite{ACMC}. 

Le rapport technique de 2001, disponible en \url{https://arxiv.org/abs/1712.04728v1}, était une version préliminaire de l'article publié \cite[Coquand \& Lombardi, 2003]{cl03} \url{}. 

L'article publié \cite{cl03}
a de larges parties communes avec le rapport technique de 2001, mais chacun des deux contient des thèmes non traités dans l'autre.

Par rapport au texte de 2001, nous avons mieux organisé l'ensemble, corrigé quelques erreurs et fait le lien avec le livre \cite{ACMC} quand cela nous a semblé utile.

La théorie \cov de la dimension de Krull via la notion de bord, telle qu'elle est exposée dans \cite{ACMC} est apparue en 2005 dans l'article \cite[Coquand, Lombardi, Roy]{clr05}. Dans le rapport technique de 2001 comme dans la note présente on développe une approche antérieure basée sur \cite[Joyal]{joy71}, \cite[Español]{esp}, \cite[Lombardi]{lom} et \cite[Cederquist, Coquand]{cc}.  

\newpage

\tableofcontents

\newpage
\markboth{Introduction}{Introduction}
\section*{Introduction} \label{sec Introduction}
Nous rappelons des versions \covs de la théorie de la dimension de Krull 
dans les anneaux commutatifs et dans les \trdis, dont les bases ont 
été posées par Joyal, Espan\~ol, et les deux auteurs 
(\cite{cc,esp,joy,lom}). Nous montrons sur les exemples du Going Up, du Going 
Down, de la dimension des variétés algébriques,  que cela nous 
permet de donner une version \cov des grands théorèmes classiques, et 
surtout de récupérer un contenu calculatoire explicite lorsque ces 
théorèmes abstraits sont utilisés pour démontrer l'existence 
d'objets concrets. Nous pensons ainsi mettre en oeuvre une réalisation 
partielle du programme de Hilbert pour l'algèbre abstraite classique. 

Notre exposé est dans le style des \coma à la Bishop 
(cf. en algèbre \cite{MRR}). Nous avons réduit au minimum l'appel aux 
notions de logique de manière à obtenir un texte dont la forme ne 
rebute pas les algébristes.
Lorsque nous disons que nous avons une version \cov d'un théorème 
d'algèbre abstraite, c'est que nous avons un théorème dont la démonstration 
est \cov, dont la signification calculatoire est claire, et à partir 
duquel nous pouvons récupérer le théorème classique correspondant 
par une utilisation immédiate d'un principe non \cof bien 
répertorié. Un théorème abstrait classique peut avoir plusieurs 
versions \covs intéressantes. 

Dans le cas des théorèmes d'algèbre abstraite classique, un principe 
non \cof qui permet de faire le travail en question est en général le 
\tcg, qui permet de construire un modèle pour une théorie formelle 
cohérente. Lorsqu'il s'applique à des structures algébriques de 
présentation énumérable (en un sens convenable) le \tcg est très voisin du $\LLPO$ de Bishop (tout nombre réel est $\geq 0$  
ou $\leq 0$). 

Notre volonté de ne pas utiliser le \tcg n'est pas au premier chef un 
choix philosophique mais un choix pratique. L'usage du \tcg conduit en 
effet à remplacer des démonstrations directes (cachées) par des démonstrations 
indirectes qui ne sont au fond qu'un double renversement par l'absurde de 
la démonstration directe, d'où la perte de son contenu calculatoire. Par 
exemple dans la démonstration abstraite du \disept on dit: si le polynôme $P$ 
n'était pas une somme de carrés de fractions rationnelles, on aurait 
un corps~$K$ dans lequel on décèlerait une absurdité en lisant la 
démonstration (\cov) que le polynôme est partout positif ou nul. La remise sur 
pied de cette démonstration abstraite est de dire: partons de la démonstration (\cov) 
que le polynôme est partout positif ou nul, et montrons (par les arguments 
explicites présents dans les démonstrations classiques) que la tentative de 
construire $K$ échoue à coup s\^ur. Cela nous donne la somme de 
carrés que nous cherchions. Entretemps il a fallu remplacer le 
théorème abstrait non \cof: \gui{tout corps réel peut être 
ordonné} par le théorème \cof: \gui{dans un corps qui refuse toute 
tentative d'être ordonné, $-1$  est une somme de carrés}. Et on 
passe du second au premier par le \tcg, tandis que la \prco du second est 
cachée dans les manipulations algébriques contenues dans la démonstration 
classique du premier (cf. \cite{lom98}).  

\mni Voici un plan rapide du papier.

\paragraph{Définition \cov de la dimension de Krull des anneaux 
commutatifs}~

\noindent Dans la section \ref{secKrA} nous donnons des démonstrations plus 
lisibles (pour le mathématicien classique) de résultats contenus dans \cite{lom}, lesquels y 
étaient démontrés en utilisant la notion de \sad introduite dans \cite[Coste, Lombardi, Roy]{clr}. 
La notion centrale qui est mise en oeuvre est celle de {\em 
spécification partielle pour une chaîne  d'idéaux premiers}. 
Les constructions abstraites de chaînes  d'idéaux premiers ont leur 
contrepartie \cov sous forme d'un théorème de collapsus simultané
(théorème \ref{ThColSimKrA}). Nous développons la notion de suite 
singulière, qui est une variation \cov autour de la notion   de suite 
régulière. Cette notion permet de caractériser \cot la dimension de Krull.
Nous montrons ainsi que la notion de dimension de Krull a un contenu 
calculatoire explicite en terme d'existence (ou non) de certains types 
d'identités algébriques. Ceci confirme l'adage selon lequel 
l'algèbre commutative admet une interprétation calculatoire sous forme de 
machineries de productions d'identités algébriques (certaines, parmi les 
plus fameuses, sont appelées des \nsts.) Enfin nous donnons une version 
\cov du théorème qui affirme que la dimension de Krull est la borne 
supérieure
de la dimension des localisés en les \idemas. 

\paragraph{Treillis  distributifs, \entrel et dimension de Krull}~

\noindent Dans la section \ref{secKrullTreil}  nous développons la 
théorie \cov de la dimension de Krull des \trdis,  dans le 
style de la section  \ref{secKrA}. Une grande simplification des démonstrations et des calculs est obtenue
en mettant systématiquement en avant la notion de \entrel (en anglais 
entailment relation), développée dans \cite{cc}, qui a son origine dans la règle de coupure du 
calcul des séquents de Gentzen
. 

\paragraph{La théorie de Joyal}~

\noindent Dans la section \ref{secJoyal} nous expliquons  la 
théorie de Joyal pour la dimension de Krull. Nous faisons le lien avec l'approche de la section \ref{secKrullTreil}. La méthode de Joyal donne une démonstration plus simple du résultat principal, à savoir que la définition constructive est équivalente en \clama à la définition classique. 
Nous faisons également le lien avec les
développements apportés par Espa\~nol à la théorie de Joyal.

\paragraph{Treillis de Zariski et de Krull dans un anneau commutatif}~

\noindent Dans la section \ref{secZariKrull}  nous introduisons  le
treillis de Zariski d'un anneau commutatif. Ses éléments peuvent être identifiés aux 
radicaux d'\itfs. Ce \trdi est la contrepartie \cov du spectre de Zariski: 
les points de ce dernier sont les \ideps du treillis de Zariski, et les 
parties constructibles du spectre de Zariski sont les éléments du 
treillis booléen engendré par le  treillis de Zariski.
L'idée de Joyal est de définir
la dimension de Krull d'un anneau commutatif comme étant celle de son
treillis de Zariski, ce qui évite tout recours aux \ideps.
Nous établissons ensuite l'équivalence entre le point de vue (\cof) 
de Joyal et le point de vue (\cof) donné à la section  \ref{secKrA} 
pour la dimension de Krull d'un anneau commutatif.

\paragraph{Dimension de Krull relative d'une extension}~

\noindent Nous donnons dans  la section \ref{secKrullRel} la version  \cov
pour la \emph{dimension de Krull relative} d'une extension d'anneaux
ou de \trdis.
Cet analogue \cof qui semble au premier abord un peu étrange 
(définition~\ref{defColRel}) n'a pas été tiré du chapeau, mais a 
été fourni comme un outil qui s'impose
de lui-même dans le déchiffrage des démonstrations classiques.
Le théorème \ref{th.nstformelRel} démontre l'équivalence (en \clama) avec la notion classique 
définie de manière usuelle et l'on obtient \cot le théorème principal
concernant la dimension de Krull d'une extension (théorème \ref{thDim1}).

\paragraph{Going Up et Going Down}~

\noindent  La section \ref{secGUGD} donne la définition \cov des 
notions classiques de Lying over, Going up et Going down, aussi bien pour les \trdis que pour les anneaux commutatifs. L'équivalence avec les définitions usuelles en \clama est établie. Comme application on obtient la version \cov  des 
 théorèmes de montée et de descente 
dans les extensions entières d'anneaux commutaitfs, ainsi que celle du Going down pour les 
extensions plates.
Nous montrons que ces théorèmes apparemment très abstraits (qui 
affirment l'existence de certains \ideps) ont une signification tout à
fait concrète en terme de machinerie de construction d'identités 
algébriques. Cette version concrète implique la version abstraite par 
une utilisation simple du \tcg. Et le plus important, la démonstration de la 
machinerie concrète est déjà à l'oeuvre de manière cachée dans 
la démonstration abstraite du théorème absrait.

\paragraph{Annexe}~

\noindent L'annexe explique le matériel logique de base pour le théorème de complétude, le \tcg, les théories géométriques, les théories dynamiques et
le jeu qui s'y déploie entre \clama et \coma. 

\paragraph{Conclusion}~

\noindent On notera que les théorèmes \cofs établis dans cet article
concernant la dimension des anneaux de polynômes, celle des algèbres \pf
sur un corps, le Going up et le Going down sont nouveaux (ils n'avaient pu
être obtenus dans le cadre de la théorie de Joyal tant qu'elle se
limitait à parler du treillis de Zariski des anneaux commutatifs sans aller voir de plus près
les calculs en jeu en termes d'identités algébriques).

\smallskip  Cet article constitue une confirmation de la possiblité concrète de 
réaliser le programme de Hilbert pour de larges pans de l'algèbre 
abstraite 
(cf.  \cite{cp,clr,kl,lom95,lom97,lom98,lom,lom99,lom99a,lq99}). 
L'idée générale est de remplacer les structures abstraites idéales 
par des spécifications partielles de ces structures.
La jolie et très courte démonstration abstraite qui utilise les objets idéaux 
a une contrepartie purement calculatoire au niveau des spécifications 
partielles de ces objets idéaux.
La plupart des théorèmes d'algèbre abstraite, dans la démonstration 
desquels l'axiome du choix et le principe du tiers exclu semblent offrir 
un obstacle sérieux à une interprétation explicite en termes de 
calculs algébriques, auraient ainsi une version \cov, à partir de 
laquelle l'utilisation du \tcg fournirait la version classique abstraite. 
Plus important encore, la démonstration abstraite idéale contiendrait toujours, 
de manière plus ou moins cachée, la \prco du théorème \cof 
correspondant.

\ms 
Signalons enfin qu'un traitement constructif du Principal ideal theorem de 
Krull est donné
dans l'article~\cite{cl}.

\patchcmd{\sectionmark}{\MakeUppercase}{}{}{}

\section{Définition \cov de la dimension de Krull des anneaux 
commutatifs} 
\label{secKrA}

Cette section est basée sur l'article \cite{lom}.

Soit un anneau commutatif $A$, on note $\gen{J}$  ou s'il y a 
ambigüité $\gen{J}_A$  l'idéal de $A$ engendré par la partie 
$J\subseteq A$. On note
$\cM(U)$ le \mo\footnote{Un \mo sera toujours un \mo multiplicatif.} 
engendré par la partie 
$U\subseteq A$.

\subsection{Premiers idéaux} 
\label{subsecProi}

\begin{definition}{\rm  
\label{defproi1} Soit un anneau commutatif $A$. 
\begin{itemize}
\item 
Une {\em spécification partielle pour un idéal premier} (ce que nous 
abrégeons en {\em \proi}, en anglais {\em idealistic prime}) est un 
couple $\cP=(J,U)$ de parties de $A$. 
\item Un \proi $\cP=(J,U)$ est dit {\em complet} si $J$ est un  idéal,  
$U$ est un \mo et $J+U=U$. 
\item Si $P$ est  un \idep et $S=A\setminus P$, on note $\wh{P}$ le \proi $(P,S)$.
\item Soient $\cP_1=(J_1,U_1)$ et $\cP_2=(J_2,U_2)$  deux \prois et $P$ un \idep.
\begin{itemize}
\item On  dit que {\em $\cP_1$ est contenu dans $\cP_2$} et on écrit 
$\cP_1\subseteq \cP_2$ si $J_1\subseteq J_2$ et  $U_2\subseteq U_1$.
\item  On  dit que  \emph{$\cP_2$ est un raffinement de $\cP_1$} et on 
écrit 
$\cP_1\leq \cP_2$ si $J_1\subseteq J_2$ et  $U_1\subseteq U_2$.
\item  On  dit que  \emph{l'\idep $P$ raffine le \proi $\cP$} si $\cP\leq \wh{P}$.
\end{itemize}
\item  On note $\kappa(A,\cP)$, ou $\kappa(\cP)$ si le contexte est clair, l'anneau  obtenu à partir de $A$ en forçant \gui{$j=0$} pour $j\in J$ et \gui{$u$ inversible} \hbox{pour $u\in U$}. Si $J$ est un idéal et $U$ un \mo c'est l'anneau $U^{-1}(A/J)$. 
\end{itemize}
}
\end{definition}

Nous regardons un \proi $\cP$  comme une 
spécification partielle pour une idéal premier
(au sens usuel)  $P$ qui raffine 
$\cP$.  

Un premier idéal $\cP=(J,U)$ engendre un premier idéal complet minimum pour la relation de raffinement, à savoir $(\gen{J},\gen{J}+\cM(U))$. 
On note $\cS(J,U)$ le \mo \hbox{$\gen{J}+\cM(U)$} et $A_\cP$ l'anneau $\cS(J,U)^{-1}A$.

\begin{definition} \label{defiproicollaps}{\rm
On dit que le \proi $\cP=(J,U)$ \emph{s'effondre}, ou \emph{collapse}, si l'idéal~$\gen{J}$ coupe le \mo $\cM(U)$, autrement dit si $0\in \gen{J}+\cM(U)$}. Il revient au même de dire que l'anneau~$\kappa(\cP)$ est trivial,
ou encore que l'anneau $A_\cP$ est trivial.
\end{definition}

Notez que le \proi $(0,1)$ s'effondre \ssi $1=_A0$.

\begin{lemma} \label{lemfinicol}
Un \proi qui s'effondre raffine un \proi fini qui s'effondre. 
\end{lemma}

\Subsubsection{Collapsus simultané et lemme de Krull}

Dans la suite nous noterons $(J\cup\so {x} , U)$ sous la forme abrégée $(J,x;U)$ et $(J, \so {x} \cup U)$ sous la forme $(J;x,U)$.

\begin{theorem} 
\label{ThColSim0}{\em (Collapsus simultané pour les \prois)}\\
Soit un  \proi $\cP=(J,U)$ dans un anneau commutatif $A$
et soit $x\in A$. Supposons que les \prois 
$(J,x;U)$ et $(J;x,U)$ 
s'effondrent tous les deux, alors $\cP$ s'effondre également.
\end{theorem}
\begin{proof}
La démonstration est le truc de Rabinovitch. A partir de deux 
égalités
$u_1+j_1+ax=0$ et $u_2x^m+j_2=0$ (avec $u_i\in \cM(U)$, $j_i\in \gen{J}$, 
$a\in A$,
$n\in\N$) on en fabrique une troisième, $u_3+j_3=0$ en éliminant $x$:
on obtient $u_2(u_1+j_1)^m+(-a)^mj_2=0$, avec $u_3=u_2u_1^m$. 
\end{proof}

\begin{corollary} 
\label{corNstformHilbert} {\em (Lemme de Krull ou \nst de Hilbert 
formel)} \\
Soit $\cP=(J,U)$  un \proi dans un anneau commutatif $A$. Le \tcg implique 
que  \propeq
\begin{enumerate}
\item Le premier idéal $\cP$ ne s'effondre pas.
\item Il existe un idéal premier détachable\footnote{Rappelons qu'une 
partie d'un ensemble est dite détachable lorsqu'on a un test pour 
l'appartenance à la partie. Par exemple les \itfs d'un anneau de 
polynômes à coefficients entiers sont détachables.} $Q$ qui raffine $\cP$.
\item Il existe un morphisme $\psi$ de $A$ vers un anneau intègre $B$ 
vérifiant $\psi(J)=0$ \hbox{et $0\notin\psi(U)$}.
\item Il existe un morphisme $\phi$ de $A$ vers un corps algébriquement clos $L$ 
vérifiant $\phi(J)=0$ \hbox{et $0\notin\phi(U)$}.
\end{enumerate}
\end{corollary}
\begin{proof} 
Les points \emph{2}, \emph{3} et \emph{4} sont clairement équivalents (au moins en \clama où tout corps possède une clôture algébrique\footnote{En fait cela résulte aussi du \tcg.}). Le point \emph{2} implique clairement le point \emph{1}. Voyons la réciproque.\\
Commençons par la démonstration classique usuelle qui s'appuie, non sur le \tcg, mais sur le 
\pte et le lemme de Zorn. On considère un \proi 
$\cQ=(P,S),$ qui est 
maximal, pour la relation de raffinement, parmi les \prois qui raffinent 
$\cP$ et qui ne s'effondrent pas (le lemme de Zorn s'applique grâce 
au lemme~\ref{lemfinicol}). Tout d'abord, il est clair que~$\cQ$ est 
complet, puisqu'on ne change pas le collapsus en complétant un \proi.
Si $\cQ$ n'était pas un idéal premier avec son 
complément, on aurait un $x\in A\setminus (P\cup S)$. 
Mais puisque le \proi $(P,S)$ ne s'effondre pas l'un des deux
\prois $(P,{x};S)$ et $(P;{x},S)$ doit ne pas s'effondrer
(théorème \ref{ThColSim0}). Et ceci contredit la maximalité.\\
Voyons la démonstration basée sur le \tcg (qui est moins fort que le principe du tiers exclu avec le lemme de Zorn).
Nous considérons pour cela une théorie propositionnelle qui décrit le \proi   
$(P,S)$ associé à un \idep $P$ dans $A$, avec les contraintes $J\subseteq P$ \hbox{et $U\subseteq S$}.
Dans cette théorie nous avons un prédicat 
pour $x\in P$ (on devrait le noter $\mathrm{P}(x)$, mais nous garderons $x\in P$ pour plus de lisibilité), un autre pour $x\in S$
et les axiomes sont les suivants ($P$ est un idéal, $S$ est un monoïde, $P\cap S=\emptyset$, $P\cup S=A$)
\begin{enumerate}{\it 
\item $a\in P$  pour chaque $a\in J$
\item $u\in S$  pour chaque $u\in U$
\item $0\in P$
\item $(a\in P\;\&\; b\in P)\Rightarrow\, a+b\in P$
\item $a\in P\Rightarrow ab\in P$
\item $1\in S$
\item $(a\in S\;\&\; b\in S)\Rightarrow ab\in S$, 
\item $(a\in P\;\&\; a\in S)\Rightarrow \bot$  (collapsus) 
\item $a\in P \hbox{ ou } a\in S$.}
\end{enumerate}
Pour la théorie avec les seuls axiomes \emph{1} à \emph{8} l'inconsistance est facilement équivalente au collapsus du \proi $(J,U)$. 
Le point \emph{1} de l'énoncé dit donc que cette théorie est consistante.
D'après le théorème~\ref{ThColSim0} la théorie  où l'on ajoute l'axiome \emph{9} est également consistante (voir l'annexe concernant les théories géométriques).
Par le \tcg, cette théorie a un modèle. Ce modèle nous donne l'\idep 
voulu.
\end{proof}

\Subsubsection{Saturation}

\begin{definition} 
\label{defSat1}{\rm
On dit que le  \proi $\cP=(J,U)$  est \emph{saturé} si:
\begin{itemize}
\item tout élément $x\in A$  tel que le \proi
$(J, {x} ; U)$ s'effondre est dans $U$,
\item tout élément $x\in A$  tel que le \proi
$(J; {x} , U)$ s'effondre est dans $J$. 
\end{itemize}}
\end{definition}

Par exemple pour un idéal premier détachable $P$, le \proi $\wh P$ est un \proi saturé. 

Lorsque le \proi $(J,U)$ est saturé  
on a les équivalences
$$ 1\in J\quad \Longleftrightarrow\quad  0\in U\quad 
\Longleftrightarrow\quad  (J,U)=(A,A).
$$

\begin{lemma} \label{lemSat1}~
\begin{enumerate}
\item Si $\cP=(J,U)$ est saturé, il est complet: $J$ est un idéal, $U$ est un \mo \hbox{et $J+U=U$}.
\item 
Un \proi $\cP=(J,U)$ engendre un \proi saturé minimum (pour la relation de raffinement)
$\cP_1=(I,V)$ obtenu comme suit.
\begin{itemize}
\item $I$ est l'ensemble des $x\in A$  tels que le \proi
$(J,x;U)$ s'effondre,
\item $V$ est l'ensemble des $x\in A$  tels que le \proi
$(J;x,U)$ s'effondre. 
\end{itemize}
\end{enumerate}

\end{lemma}
%
\begin{Proof}{Démonstration sans calcul} \emph{1}. 
Montrons que $J$ est un idéal. Soit $z\in \gen{J}$. Il suffit de montrer que $(J;z,U)$ s'effondre, i.e. que~$\gen{J}$ coupe $\cM(z,U)$, ce qui est évident.  
Montrons \hbox{que $J+U\subseteq U$}. Soit $a\in J$ et $u\in U$. Il suffit de montrer que $(J,a+u;U)$ s'effondre, i.e. que $\gen{J,a+u}$ coupe $\cM(U)$, ce qui est évident.

\sni \emph{2}.
Il suffit de montrer que $\cP_1$ est saturé. Soit par exemple un $x$
tel que $(I,x;V)$ s'effondre. On doit montrer que $x\in V$.
On a par exemple $(J,y_1,y_2,x;U,z)$ qui s'effondre avec $y_1,y_2\in I$ et $z\in V$. Autrement dit $(J,y_1,y_2,x;U,z)$, $(J;U,y_1)$,  $(J;U,y_2)$
et $(J,z;U)$ s'effondrent tous les quatre. Alors en appliquant plusieurs fois le collapsus simultané, on obtient bien que $(J,x;U)$ s'effondre, \cad que $x\in V$.  
\end{Proof}
%

\begin{corollary} 
\label{corNstfH2} 
Soit $\cP=(J,U)$  un \proi dans un anneau commutatif $A$ et
 $(I,V)$ le saturé de $\cP$. Le \tcg implique que $I$ est l'intersection 
des idéaux premiers contenant~$J$ et ne coupant pas $U$,
tandis que $V$ est l'intersection des complémentaires de ces mêmes 
idéaux premiers. En particulier le nilradical d'un idéal de $A$ est l'intersection des \ideps qui le contiennent.
\end{corollary}
%
\begin{Proof}{Démonstration sans calcul}
Par définition tout \idep $P$ qui raffine $\cP$ contient $I$ et son complémentaire contient~$V$. Soit maintenant un $x$ qui n'est pas dans $V$.
Le \proi $(J,x;U)$ ne s'effondre pas. Il existe donc (corolaire \ref{corNstformHilbert}) un \idep $P$ qui raffine $\cP$ et qui contient $x$. Cet $x$ n'est donc pas dans l'intersection des complémentaires des \ideps qui raffinent $\cP$.    
\end{Proof}
%

\begin{propdef}  
\label{defConjug}
Un idéal $I$  et un \mo $S$ sont dits {\em conjugués} lorsque \hbox{l'on a}:
\begin{itemize}
\item  $I+S\subseteq S$,
\item   $s  t\in S     \Rightarrow  s \in S$,
\item $(s a\in I, \; s\in S)   \Rightarrow  a\in I$,
\item   $a^n\in I       \Rightarrow  a\in I  \quad (n\in\N,n>0)$. 
\end{itemize}
Le \proi $(I,S)$ est saturé \ssi l'idéal $I$ et le \mo $S$ sont conjugués.
\end{propdef}

\begin{proof}L'implication directe est facile. La réciproque se fait sans calcul comme pour le point \emph{1} du lemme~\ref{lemSat1}.
\end{proof}
Ainsi le saturé du \proi $(0,1)$ est $(\cN,A^{\times })$ où $\cN$ est le 
nilradical de $A$ et $A^{\times}$ le groupe des unités.

\Subsubsection{Exemples de versions constructives de théorèmes classiques} 
Nous donnons maintenant des exemples de versions constructives de 
théorèmes usuels en \clama concernant les idéaux premiers.

Notre premier exemple concerne les anneaux arithmétiques, 
\cad les anneaux dont les idéaux de type fini sont localement 
principaux.
Par exemple un domaine de Prüfer ou de Dedekind est un anneau 
arithmétique. 

La \emph{définition constructive d'un anneau 
arithmétique $A$} est la suivante:  pour tous $x,y\in A$ on peut 
trouver un élément $s\in A$ tel que dans $A_s$ on a  
$\gen{x,y}=\gen{x}$
et dans  $A_{1-s}$ on a  $\gen{x,y}=\gen{y}$. Cela revient encore à dire 
qu'on peut trouver $s$, $t$, $v,$ $w$ tels que $sx=vy$, $ty=wx$ et $s+t=1$.

\begin{lemma} 
\label{lemAritCom}~ 
\\
\emph{\'Enoncé classique.} Dans un anneau arithmétique deux idéaux 
premiers incomparables sont comaximaux.\\
\emph{\'Enoncé constructif.} Dans un anneau arithmétique deux \prois 
\emph{incomparables}\footnote{Définition ci-après.} sont comaximaux.
\end{lemma}
\begin{proof}
Le substitut constructif pour deux idéaux premiers incomparables est 
donné par deux \prois $\cP=(I,U)$ et $\cQ=(J,V)$ qui sont \emph{incomparables} au 
sens suivant: on connaît deux éléments~$x$ et~$y$  de $A$ qui 
témoignent que $\cP$ et $\cQ$ ne peuvent pas être raffinés en des 
idéaux premiers comparables. 
Plus précisément, si~$\cP'=(I',U')$ et~$\cQ'=(J',V')$ sont les 
saturés de
$\cP$ et $\cQ$ \hbox{on a $x\in I'\cap V'$} \hbox{et $y\in J'\cap U'$}. 
La conclusion doit être que $I'+J'=\gen{1}$. La démonstration de l'énoncé constructif est immédiate:
on considère  $s$, $t$, $v,$ $w$ comme ci-dessus 
(définition constructive d'un anneau arithmétique),
puisque $x\in I'$, $y\in U'$,  $ty=wx$ et puisque $\cP'$ est saturé, 
on a $t\in I'$. Et de façon symétrique on a $s\in J'$. 
\end{proof}

Il est clair que le théorème classique découle de sa version 
constructive. Et que le théorème classique implique la version 
constructive en \clama. Mais la version constructive n'a pas besoin 
d'idéaux premiers et elle a une signifcation calculatoire claire.

\ms Notre deuxième exemple est le célèbre lemme d'évitement des 
idéaux premiers. Une version constructive naturelle (qui implique 
l'énoncé classique lorsque l'idéal $I$ est \tf) est le lemme \ref{lemEvit} suivant.
On peut trouver une autre version constructive dans \cite{MRR} (section II.2, théorème 2.3, exercices 12 et 13). 
Il est possible de démontrer que dans un anneau noethérien cohérent fortement 
discret, si $\cP$ est l'un des \prois obtenus en saturant une \proc finie
(voir la définition en section \ref{subsecProc}), 
alors  $\cP$ vérifie l'hypothèse de décidabilité requise pour les 
\prois dans le lemme~\ref{lemEvit}.
\begin{lemma} 
\label{lemEvit} {\rm (Lemme d'évitement des idéaux premiers)}.
Soit $I$ un \itf \hbox{et $\cP_1,\ldots,\cP_n$} des \prois 
d'un anneau $A$. Supposons qu'on sache tester le collapsus des 
raffinements
finis des $\cP_k$.
Alors on sait construire des raffinements finis $\cQ_1,\ldots,\cQ_n$  \hbox{de  
$\cP_1,\ldots,\cP_n$} tels que (en notant~\hbox{$\cQ_k=(I_k,V_k)$}):  
\begin{itemize}
\item ou bien pour un $k$ on a $I\subseteq  J_k$,
\item ou 
bien $\exists x\in I$ tel que pour chaque $k$, on a $x\in V_k$. 
%
\end{itemize}
\end{lemma}
\begin{proof}
Il suffit de recopier la démonstration classique donnée dans \cite{Ei}.
\end{proof}

\subsection{Chaînes partiellement spécifiées d'idéaux
premiers} 
\label{subsecProc}

\begin{definition}{\rm  
\label{defproch1} Soit un anneau commutatif $A$. 
\begin{itemize}
\item 
Une \emph{spécification partielle pour une chaîne  d'idéaux 
premiers} (ce que nous abrégeons en {\em \proc}, en anglais {\em 
idealistic chain}) est définie comme suit. Une {\em \prolo $\ell$}  est 
une liste de $\ell+1$ \prois:  $\cC=(\cP_0,\ldots,\cP_\ell)$. \\
On note $\cC(i)$ pour $\cP_i$. La \proc est dite 
{\em finie}  si toutes les parties $J_i$ et $U_i$ sont finies. 
\item  
Une \proc $\cC=(\cP_0,\ldots,\cP_\ell)$ est dite {\em complète} si les 
\prois $\cP_i$ sont complets et si on a les inclusions $\cP_i\subseteq 
\cP_{i+1}$ $(i=0,\ldots ,\ell-1)$. 
\item   
Soient deux \prolos $\ell$, $\cC=(\cP_0,\ldots,\cP_\ell)$ et
$\cC'=(\cP'_0,\ldots,\cP'_\ell)$. \\
On dit que {\em $\cC'$ est un 
raffinement de $\cC$} et on écrit 
$\cC\leq \cC'$ si $\cP_i\leq \cP'_i$ pour $i=0,\ldots ,\ell$. 
\end{itemize}
}
\end{definition}
Nous regardons une \proc $\cC$ de longueur $\ell$  dans $A$ comme une 
spécification partielle pour une chaîne croissante d'idéaux 
premiers
(au sens usuel)  $P_0,\ldots,P_\ell$. Autement dit  
$\cC\leq (\wh{P_0},\ldots, \wh{P_\ell})$. 
 
\begin{lemma} 
\label{factProcComp} 
Toute \proc $\cC=\big((J_0,U_0),\ldots,(J_\ell,U_\ell)\big)$ engendre une \proc 
complète minimale $\cC'=\big((I_0,V_0),\ldots,(I_\ell,V_\ell)\big) $  
définie par les égalités suivantes
\begin{itemize}
\item $ I_0=\gen{J_0}$, $ I_{h+1}=\gen{J_{h+1}\cup I_h}$ pour $0\leq h<\ell$.
\item  $U'_i=\cM(U_i)$ pour $0\leq i<\ell$. 
\item   $V_\ell=U'_\ell+I_\ell$,   $V_{h-1} = U'_{h-1} V_h + I_{h-1}$ pour $1\leq h\leq \ell$. En particulier tout \elt de $V_0$ se réécrit sous la forme
$$u_0\cdot(u_1\cdot(\cdots(u_\ell+j_\ell)+\cdots)+j_1)+j_0=u_0\cdots 
u_\ell+
u_0\cdots u_{\ell-1}\cdot j_\ell+\cdots+u_0\cdot j_1+j_0 
$$
avec les $j_i\in \gen{J_i}$ et $u_i\in\cM(U_i)$. 
\end{itemize}
\end{lemma}

\begin{definition}{\rm  
\label{defproch2}  {\em (Collapsus)} 
Soit  $\cC=\big((J_0,U_0),\ldots,(J_\ell,U_\ell)\big)$ une \proc.
\begin{itemize}
\item  On dit que la \proc $\cC$ {\em collapse} ou encore {\em s'effondre} 
si le dernier élément de la \proc complétée s'effondre, i.e.  s'il existe $j_i\in 
\gen{J_i}$, $u_i\in\cM(U_i)$, $(i=0,\ldots,\ell)$, vérifiant 
l'égalité  
\begin{equation} \label {eqdefproch2}
u_0\cdot(u_1\cdot(\cdots(u_\ell+j_\ell)+\cdots)+j_1)+j_0=0 
\end{equation}
%
\item  Une \proc est dite {\em saturée} si elle est complète
et si les \prois $(J_i,U_i)$ sont saturés. 
\item  La  \proc 
$\big((A,A),\ldots,(A,A)\big)$ est dite {\em triviale}: une chaîne saturée 
qui s'effondre est triviale. 
\end{itemize}
}
\end{definition}

\begin{lemma} \label{lemfinicol2}
Une \proc qui s'effondre raffine un \proc finie qui s'effondre. 
\end{lemma}

\begin{lemma} 
\label{lemma} 
\label{factColAc} 
Une \proc $\cC=(\cP_0,\ldots,\cP_\ell)$ s'effondre  si l'un des $\cP_h$  s'effondre. Plus généralement si une \proc $\cC'$ extraite d'une \proc $\cC$ s'effondre, alors $\cC$ s'effondre. De même si  $\cC$ s'effondre, tout raffinement de $\cC$ s'effondre.  
\end{lemma}
Les démonstrations des propriétés suivantes sont instructives.
\begin{lemma} 
\label{factSatur} 
Soit deux  \procs $\cC^1=(\cP_0,\ldots,\cP_\ell)$ 
et $\cC^2=(\cP_{\ell+1},\ldots,\cP_{\ell+r})$
d'un  anneau $A$ (avec $\cP_i=(J_i,U_i)$ pour chaque $i$) et notons $\cC=\cC^1\bullet 
\cC^2=(\cP_0,\ldots,\cP_{\ell+r})$.
\begin{enumerate}
\item  Supposons $\cC^1$ saturée. Alors $\cC$ s'effondre dans $A$
\ssi $\cP_\ell\bullet \cC^2$   s'effondre dans~$A$
\ssi $\cC^2$   s'effondre dans le quotient $A/J_\ell$. 
\item  Supposons $\cC^2$ complète. Alors $\cC$ s'effondre dans $A$
\ssi $ \cC^1\bullet\cP_{\ell+1}$   s'effondre dans $A$
\ssi $\cC^2$   s'effondre dans le localisé $A_{U_{\ell+1}}$.
\item  Supposons $\cC^1$ saturée et $\cC^2$ complète. 
Alors $\cC$ s'effondre dans $A$  \ssi 
$(\cP_\ell,\cP_{\ell+1})$   s'effondre dans $A$
\ssi $J_\ell\cap U_{\ell+1}\not= \emptyset$.
\end{enumerate}
\end{lemma}
\begin{proof}
Laissée à la lectrice.
\end{proof}
\subsection{Collapsus simultané et \nst formel pour les \procs} 
\label{subsecColsim}
\begin{notations} 
\label{notaAjou}
{\rm  Dans la suite nous utiliserons les notations suivantes
pour certains raffinements élémentaires.
Si $\cP=(J,U)$, $\cC=(\cP_1,\ldots \cP_n)$, nous notons
\begin{itemize}
\item  $\cP\avec{x\in \cP}$ ou encore $(J,x;U)$ 
pour $(J\cup \left\{x\right\}, U)$ 
\item $\cP\avec{x\notin\cP}$  ou encore $(J;x,U)$
pour $(J,U\cup \left\{x\right\})$ 
%
\item  $(a_1,\ldots ,a_m;v_1,\ldots,v_p)$ pour 
$(\left\{a_1,\ldots ,a_m\right\},\left\{v_1,\ldots,v_p\right\})$ 
\item $\cC\avec{x\in \cP_i}$ pour $(\cP_1,\ldots,\cP_i\avec{x\in \cP_i},\ldots, \cP_n)$  
\item   $\cC\avec{x\notin \cP_i}$ pour $(\cP_1,\ldots,\cP_i\avec{x\notin \cP_i},\ldots, \cP_n)$
\end{itemize}
} 
\end{notations}

\Subsubsection{Collapsus simultané pour les \procs} 

\begin{theorem} 
\label{ThColSimKrA}{\em (Collapsus simultané pour les \procs)}\\
Soit une  \proc $\cC=\big((J_0,U_0),\ldots,(J_\ell,U_\ell)\big)$ dans un anneau commutatif $A$. 
\begin{itemize}
\item [$(1)$] Soit $x\in A$ et $i\in\left\{0,\ldots ,\ell\right\}$.
Supposons que les \procs 
$\cC\avec{x\in \cC(i)} $
et 
$\cC\avec{x\notin \cC(i)} $
s'effondrent toutes les deux, alors $\cC$ s'effondre également.
\item [$(2)$] La \proc $\cC$ engendre une \proc saturée minimum.
Celle-ci s'obtient en ajoutant dans 
$U_i$ (resp. $J_i$) tout élément $x\in A$  tel que la \proc
$\cC\avec{x\in \cC(i)} $
(resp. $\cC\avec{x\notin \cC(i)} $) 
s'effondre. 
\end{itemize}
\end{theorem}
\begin{proof}
Notons $(1)_\ell$ et $(2)_\ell$ les énoncés avec $\ell$ fixé.
Remarquons que $(1)_0$ et $(2)_0$ ont fait l'objet du 
théorème~\ref{ThColSim0}.  Nous allons faire une démonstration par récurrence sur 
$\ell$.
Nous pouvons supposer la \proc~$\cC$ complète (cela n'influe pas sur 
l'existence des collapsus).\\
Le fait que $(1)_\ell\Rightarrow (2)_\ell$ n'offre pas de difficulté, on 
raisonne comme dans le théorème \ref{ThColSim0}. \\
Il reste donc à montrer 
$((1)_{\ell-1}\; {\rm et}\; (2)_{\ell-1})\,\Rightarrow\, (1)_{\ell}$ (pour 
$\ell>0$).\\ 
Si $i<\ell$
on utilise le fait \ref{factSatur} qui nous ramène à des \prolos $i$
dans l'anneau localisé $A_{U_{i+1}}$ et on applique l'\hdr.\\
Si $i=\ell$ on considère la \prolo $\ell-1$ 
$\big((K_0,S_0),\alb\ldots,\alb(K_{\ell-1},S_{\ell-1})\big)$ obtenue en saturant 
$\big((J_0,U_0),\ldots,(J_{\ell-1},U_{\ell-1})\big)$ (on applique $(2)_{\ell-1}$). 
Pour des $j_i\in \gen{J_i}$ et 
$u_i\in\cM(U_i)$ arbitraires ($0\le i\le\ell$),
considérons les affirmations suivantes:
$$ u_0\cdot(u_1\cdot(\cdots(u_{\ell-1}\cdot (u_\ell+j_\ell)+ j_{\ell-1}) 
+\cdots)+j_1)+j_0=0\eqno (\alpha)
$$
$$ (u_\ell+j_\ell)\in K_{\ell-1}\eqno (\beta)
$$
$$ \exists n\in\N\;\; u_0\cdot(u_1\cdot(\cdots(u_{\ell-1}\cdot 
(u_\ell+j_\ell)^n+ j_{\ell-1}) +\cdots)+j_1)+j_0=0\eqno (\gamma)
$$
On a $(\alpha) \Rightarrow (\beta) \Rightarrow (\gamma)$.
On a donc les propriétés équivalentes suivantes: 
\begin{itemize}
\item la \proc~$\cC$ s'effondre dans $A$ (certifié par une égalité de type $(\alpha)$ 
ou de type $(\gamma)$), 
\item  le \proi $(J_\ell,U_\ell)$ s'effondre dans
$A/K_{\ell-1}$ (certifié par une égalité de type $(\beta)$).
\end{itemize}
On est 
donc ramené (sur l'anneau $A/K_{\ell-1}$) au cas $(1)_0$ traité au 
théorème
\ref{ThColSim0}. 
\end{proof}
Remarquons que nous ne nous sommes pas appuyés, dans la fin du 
raisonnement, sur le lemme \ref{factSatur}, qui n'est pas assez fort dans cette situation.
Nous donnons dans le corolaire qui suit des conséquences simples du 
théorème \ref{ThColSimKrA}. Le deuxième point renforce le lemme \ref{factSatur}. 
\begin{corollary} 
\label{corColSimKrA}~
\begin{enumerate}
\item  Une \proc $\cC$  s'effondre \ssi toute \proc saturée qui raffine 
$\cC$ est triviale.
\item  On ne change pas le collapsus d'une \proc $\cC$ si on la remplace par sa saturée $\cC'$ ou par n'importe quelle \proc
 $\cC''$ telle que $\cC\leq \cC''\leq \cC'$.
\item  Soit $x_1,\ldots,x_k\in A$ et $i\in\{0,\dots,\ell\}$. Supposons que les \procs 
$$\big((J_0,U_0),\ldots,(J_i\cup\{(x_h)_{h\in H}\}, 
U_i\cup\{(x_h)_{h\in H'}\}), \ldots,(J_\ell,U_\ell)\big)
$$
s'effondrent pour tout couple de parties 
complémentaires $(H,H')$ de $\{1,\ldots,k\}$, 
alors la chaîne  $\cC$ s'effondre également.
\item  Soient $X,\,Y\subseteq A$, $\cC=(\cP_0,\ldots,\cP_\ell)$ une \proc et   
$\cC'=(\cQ_0,\ldots,\cQ_\ell)$ sa saturée. \\
Si $\cP_k=(J_k,U_k)$ 
et
 $\cQ_k=(I_k,V_k)$  la \proc 
$$(\cP_0,\ldots,(J_k,X;Y,U_k),\ldots,\cP_\ell)$$
s'effondre \ssi le \proi
$(I_k,X;Y,V_k)$ s'effondre.
\item  Soient $x\in A$ et $\cC=(\cP_0,\ldots,\cP_\ell)$,  si chacune des 
\procs 
$$
\cC\avec{x\in\cP_0},\quad  
\cC\avec{x\notin\cP_\ell},\quad  
\cC\avec{x\in\cP_i}\avec{x\notin\cP_{i-1}} \,(1\leq i\leq \ell)
$$
s'effondrent, alors $\cC$ s'effondre.
\end{enumerate}
\end{corollary}
\begin{definition}{\rm  
\label{defproch3} 
 Deux \procs qui engendrent la même \proc saturée sont dites 
{\em équivalentes}.
 Une \proc équivalente à une \proc finie est appelée une {\em 
\proc de type fini}.
}
\end{definition}

Notre idée est que le théorème \ref{ThColSimKrA} révèle un 
contenu calculatoire \gui{caché} dans les démonstrations classiques concernant les 
chaînes croissantes d'idéaux premiers. Cette idée est illustrée
par le théorème suivant qui, en \clama, donne une caractérisation
concrète des \procs qui spécifient incomplètement des chaînes 
croissantes d'idéaux premiers.
 
\begin{theorem} 
\label{th.nstformel} {\em (\nst formel pour les \procs dans un anneau commutatif)}  
Soit $A$ un anneau et 
$\cC=\big((J_0,U_0),\ldots,(J_\ell,U_\ell)\big)$ une \proc dans $A$. Le \tcg implique 
que \propeq
\begin{itemize}
\item [$(a)$] Il existe une chaîne de $\ell+1$ idéaux premiers détachables 
$P_0\subseteq \cdots\subseteq P_\ell$  qui raffine $\cC$, i.e. telle que
 $J_i\subseteq P_i$, $U_i\cap P_i=\emptyset $, $(i=0,\ldots,\ell)$. 
\item [$(b)$] La \proc $\cC$ ne s'effondre pas, i.e. pour tous $j_i\in \gen{J_i }$ et $u_i\in\cM(U_i)$, 
$(i=0,\ldots,\ell)$  
$$u_0\cdot(u_1\cdot(\cdots(u_\ell+j_\ell)+\cdots)+j_1)+j_0\neq 0 
$$
\end{itemize}
\end{theorem}
 
\begin{proof}
Seul $(b)\Rightarrow (a)$ pose problème. 
Commençons par une démonstration qui s'appuie, non sur le \tcg, mais sur le 
\pte et le lemme de Zorn. On considère une \proc 
$\cC^1=\big((P_0,S_0),\ldots,(P_\ell,S_\ell)\big)$ 
maximale (pour la relation de raffinement) parmi les \procs qui raffinent 
$\cC$ et qui ne s'effondrent pas. Tout d'abord, il est clair que $\cC^1$ est 
complète, puisqu'on ne change pas le collapsus en la complétant.
Si ce n'était pas une chaîne d'idéaux premiers avec leurs 
compléments, on aurait \hbox{$S_i\cup P_i\neq A$} pour un indice~$i$. 
Dans ce cas, soit $x\in A\setminus (S_i\cup P_i)$. 
Alors  $\big((P_0,S_0),\ldots,(P_i,x;S_i),\ldots,(P_\ell,S_\ell)\big)$
doit s'effondrer (par maximalité). 
La même chose pour 
$\big((P_0,S_0),\ldots,(P_i;x,S_i),\ldots,(P_\ell,S_\ell)\big)$.
Le théorème \ref{ThColSimKrA} nous dit donc que la \proc 
$\cC^1$ s'effondre, ce qui est absurde.\\
Voyons maintenant une démonstration dans laquelle intervient seulement le \tcg, plus élémentaire que les outils utilisés dans la première démonstration.\\
Nous considérons la théorie propositionnelle qui décrit une chaîne croissante de \prois de longueur $\ell$ dans $A$.
Dans cette théorie nous avons pour chaque $i\in\so{0,\dots,\ell}$ un prédicat 
\hbox{pour $x\in P_i$} et un prédicat pour $x\in S_i$.
Les axiomes sont comme les axiomes \emph{1} à~\emph{7} décrits dans la démonstration du corolaire
\ref{corNstformHilbert}. Il faut ajouter les trois inclusions 
$$
(x\in P_i,\; u \in S_i)\Rightarrow x+u\in S_{i},\quad   x\in P_i\Rightarrow x\in P_{i+1}\quad \hbox{  et  }\quad u\in S_i\Rightarrow u\in S_{i-1}.
$$ 
Le seul collapsus est $0\in S_0\Rightarrow\bot$. La théorie obtenue collapse donc exactement quand la chaîne~$\cC$ s'effondre. L'hypothèse $(b)$ est donc que la théorie ne collapse pas.
D'après le théorème \ref{ThColSimKrA} de collapsus simultané, la théorie dans laquelle on ajoute les axiomes disjonctifs $x\in P_i \vee x\in S_i$ est consistante.
Par le \tcg, elle a un modèle. Ce modèle nous donne la chaîne stricte de $\ell+1$
\ideps voulue.
\end{proof}

On notera également que le théorème \ref{th.nstformel}  implique (en 
deux lignes) le théorème de collapsus simultané \ref{ThColSimKrA}. 
Ce dernier peut donc être légitiment considéré comme la version 
\cov du premier.
Un corolaire du théorème \ref{th.nstformel}  serait une 
caractérisation (en \clama) de la \proc saturée engendrée par une \proc $\cC$
(décrite au point 2 du théorème \ref{ThColSimKrA}) au moyen de la famille des chaînes d'idéaux premiers qui sont des  
raffinements de $\cC$, comme dans le corolaire~\ref{corNstfH2} pour le  saturé d'un \proi donné.
\subsection{Suites singulières et dimension de Krull} 
\label{subsecPsr}

\vspace{1em}
\Subsubsection{Définition \cov de la dimension de Krull d'un anneau commutatif} 
\label{nbpIneq}Dans un cadre \cof, il est parfois préférable de 
considérer une relation d'inégalité $x\neq 0$ qui ne soit pas 
simplement l'impossibilité de $x=0$. Par exemple un nombre réel est 
dit $\neq 0$ lorsqu'il est inversible, \cad clairement non nul. Chaque 
fois que nous mentionnons une relation d'inégalité $x\neq 0$ dans un 
anneau, nous supposons donc toujours implicitement que cette relation a 
été définie au préalable dans l'anneau que nous considérons. 
Nous demandons que cette relation soit une inégalité standard, \cad 
qu'elle puisse être démontrée équivalente à $\lnot(x=0)$ en 
utilisant le principe du tiers exclu. Nous demandons en outre que l'on ait 
\cot  
$\; (x\neq 0,\; y=0)\; \Rightarrow\;  x+y\neq 0$,  
$\; xy\neq 0\Rightarrow \; x\neq 0,$ et
$\lnot(0\neq 0)$. Enfin  $x\neq y$ est défini par  $x-y\neq 0$.
En l'absence de précisions concernant  $x\neq 0$, on peut toujours 
considérer qu'il s'agit de la relation $\lnot(x=0)$. Lorsque l'anneau 
est un ensemble discret, \cad lorsqu'il possède un test d'égalité 
à zéro, on choisit toujours l'inégalité $\lnot(x=0)$. Néanmoins 
ce serait une erreur de principe grave de considérer que l'algèbre 
commutative ne doit travailler qu'avec des ensembles discrets.
\begin{definition}{\rm  
\label{defproch} 
Soit  $(x_1,\ldots,x_\ell)$ une suite de longueur $\ell$ dans un anneau 
commutatif $A$. 
\begin{itemize}
\item  La \proc $\big((0;x_1),(x_1;x_2),\alb\ldots, 
(x_{\ell-1};x_\ell),(x_\ell;1)\big)$ est appelée une {\em \proel}. 
On dit qu'elle est {\em associée à la suite $(x_1,\ldots,x_\ell)$}.
On la note $\overline{(x_1,\ldots,x_\ell)}$.
\item  On dit que la suite 
$(x_1,\ldots,x_\ell)$ est une {\em suite singu\-lière}\footnote{Dans la version originale de l'article, on disait suite pseudo  singu\-lière.} 
lorsque
la \proel associée $\overline{(x_1,\ldots,x_\ell)}$ s'effondre. 
Précisément, cela signifie qu'il existe $a_1,\ldots,a_\ell\in A$ et 
$m_1,\ldots,m_\ell\in \N$ tels que
\begin{equation} \label {eqdefproch4}
 x_1^{m_1}(x_2^{m_2}\cdots(x_\ell^{m_\ell} (1+a_\ell x_\ell) + 
\cdots+a_2x_2) + a_1x_1) =  0
\end{equation}
ou encore qu'il existe $a_1,\ldots,a_\ell\in A$ et $m\in \N$ tels que
$$ (x_1x_2\cdots x_\ell)^{m} +
a_\ell (x_1\cdots x_{\ell-1})^{m}x_\ell^{m+1} + 
a_{\ell-1} (x_1\cdots x_{\ell-2})^{m}x_{\ell-1}^{m+1} + 
\cdots+
a_{1} x_1^{m+1} =  0
$$
\item   On dit que la suite 
$(x_1,\ldots,x_\ell)$ est une {\em suite pseudo ré\-gu\-liè\-re} 
lorsque
la \proel associée ne s'effondre pas. \\
Précisément,
pour tous $a_1,\ldots,a_\ell\in A$ et tous $m_1,\ldots,m_\ell\in \N$ on 
a  
$$ x_1^{m_1}(x_2^{m_2}\cdots(x_\ell^{m_\ell} (1+a_\ell x_\ell) + 
\cdots+a_2x_2) + a_1x_1) \neq  0
$$
\end{itemize}
}
\end{definition}
On notera que la longueur de la \proel associée à une suite est égal 
au nombre d'éléments de la suite.

Le lien avec les suites régulières\footnote{La définition le meilleure 
nous semble être celle de Bourbaki: $(x_1,\dots,x_n)$ est dite régulière si chaque $x_i$ est régulier modulo l'idéal $\gen{x_j;j<i}$. Cette définition évite l'usage de la négation. On ne demande pas que $\gen{x_1,\dots,x_n}\neq \gen{1}$.} est donné par la proposition 
suivante, qui est immédiate. 

\begin{proposition} 
\label{prop.regseq} Dans un anneau commutatif, une suite 
régulière qui n'engendre pas l'idéal $\gen{1}$ est pseudo régulière.
\end{proposition}

Le lemme suivant est parfois utile.

\begin{lemma} 
\label{lem.pseudoreg} 
Soient   $(x_1,\ldots,x_\ell)$ et  $(y_1,\ldots,y_\ell)$ dans un anneau 
commutatif $A$. 
Supposons que pour chaque $j$, $x_j$ divise une puissance de $y_j$ 
et $y_j$ divise une puissance de $x_j$. 
Alors la suite  $(x_1,\ldots,x_\ell)$ est singulière \ssi
la suite  $(y_1,\ldots,y_\ell)$ est singulière
\end{lemma}
\begin{proof}
En effet, si   $x$ divise une puissance de $y$  et  $y$
divise une puissance de $x$, 
on établit les relations de raffinement suivantes:
$$(a;x)(x;b)  \leq  (a;x,y)(x,y;b) \leq  
         {\rm \;la\; satur\aigu ee\; de\; la\; suite\; }   (a;x)(x;b)$$
on ajoute le premier  $y$ en disant que  $yc=x^k$
  ($y$ est donc dans le saturé du \mo engendré par~$x$),  
on ajoute le second en disant que      $y^m=dx$
  ($y$  est donc dans le radical de l'idéal engendré par~$x$).
On en déduit par symétrie que $(a;x)(x;b)$ et $(a;y)(y;b)$ ont la 
même saturation. 
\end{proof}

Un corolaire immédiat du théorème \ref{th.nstformel} est le 
théorème suivant \ref{th.pseudoreg}.

\begin{theorem} 
\label{th.pseudoreg} 
{\em (suites pseudo régulières et chaînes croissantes d'idéaux 
premiers)}\\ \Tcgi
Dans un anneau $A$ une suite $(x_1,\ldots,x_\ell)$ est pseudo 
régulière \ssi il existe  $\ell+1$ idéaux premiers 
$P_0\subseteq \cdots\subseteq P_\ell$ 
avec $x_1\in P_1\setminus P_0$, $x_2\in P_2\setminus P_1, \ldots$  
$x_\ell\in P_\ell\setminus P_{\ell-1}$.
\end{theorem}

Ceci conduit à la définition suivante, qui assure un contenu \cof 
explicite à la notion de {\em dimension de Krull d'un anneau}.
\begin{definition}{\rm  
\label{def.dimKrull}~{\em (Dimension de Krull d'un anneau)} 
\begin{itemize}
\item  Un anneau $A$ est dit {\em de dimension $-1$} lorsque $1=_A0$. 
\item  Pour $\ell\geq 1$ l'anneau est dit {\em de dimension $\leq \ell-1$} si toute \proelo  $\ell$ s'effondre. On notera $\dim A\leq l-1$
\end{itemize}
Les définitions suivantes utilisent des négations plus ou moins compliquées.
\begin{itemize}
\item  Il 
est dit de dimension $\geq 0$ si  $1\neq_A0$, 
\item  Il 
est dit de dimension $>-1$ si  $\lnot(1=_A0)$
\item  Il 
est dit de dimension  
$<0$ si  $\lnot(1\neq_A0)$.  
\item  L'anneau est dit {\em de dimension $\geq \ell$} s'il existe une 
suite pseudo régulière de longueur~$\ell$. 
\item  L'anneau est dit {\em de dimension $\ell$} s'il est à la fois de 
dimension $\ge\ell$ et $\le\ell$.  
\item  Il est dit {\em  de dimension $<\ell$ } lorsqu'il est impossible 
qu'il soit de dimension 
$\geq \ell$ 
\item  Il est dit {\em  de dimension $>\ell$ } lorsqu'il est impossible 
qu'il soit de dimension 
$\leq \ell$\footnote{En fait, il existe une et une seule \proel de longueur 
$0$: $(0,1)$, donc il n'y avait pas besoin de commencer
par une définition particulière pour les anneaux de dimension $-1$. 
Dans ce cadre, on retrouve les distinctions entre dimension $\ge 0$ et 
dimension $>-1$, ainsi que  entre dimension $\le -1$ et dimension 
$<0$.}.
\end{itemize}
}
\end{definition}


Un anneau est donc de dimension (de Krull) $\leq \ell-1$ si pour toute 
suite  $(x_1,\ldots,x_\ell)$ dans~$A$, on peut trouver   
$a_1,\ldots,a_\ell\in A$ et  $m_1,\ldots,m_\ell\in \N$ tels que
$$ x_1^{m_1}(\cdots(x_\ell^{m_\ell}(1+a_\ell x_\ell)+\cdots)+a_1x_1)=  0
$$

En particulier un anneau est de dimension $\leq 0$  \ssi pour tout $x\in 
A$ il existe $n\in\N$ et $a\in A$ tels que $x^n=ax^{n+1}$. Tandis qu'il 
est de dimension $<1$ \ssi il est absurde qu'on puisse trouver $x\in A$ 
avec, pour tout  $n\in\N$ et tout $a\in A$, $x^n\neq ax^{n+1}$.

Notez que $\R$ est un anneau local de dimension $<1$, mais on ne
peut pas prouver \cot qu'il est de dimension $\leq 0$.

Notez aussi qu'un anneau est local zéro-dimensionnel \ssi on a
$$\forall x\in A\; \; \; x \; \; {\rm est\; inversible\; ou\;
nilpotent}
$$

\Subsubsection{Dimension de Krull d'un anneau de \pols sur un corps
discret} 
Nous avons tout d'abord.
\begin{proposition} 
\label{propKrDimetDegTr} 
Soit $K$ un corps discret, $A$ une $K$-algèbre commutative, et $x_1$, 
\ldots, $x_\ell$ dans $A$  algébriquement dépendants sur $K$. Alors la 
suite
$(x_1,\ldots,x_\ell)$ est singulière.
\end{proposition}
\begin{proof}
Soit $Q(x_1,\ldots,x_\ell)=0$ une relation de dépendance algébrique
sur $K$. Ordonnons les monômes de $Q$ dont le coefficient est non nul 
selon l'ordre lexicographique. On peut supposer \spdg que le coefficient 
du premier monôme non nul (pour cet ordre) est égal à $1$. Si 
$x_1^{m_1}x_2^{m_2}\cdots x_\ell^{m_\ell}$ est ce monôme, il est clair que 
$Q$ s'écrit sous forme
$$ Q=x_1^{m_1}\cdots x_\ell^{m_\ell}+ 
x_1^{m_1}\cdots x_\ell^{1+m_\ell}R_\ell+
x_1^{m_1}\cdots x_{\ell-1}^{1+m_{\ell-1}}R_{\ell-1}+\cdots+
x_1^{m_1}x_2^{1+m_2}R_2+ x_1^{1+m_1}R_1
$$
ce qui est le collapsus recherché.
\end{proof}

On en déduit.
\begin{theorem} 
\label{thKDP} Soit $K$ un corps discret. La dimension de Krull de l'anneau 
$K[X_1,\ldots,X_\ell]$ est égale à $\ell$.
\end{theorem}
\begin{proof}
Vue la proposition \ref{propKrDimetDegTr} il suffit de vérifier que la 
suite $(X_1,\ldots,X_\ell)$ est pseudo-régulière. Or elle est 
régulière. 
\end{proof}

On notera que nous avons obtenu ce résultat de base avec un minimum 
d'effort, et que notre démonstration vaut évidemment en mathématiques 
classiques (avec la définition classique de la dimension de Krull dès 
qu'on admet le \tcg). Ceci infirme l'opi\-nion couramment admise que les 
\prcos sont nécessairement plus compliquées que les démonstrations 
classiques.
\subsection{Dimension de Krull et principe local-global} 
\label{subsecKrLocGlob}
\Subsubsection{Monoïdes \com} 
\begin{definition}  
\label{def.moco}{\rm ~
\begin{itemize}
\item [$(1)$] Des \mos $S_1,\ldots ,S_n$ de l'anneau $A$ sont dits
 {\em \com} si un idéal de $A$ qui coupe chacun des $S_i$ contient
 toujours $1$, autrement dit si on a:
$$ \forall s_1\in S_1 \;\cdots\; \forall s_n\in S_n \;\; 
\exists a_1,\ldots, a_n\in A\quad \som_{i=1}^{n} a_i s_i =1
$$
\item [$(2)$] On dit que 
{\em les \mos $S_1,\ldots ,S_n$ de l'anneau $A$ recouvrent le \mo $S$} 
si $S$ est contenu dans les $S_i$ et si un idéal de $A$ qui coupe 
chacun des $S_i$  coupe toujours $S$, autrement dit si on a:
$$ \forall s_1\in S_1 \;\cdots\; \forall s_n\in S_n \;\; 
\exists a_1,\ldots, a_n\in A\quad \som_{i=1}^{n} a_i s_i \in S
$$
\end{itemize}
}
\end{definition}

L'exemple fondamental de \moco est le suivant:
 lorsque $s_1,\ldots, s_n\in A$ vérifient 
$\gen{s_1,\ldots, s_n}=\gen{1}$,  les \mos $\cM(s_i)$ sont \com.

Les deux lemmes suivants sont aussi très utiles pour construire des 
\moco.
\begin{lemma} 
\label{lemAssoc} (Recouvrements, calculs immédiats)
\begin{itemize}
\item [$(1)$] (associativité) Si les \mos $S_1,\ldots ,S_n$ de l'anneau 
$A$ recouvrent le \mo $S$  et si chaque $S_\ell$ est recouvert par des 
\mos  
$S_{\ell,1},\ldots ,S_{\ell,m_\ell}$, alors les $S_{\ell,j}$ recouvrent~$S$.
\item [$(2)$] (transitivité) Soit  $S$ un \mo de  l'anneau $A$ et 
$S_1,\ldots ,S_n$ des \moco de  l'anneau localisé $A_S$. 
Pour $\ell=1,\ldots,n$ soit $V_\ell$ le \mo
de  $A$ formé par les numérateurs des éléments de  $S_\ell$. 
Alors les \mos $V_1,\ldots ,V_n$ recouvrent  $S$. 
\end{itemize}
\end{lemma}

Rappelons que  pour  $I,U\subseteq A$, nous notons $\cS(I,U)$ le \mo 
\hbox{$\cM(U)+\gen{I}$}  (obtenu en complétant le \proi $(I,U)$).

\begin{lemma} 
\label{lemRecouvre} 
Les \mos
$\cS(I,a;U)$ et $\cS(I;a,U)$ recouvrent le \mo $\cS(I,U)$. 
\end{lemma}
\begin{proof}
Pour  $x\in \cS(I;U,a)$ et  $y\in \cS(I,a;U)$ 
on doit trouver un $x_1x+y_1y\in \cS(I;U)$ ($x_1,y_1\in A$). 
On écrit $x=u_1a^k+j_1$, $y=u_2+j_2-az$ avec $u_1,u_2\in \cM(U)$,
   $j_1,j_2\in \cI(I)$,   $z\in A$. L'identité fondamentale 
$\, c^k-d^k=(c-d)\times \cdots\, $ avec $c=u_2+j_2$ et $d=az$ donne un 
$y_2\in A$  tel que $$y_2y=(u_2+j_2)^k-(az)^k=(u_3+j_3)-(az)^k$$
et on écrit $z^kx+u_1y_2y=u_1u_3+u_1j_3+j_1z^k=u_4+j_4$.
\end{proof}

\begin{corollary} 
\label{corS} Soient $u_1,\ldots,u_n\in A$. 
Notons  
$$
S_1=\cS(0;u_1) ,\, S_2=\cS(u_1;u_2) ,\, \dots, S_n=\cS(u_1,\dots,u_{n-1};u_n)\hbox{ et }
S_{n+1}=\cS(u_1,\dots,u_n;1).
$$ 
Alors les \mos $S_1,\ldots,S_{n+1}$ sont \com.
\end{corollary}

Les \moco constituent un outil constructif qui permet en général de 
remplacer des arguments abstraits de type local-global par des calculs 
explicites.
Si  $S_1,\ldots ,S_n$ sont des \moco de l'anneau $A$, l'anneau produit des 
localisés $A_{S_i}$ est une $A$-algèbre fidèlement plate.
Donc de très nombreuses propriétés sont vraies avec $A$ \ssi elles 
le sont avec chacun des $A_{S_i}$.

Nous le vérifions pour la dimension de Krull dans le paragraphe suivant.
\Subsubsection{Caractère local de la dimension de Krull} 
La proposition suivante est immédiate.
\begin{proposition} 
\label{propKrLoc} 
Soit $A$ un anneau. 
\begin{enumerate}
\item La dimension de Krull de $A$ est toujours supérieure ou 
égale à celle d'un quotient ou d'un localisé de $A$. 
\item Plus précisément toute \proel qui s'effondre dans $A$ s'effondre 
dans tout quotient et dans tout localisé de $A$, et toute \proel  
dans un localisé de $A$ est équivalente 
(dans le localisé) à une \proel 
écrite dans $A$.
\item Enfin, si une \proel $\ov{(a_1,\ldots,a_\ell)}$ de $A$ s'effondre dans un 
localisé $A_S$, il existe un $m$ dans $S$ tel que  
$\ov{(a_1,\ldots,a_\ell)}$ 
s'effondre dans $A[1/m]$.
\end{enumerate}
\end{proposition}

\begin{proposition} 
\label{propKrLocGlob} 
Soient  $S_1,\ldots ,S_n$  des \moco  de l'anneau $A$ et $\cC$ une \proc 
dans $A$. Alors $\cC$ s'effondre dans $A$ \ssi elle  s'effondre dans chacun 
des
 $A_{S_i}$. En particulier la dimension de Krull de $A$ est $\leq \ell$ 
\ssi
la dimension de Krull de chacun des  $A_{S_i}$ est $\leq \ell$.
\end{proposition}
\begin{proof}
Il faut montrer qu'une \proc $\cC$ s'effondre dans $A$ si elle  s'effondre 
dans chacun des $A_{S_i}$. Pour simplifier prenons une chaîne de 
longueur $2$: $\big((J_0,U_0),(J_1,U_1),(J_2,U_2)\big)$ avec des idéaux $J_k$  
et des \mos $U_k$. Dans chaque  $A_{S_i}$ on a une égalité
$$ u_{0,i}\,u_{1,i}\,u_{2,i}+u_{0,i}\,u_{1,i}\,j_{2,i} 
+u_{0,i}\,j_{1,i}+j_{0,i}=0
$$
avec $u_{k,i}\in U_k$ et  $j_{k,i}\in J_kA_{S_i}$.
Après avoir chassé les dénominateurs et multiplié par un 
élément convenable de $S_i$ on obtient une égalité dans $A$ du 
type suivant
$$ s_i\,u_{0,i}\,u_{1,i}\,u_{2,i}+u_{0,i}\,u_{1,i}\,j'_{2,i} 
+u_{0,i}\,j'_{1,i}+j'_{0,i}=0
$$
avec $s_i\in S_i$, $u_{k,i}\in U_k$ et  $j'_{k,i}\in J_k$.
On pose $u_k=\prod_iu_{k,i}$. En multipliant l'égalité précédente 
par un produit convenable on obtient une égalité
$$ s_i\,u_{0}\,u_{1}\,u_{2}+u_{0}\,u_{1}\,j''_{2,i} 
+u_{0}\,j''_{1,i}+j''_{0,i}=0
$$
avec $s_i\in S_i$, $u_{k}\in U_k$ et  $j''_{k,i}\in J_k$.
Il reste à écrire $\sum_i a_is_i=1$, à multiplier l'égalité
correspondant à $S_i$ par $a_i$ et à faire la somme.
\end{proof}

\Subsubsection{Exemple d'application} 
En \clama la dimension de Krull d'un anneau est la borne supérieure
des dimensions de Krull des localisés en tous les idéaux maximaux.
Cela résulte facilement (en \clama) des propositions \ref{propKrLoc} et 
\ref{propKrLocGlob}.

La proposition \ref{propKrLocGlob}  devrait permettre d'obtenir \cot les 
mêmes conséquences concrètes (que celles qu'on obtient non \cot en 
\clama en appliquant la propriété ci-dessus) même lorsqu'on n'a pas 
accès aux \ideps.  

\ss Nous nous contentons ici d'un exemple simple, dans lequel on a accès 
aux \ideps.
Supposons que nous ayions une \prco simple
que la dimension de Krull de $\Z_{(p)}[x_1,\ldots,x_\ell]$ est 
 $\leq\ell+1$ 
($p$ est un nombre premier arbitraire, et $\Z_{(p)}$ est le localisé de 
$\Z$ en $p\Z$). Alors nous pouvons en déduire le même résultat pour 
 $A=\Z[x_1,\ldots,x_\ell]$ en appliquant le principe local-global 
précédent.

Soit en effet une liste $(a_1,\ldots,a_{\ell+2})$ dans $A$. Le collapsus 
de la \proel $\ov{(a_1,\ldots,a_{\ell+2})}$ dans
$\Z_{(2)}[x_1,\ldots,x_\ell]$ se relit comme un collapsus  dans
$\Z[1/m_0][x_1,\ldots,x_\ell]$ pour un certain $m_0$ impair.
Pour chacun des diviseurs premiers $p_i$ de $m_0$ ($i=1,\ldots,k$), le 
collapsus de la \proel $\ov{(a_1,\ldots,a_{\ell+2})}$ dans
$\Z_{(p_i)}$ se relit comme un collapsus  dans
$\Z[1/m_i][x_1,\ldots,x_\ell]$ pour un certain $m_i$ étranger à $p_i$.
Les entiers $m_i$ ($i=0,\ldots,k$) engendrent l'idéal $\gen{1}$, donc 
les \mos $\cM(m_i)$ sont \com et on peut appliquer la proposition
\ref{propKrLocGlob}.

\section{Treillis  distributifs, \entrels et dimension de Krull}
\label{secKrullTreil}
\subsection{Treillis  distributifs, idéaux et filtres}
\label{subsecTrd1}
Un \trdi est un ensemble ordonné avec sup et inf finis, un élément
minimum (noté~$0$) et un élément maximum
(noté $1$). On demande que les lois sup et inf soient distributives
l'une par rapport à l'autre (une distributivité implique l'autre). On
note  ces lois $\vu$ et $\vi$.
La relation $a\leq b$  peut être définie par $a\vu b=b$. La
théorie des \trdis est alors purement équationnelle.
Il y a donc des \trdis définis par générateurs et relations.

Une règle particulièrement importante, dite {\em coupure}, est la
suivante
$$ \left( ((x\vi a)\; \leq\;  b)\quad\&\quad  (a\; \leq\; (x\vu  b))
\right)\; \Longrightarrow \; a \leq\;  b
$$
pour la démontrer on écrit $ x\vi a\vi b=x\vi a$  et
$a= a\vi(x\vu b)$ donc
$$ a=(a\vi x)\vu(a\vi b)=(a\vi x\vi b)\vu(a\vi b)=a\vi b
$$

Un ensemble totalement ordonné est un \trdi s'il possède un maximum
et un minimum. On note ${\bf n}$ un ensemble totalement ordonné
à $n$ éléments (c'est un \trdi si $n> 0$).
Tout produit de \trdis est un \trdi.
Les entiers naturels, munis de la relation de divisibilité,
forment un \trdi (le minimum absolu est $1$  et le maximum absolu est
$0$).
Si $T$  et $T'$ sont deux \trdis, l'ensemble $\Hom(T,T')$ des \homos
(applications conservant, sup, inf, 0 et 1) de  $T$ vers $T'$ est muni
d'une structure d'ordre naturelle donnée par
$$ \varphi \leq \psi \equidef  \forall x\in T\;
\;
\varphi(x) \leq \psi(x)
$$
Une application entre deux treillis totalement ordonnés $T$ et $S$
est un \homo \ssi elle est croissante et si $0_T$  et $1_T$  ont pour images
$0_S$ et $1_S$.

\begin{notations}
{\rm  On note $\Pfe(X)$ l'ensemble des parties finiment énuméréess de l'ensemble~$X$.
Si $A$  est une partie finiment énumérée d'un \trdi $T$ on notera
$$ \Vu A:=\Vu\nolimits_{\!x\in A}x\qquad {\rm et}\qquad \Vi A:=\Vi\nolimits_{x\in A}x
$$
Pour un \trdi $T$, on note $A \vda B$ ou $A \vdash_T B$ la relation définie comme suit sur~$\Pfe(T)$.
$$ A \vda B \; \; \equidef\; \; \Vi A\;\leq \;
\Vu B
$$
}
\end{notations}

Dans la suite nous dirons \gui{partie finie} au lieu de \gui{partie finiment énumérée}.
 
Notez que la relation  $A \vda B$ est bien définie sur l'ensemble des
parties finies parce que les lois~$\vi$  et~$\vu$  sont associatives,
commutatives et idempotentes.

Notons que $\; \emptyset  \vda \{x\}$ signifie $x=1$, et
$ \{y\} \vda \emptyset$ signifie $y=0$.

La relation ainsi définie sur les parties finies vérifie les axiomes suivants, dans lesquels on
écrit~$x$ pour $\{x\}$ et $A,B$  pour $ A\cup B$.
$$\begin{array}{rcrclll}
&    & a  &\vda& a    &\; &(R)     \\[.3em]
A \vda B &   \; \Longrightarrow \;  & A,A' &\vda& B,B'   &\; &(M)     \\[.3em]
(A,x \vda B)\;
\&
\;(A \vda B,x)  &   \Longrightarrow  & A &\vda& B &\;
&(T)
\end{array}$$
on dit que la relation est \emph{réflexive}, \label{remotr} \emph{monotone} et
\emph{transitive}. La troisième règle (transitivité) s'appelle aussi la
règle de \emph{coupure} (elle généralise la règle de coupure décrite au début de la section).

Signalons aussi les deux règles suivantes de \emph{distributivité}:
$$\begin{array}{rcl}
(A,\;x \vda B)\;\& \;(A,\;y \vda B)  &  \;  \Longleftrightarrow  \; &
A,\;x\vu y \vda B  \\[.3em]
(A\vda B,\;x )\;\&\;(A \vda B,\;y)  &   \Longleftrightarrow  &
A\vda B,\;x\vi y
\end{array}$$

La proposition suivante est facile.
\begin{propdef}
\label{propIdealFiltre} \emph{(Quotients particuliers d'un \trdi)}
Soit $T$ un \trdi et
$(J,U)$ un couple de parties de $T$.
Nous noterons $T/(J=0,U=1)$ \hbox{ou $T/(J,U)$} le treillis quotient
 de $T$ par les
relations $x=0$ pour les $x\in J$ et $y=1$ pour \hbox{les $y\in U$}. Ce quotient est décrit précisément comme suit.
\begin{enumerate}
\item  La classe d'équivalence de $0$ est l'ensemble des $a$ qui
vérifient
$$ \exists J_0\in\Pfe(J), U_0\in\Pfe(U) \qquad
a,\; U_0  \vdi T    J_0.
$$
\item  La classe d'équivalence de $1$  est l'ensemble des $b$ qui
vérifient:
$$  \exists J_0\in\Pfe(J), U_0\in\Pfe(U)\qquad
 U_0  \vdi T   b,\, J_0.
$$
\item  Plus généralement, on a $a\leq_{T/(J,U)}b$ \ssi
il existe une partie finie $J_0$ de $J$ et une partie finie $U_0$ de
$U$ telles que
$$  a,\; U_0 \vdi T  b,\; J_0.
$$
Dans le cas où $J$ et $U$ sont finies, en posant $j=\Vi J$ et $u=\Vu U$ on obtient simplement la condition 
$$  a,\, u \vdi T  b,\,j .
$$
\end{enumerate}
\end{propdef}

Dans cette proposition, si $U=\emptyset$, la classe d'équivalence de $0$ est appelée un {\em idéal} du treillis, c'est l'idéal engendré par $J$. On le note $\gen{J}_T$. Un idéal $I$ est soumis aux seules contraintes suivantes\footnote{Si on remplace $T$ par un anneau commutatif $A$, $\vu$ par l'addition et $\vi$ par la multiplication, on trouve la notion d'idéal de l'anneau $A$.}:
$$\begin{array}{rcl}
  & &  0 \in I   \\[.3em]
x,y\in I& \Longrightarrow   &  x\vu y \in I   \\[.3em]
x\in I,\; z\in T& \Longrightarrow   &  x\vi z \in I   
\end{array}$$
(la dernière se réécrit $(x\in I,\;y\leq x)\Rightarrow y\in I$).

En outre pour tout homomorphisme de \trdis $\varphi :T\rightarrow T'$,
$\varphi^{-1}(0)$ est un idéal de~$T$.

Un {\em idéal principal} est un idéal engendré par un élément
$a$. On a  $\gen{a}_T=\{x\in T\; ;\; x\leq a \}$. S'il n'est pas nécessaire de préciser $T$, on le note $\dar a$. Tout \itf est principal.

La notion opposée\footnote{Par renversement de la relation d'ordre.}  de celle d'idéal est la notion de {\em filtre}.
Un filtre $F$ est l'image réciproque de~$1$ par un \homo de \trdis.
 Il est soumis aux seules contraintes
$$\begin{array}{rcl}
  & &  1 \in F   \\[.3em]
x,y\in F& \Longrightarrow   &  x\vi y \in F   \\[.3em]
x\in F,\; z\in T& \Longrightarrow   &  x\vu z \in F   
\end{array}$$

Soit $\psi:T\rightarrow T/ (J,U)$ la projection canonique.
Si $I$ est l'idéal $\psi^{-1}(0)$ et $F$ le filtre~$\psi^{-1}(1)$, 
on dit  que {\em l'idéal $I$  et le filtre $F$ sont conjugués}.
D'après la proposition précédente, un idéal~$I$  et un filtre $F$
sont conjugués \ssi on a pour tous $x\in T,\,I_0\in\Pfe(I)$, \hbox{$F_0\in\Pfe(F)$}:
$$
  (x,\; F_0 \vda I_0) \,\Rightarrow\,  x\in I  \quad{\rm  et}\quad 
(F_0 \vda x,\;  I_0)\,\Rightarrow\,  x\in F
$$
Cela peut se reformuler comme suit:
$$
(f\in F,\; x\vi f \in I) \,\Rightarrow\, x\in I
\quad {\rm et}\quad
(j\in I,\; x\vu j \in F) \,\Rightarrow\, x\in F
$$
Lorsqu'un idéal $I$ et un filtre $F$  sont conjugués, on a
$$
1\in I\; \;\Longleftrightarrow\;\;  0\in F
\;\; \Longleftrightarrow\;\;  (I,F)=(T,T)
$$

D'après les propriétés générales des structures quotients dans les théories purement équationnelles,, tout \homo de \trdis $\varphi:T\to T_1$ qui vérifie les égalités
$\varphi(J)=\{0\}$ et $\varphi(U)=\{1\}$  se factorise de manière unique
par le treillis quotient $T/(J,U)$.

On voit cependant sur l'exemple des ensembles totalement ordonnés
qu'une structure quotient d'un \trdi n'est pas en général
caractérisée par les classes d'équivalence de $0$ et $1$.

\subsection{Treillis  distributifs et \entrels}
\label{subsecTrd2}
Une manière intéressante d'aborder la question des \trdis définis
par générateurs et relations est de considérer la relation
$A \vda B$ définie sur l'ensemble $\Pfe(T)$ des parties
finiment énumérées d'un \trdi $T$.
En effet, si $S\subseteq T$ engendre $T$ comme treillis, alors la
connaissance de la relation
$\vda$ restreinte à $\Pfe(S)$ suffit à caractériser sans ambigüité
le treillis $T$, car toute formule sur $S$ peut être réécrite, au
choix, en forme normale conjonctive (inf de sups dans $S$) ou normale
disjonctive (sup de infs dans $S$). Donc si on veut comparer deux
éléments du treillis engendré par $S$ on écrit le premier en forme
normale disjonctive, le second en forme normale conjonctive, et on 
remarque
que
$$ \Vu\nolimits_{i\in I}\left(\Vi A_i \right)\; \leq \; \Vi\nolimits_{j\in J}\left(\Vu B_j
\right)
\qquad \Longleftrightarrow\qquad  \&_{(i,j)\in I\times J}\;  \left(
A_i \vda  B_j\right)
$$

\begin{definition}{\rm 
\label{defEntrel}
Pour un ensemble $S$ arbitraire, une relation sur  $\Pfe(S)$  qui est
réflexive, monotone et transitive (voir page \pageref{remotr}) est
appelée une {\em \entrel} (en anglais {\em entailment relation}).
}
\end{definition}

L'origine des \entrels est dans le calcul des séquents de Gentzen, qui
mit le premier l'accent sur la règle $(T)$ (la coupure).
Le lien avec les \trdis a été mis en valeur dans \cite{cc,cp,Lor1951}.
Le théorème suivant (cf. \cite{cc,Lor1951}) est fondamental. Il dit que les
trois propriétés des \entrels sont exactement ce qu'il faut pour que
l'interprétation en forme de \trdi soit adéquate.

\begin{theorem}
\label{thEntRel1} {\rm  (Théorème fondamental des \entrels)} On
considère un ensemble~$S$  avec une \entrel
$\vdash_S$ sur $\Pfe(S)$. On considère le \trdi $T$ défini par
générateurs et relations comme suit: les générateurs sont les
éléments de $S$ et les relations sont les
$$ A\; \vdash_T \;  B
$$
chaque fois que $A\; \vdash_S \; B$.  Alors pour toutes parties finies $A$
et $B$ de $S$  on a
$$  A\; \vdash_T \;  B
\; \Longrightarrow \; A\; \vdash_S \;  B
$$
\end{theorem}
\begin{proof}
On donne une description explicite possible du \trdi $T$. Les
éléments de $T$ sont représentés par
des ensembles finis d'ensembles finis d'éléments de $S$
$$X=\{A_1,\dots,A_n\}$$
(intuitivement $X$ représente $\Vi A_1\vu\cdots\vu\Vi A_n$).
On définit alors de manière inductive la relation
 $A\prec Y$ pour $A\in \Pfe(S)$ et $Y\in T$
(intuitivement, $\Vi A\leq \Vu_{C\in Y} \left(\Vi C\right) $) comme suit
\begin{itemize}
\item si $B\in Y$ et $B\subseteq A$ alors $A\prec Y$
\item si on a $A\vdash_S y_1,\dots,y_m$ et $A,y_j\prec Y$ pour
$j=1,\ldots,m$ alors $A\prec Y$
\end{itemize}
On montre facilement que si $A\prec Y$ et $A\subseteq A'$ alors on
a $A'\prec Y.$ On en déduit que $A\prec Z$ si $A\prec Y$ et $B\prec Z$
pour tout $B\in Y$. On peut alors définir $X\leq Y$ par
$A\prec Y$ pour tout $A\in X$ et on vérifie que $T$ est alors
un treillis distributif\footnote{$T$ est en fait le quotient de
$\Pfe(\Pfe(S))$ par la relation d'équivalence: $X\leq Y$ et $Y\leq X$.}
pour les opérations
$$0 = \emptyset~~~~1 = \{\emptyset\}~~~~~X\vee Y = X\cup Y~~~~~
X\wedge
Y = \{ A \cup B~|~A\in X,~B\in Y\}
$$
Pour ceci on montre que si $C\prec X$ et $C\prec Y$ alors on a
$C\prec X\vi Y$ par induction sur les démonstrations de $C\prec X$ et
$C\prec Y$.

On remarque alors que si $A\vdash_S y_1,\dots,y_m$ et $A,y_j\vdash_S B$
pour tout $j$ alors $A\vdash_S B$ en utilisant $m$ fois la règle de
coupure. Il en résulte que si on a $A\vdash_T B$,
c'est à dire $A\prec \{\{b\}~|~b\in B\}$, alors on a $A\vdash_S B$.
\end{proof}

Comme première application, nous pouvons citer la description du
treillis booléen engendré par un \trdi.
Rappelons qu'un treillis booléen
est un \trdi muni d'une loi $x\mapsto \overline{x}$ qui vérifie, pour 
tout
$x$: $x\vi \overline{x}=0$ et $x\vu \overline{x}=1$.
L'application  $x\mapsto \overline{x}$  est alors
un \iso du treillis sur son opposé.
\begin{proposition}
\label{propTrBoo}
Soit $T$ un \trdi ($\neq \Un$). Il existe un treillis booléen
engendré
par $T$. Il peut être décrit comme le \trdi engendré par
l'ensemble
$T_1=T\cup\overline{T}$ (où $\overline{T}$  est une copie de $T$ disjointe de
$T$)
 muni de la \entrel $\vdash_{T_1}$ définie comme
suit:
si $A,B,A',B'$  sont quatre parties finies de $T$  on a
$$
A,\overline{B}\;\vdash_{T_1}\; A',\overline{B'}\equidef A,B'\vda
A',B\quad
{\rm dans} \;  T
$$
Si on note $T_{Bool}$ ce treillis booléen, l'ensemble $T_1$ s'injecte naturellement dans $T_{Bool}$. En particulier le morphisme naturel
de $T$ dans $T_{Bool}$ est injectif, il identifie $T$ à son image 
dans~$T_{Bool}$.
\end{proposition}
\begin{proof}
Voir \cite{cc}.
\end{proof}

Une autre application concerne les constructions de treillis quotients.
Le contexte général est le suivant:
$T$ est un \trdi, $S$ un système générateur de $T$ et $A,B,X,Y\in\Pfe(S)$. On note $\vda$ la \entrel de~$T$ (restreinte à $\Pfe(S)$). 
On considère un treillis quotient $T'$. On note $\,\vdash'\,$ la \entrel de $T'$ restreinte à $\Pfe(S)$.
Notons que comme cas particulier on a les mêmes résultats pour $S=T$.

Le fait que $T'$ est un quotient de $T$ sous certaines relations imposées, et le fait que la structure de $T'$ est entièrement déterminée par la restriction de $\vdi{T'}$ à $S$ implique que la \entrel  $\,\vdash'\,$ est la plus forte des \entrels sur $S$ moins fortes que $\vda$ et pour lesquelles  les conditions imposées sont réalisées.

\begin{proposition} \label{lemquotrdi2}
On considère le treillis  $T'$ quotient de $T$ obtenu en forçant la relation $A\,\vdash'\, B$. 
 \Propeq
\begin{enumerate}
\item $X\,\vdash'\,Y$
\item Pour chaque $a\in A$ et chaque $b\in B$, on a: $X,b\vda Y$ et $X\vda Y,a$.
\end{enumerate}
\end{proposition}
\begin{proof}
On vérifie facilement que la relation $X\,\vdash'\, Y$ définie par le point \emph{2} est une \entrel moins forte que $\vda$ et qu'elle satisfait  $A\,\vdash'\, B$.
Il suffit donc de vérifier que $X\,\vdash'\, Y$ implique $X\vdi1Y$ si $\vdi1$
est une \entrel moins forte que~$\vda$ qui satisfait $A\vdi1 B$.
Or si $X\,\vdash'\, Y$ alors $X,b\vdi1Y$ et $X\vdi1Y,a$ pour chaque $a\in A$ et $b\in B$, et puisque $A\vdi1B$, cela donne après de nombreuses coupures $X\vdi1Y$.
\end{proof}

\subsection{Idéaux premiers, premiers idéaux et spectre}
 Classiquement un {\em idéal premier} $I$ d'un treillis
$T\neq \Un$  est un idéal dont
le complémentaire~$F$ est un filtre (qui est alors un {\em filtre
premier}).
En \clama il revient au même de dire que
$$ 1\notin I\qquad {\rm et}\qquad (x\vi y)\in I\; \Longrightarrow \;
(x\in I \;{\rm ou }\;y\in I)\eqno(*)
$$
ou encore de dire que $I$ est le noyau d'un homorphisme de $T$ vers le
treillis à deux éléments, noté $\Deux$.
En \coma il semble logique de choisir la définition $(*)$ qui est
moins contraignante, mais en donnant sa pleine signification \cov au
\gui{ou}. On définit la notion de filtre premier en renversant l'ordre.

On appelle {\em spectre du \trdi $T$} et on note $\Spec(T)$ l'ensemble
ordonné $\Hom(T,\Deux)$. Il revient au même de
considérer l'ensemble des \ideps détachables de $T$. La
relation d'ordre correspond alors à l'inclusion renversée des
\ideps\footnote{Pour la topologie spectrale usuelle de $\Spec(T)$, la relation d'ordre usuellement choisie est celle de l'inclusion des \ideps.
Cela ne change pas notre propos sur la dimension de Krull.}.
\hbox{On a $\Spec(\Un)= \emptyset$}, $\Spec(\Deux)\simeq \Un$,  $\Spec(\Trois)\simeq \Deux$,
$\Spec(\Quatre)\simeq \Trois$, etc\ldots

\smallskip Le lemme suivant montre la pertinence du point de vue \gui{\entrel sur une partie génératrice} pour décrire un \trdi.
\begin{lemma} \label{lemidepsurGen} (en \clama)
Soit $T$ un \trdi et $S$ une partie génératrice de $T$. Un \idep $P$ de $T$
est caractérisé par sa trace sur $S$. 
\end{lemma}
%
\begin{proof}
Notons $F$ le filtre complémentaire de $P$, $P_S=P\cap S$ et $F_S=F\cap S$.
Un \elt $z$ arbitraire de $T$ s'écrit sous forme $z=\Vu_i(\Vi A_i)$ pour un nombre fini de parties finies $A_i$ de $S$. 
Puisque~$P$ est un idéal, on a $z\in P$ \ssi chacun des $\Vi A_i\in P$.
Il nous faut donc tester si  $\Vi A\in P$ pour une partie finie $A$ de $S$.
Soit $a=\Vi A=a_1\vi \dots \vi a_k$ avec les $a_j\in S$. Si chacun des $a_j\in F_S$ alors $a\in F$ donc $a\notin P$. Si l'un des $a_j$ est dans $P_S$, alors $a\in P$. 
\end{proof}
%

\begin{definition}{\rm 
\label{defProiT} Soit $T$ un \trdi.
\begin{enumerate}
\item 
Une {\em spécification partielle pour un idéal premier} (ce que nous 
abrégeons en {\em \proi}, en anglais {\em idealistic prime}) est un 
couple $\cP=(J,U)$ de parties de $T$. 
\item Un \proi $\cP=(J,U)$ est dit {\em complet} si $J$ est un  idéal et  
$U$ est un filtre. 
\item Si $P$ est  un \idep et $S=T\setminus P$, on note $\wh{P}$ le \proi $(P,S)$.
\item Soient $\cP_1=(J_1,U_1)$ et $\cP_2=(J_2,U_2)$  deux \prois et $P$ un \idep.
\begin{itemize}
\item On  dit que {\em $\cP_1$ est contenu dans $\cP_2$} et on écrit 
$\cP_1\subseteq \cP_2$ si $J_1\subseteq J_2$ et  $U_2\subseteq U_1$.
\item  On  dit que  \emph{$\cP_2$ est un raffinement de $\cP_1$} et on 
écrit 
$\cP_1\leq \cP_2$ si $J_1\subseteq J_2$ et  $U_1\subseteq U_2$.
\item  On  dit que  \emph{l'\idep $P$ raffine le \proi $\cP$} si $\cP\leq \wh{P}$
(ceci explique le point~1). 
\end{itemize}
\item  Un \proi $(J,U)$  est dit {\em saturé} si $J$  est un idéal, 
$U$ un filtre et si $J$  et $U$  sont conjugués. Tout \proi engendre un 
\proi saturé $(I,F)$ décrit à la proposition \ref{propIdealFiltre}.
\item  On dit que le \proi $(J,U)$ {\em collapse} ou
{\em s'effondre} si  le treillis quotient $T'=T/(J,U)$ est réduit à un
point, \cad  $1 \leq_{T'}0$. Cela signifie que l'idéal engendré par $J$ coupe le filtre engendré par $U$. 
\end{enumerate}
}
\end{definition}

On peut résumer le collapsus par les deux lemmes qui suivent.
 
\begin{lemma} \label{lemfinicol3}
Un \proi qui s'effondre raffine un \proi fini qui s'effondre. 
\end{lemma}
\begin{lemma} \label{lemfinicol4}
Un \proi fini  $(J,U)$   s'effondre  \ssi $U\vda J$
dans  $T$. 
\end{lemma}

Le lemme suivant est une généralisation du lemme \ref{lemidepsurGen}. En outre sa preuve est constructive, sans utilisation du tiers ecxlu.
\begin{lemma} \label{lemproisatsurGen}
Soit $T$ un \trdi et $S$ une partie génératrice de $T$. 
\begin{enumerate}
\item Un \proi saturé  $(J,U)$ de $T$ est caractérisé par sa trace sur $S$.
\item Plus précisément, soient $J_S$ et $U_S$ deux parties de $S$ qui satisfont les conditions nécessaires suivantes:
\begin{itemize}
\item si $J_0\in\Pfe(J_S)$, $U_0\in\Pfe(U_S)$, $x\in S$ et $U_0\vda x,J_0$, alors $x\in U_S$,
\item si $J_0\in\Pfe(J_S)$, $U_0\in\Pfe(U_S)$, $x\in S$ et $U_0,x\vda U_0$, alors $x\in J_S$. 
\end{itemize}
Alors il existe un unique \proi saturé $(J,U)$ de $T$ tel que $J_S=J\cap S$
\hbox{et $U_S=U\cap S$}.
\item On peut dire la même chose sous la forme suivante, d'apparence plus générale.\\
 Si $J_S$ et $U_S$ sont deux parties de $S$ et si $(J,U)$ est le \proi saturé de $T$ engendré par le \proi $(J_S,U_S)$, alors 
\begin{itemize}
\item $U\cap S=\sotq{x\in S}{\exists U_0\in\Pfe(U_S)\;  \exists J_0\in\Pfe(J_S)\;\;\;U_0, \vda x,J_0}$,
\item $J\cap S=\sotq{x\in S}{\exists U_0\in\Pfe(U_S)\;  \exists J_0\in\Pfe(J_S)\;\;\;U_0,x \vda J_0}$. 
\end{itemize}
\end{enumerate}
\end{lemma}
%
\begin{proof} On démontre le point \emph{3}.
On considère sur $S$ la \entrel $\vdi1$ définie comme suit
\begin{equation} \label {eqlemproisatsurGen}
A\vdi1 B \equidef  \exists U_0\in\Pfe(U_S)\;  \exists J_0\in\Pfe(J_S)\;\;\; U_0, A\vda B,J_0
\end{equation}
Il est immédiat qu'il s'agit bien d'une \entrel. Le treillis $T_1$ défini par $\vdi1$ est le treillis quotient  $T/(J_S=0,U_S=1)$.
En effet, la \entrel définie par le second membre de l'équivalence $(\ref{eqlemproisatsurGen})$ est la plus petite possible parmi celles qui raffinent $\vda$ et qui forcent $J_S=0$ \hbox{et $U_S=1$}. 
En notant  $J=\sotq{z\in T}{z=_{T_1}0}$ et $U=\sotq{z\in T}{z=_{T_1}1}$ on obtient donc \hbox{que $(J,U)$} est le \proi saturé de $T$ engendré par $(J_S,U_S)$.
Et les égalités indiquées pour $J\cap S$ et $U\cap S$ sont immédiates. 
\end{proof}

On établit maintenant le théorème suivant analogue au théorème \ref{ThColSim0}.

\begin{theorem}
\label{thColSimT1} {\em (Collapsus simultané pour les \prois d'un \trdi)} Soit un
\proi $(J,U)$  et un élement $x$ dans un \trdi $T$.
\begin{enumerate}
\item  Si les \prois $(J,x;U)$ et
$(J;U,x)$  s'effondrent, alors  $(J,U)$  s'effondre également.
\item  Le \proi  $(J,U)$  engendre un \proi saturé minimum.
Celui-ci s'obtient en ajoutant dans
$U$ (resp. $J$) tout élément $x\in A$  tel que le \proi
 $(J,x;U)$
(resp. $(J;U,x)$) s'effondre.
\end{enumerate}
\end{theorem}
\begin{proof} Voyons le point \emph{1}.
On a des parties finies  $J_0,J_1$ de $J$ et des parties finies
$U_0,U_1$ de $U$ telles que
$$ x,\; U_0 \vda J_0\quad{\rm  et}\quad   U_1  \vda x,\; J_1
$$
donc
$$ x,\; U_0,\; U_1 \vda J_0,\; J_1\quad{\rm  et}\quad  U_0,\; U_1 \vda
x,\;J_0,\; J_1
$$
Donc, par coupure,
$$  U_0,\; U_1 \vda J_0,\; J_1
$$
Le point \emph{2} a déjà été démontré (avec une formultaion
légèrement différente) dans la proposition~\ref{propIdealFiltre}.
\end{proof}

On notera que pour les \prois d'un \trdi le théorème de collapsus simultané est pratiquement la même chose que la règle de coupure.

Voici l'analogue du lemme de Krull \ref{corNstformHilbert}.
\begin{corollary}
\label{propTr2} \emph{(\nst formel pour les \prois dans un \trdi)}\\
 \Tcgi
\begin{enumerate}
\item Si $(J,U)$ est un  \proi qui ne s'effondre pas, il existe
$\varphi\in\Spec(T)$
tel que $J\subseteq \varphi^{-1}(0)$  et
$U\subseteq \varphi^{-1}(1)$.
En particulier: 
\item Si $b\not\leq a$, il existe $\varphi\in\Spec(T)$ tel que
$\varphi(a)=0$ et $\varphi(b)=1$.
\item Si $a\in T$, l'idéal $\,\dar a$ est intersection des \ideps qui contiennent $a$.
\item Si $T\neq \Un$ alors~$\Spec(T)$ est non vide.
\end{enumerate}
\end{corollary}

%
\begin{proof}
Il suffit de montrer le point \emph{1}. On peut recopier les démonstrations du corolaire \ref{corNstformHilbert}
(l'une utilisant le tiers exclu et le lemme de Zorn, l'autre utilisant le \tcg) en se basant sur le théorème de collapsus simultané.
Mais dans le cas des \trdis, la démonstration par le \tcg est plus directe:
il n'y a en effet pas besoin de faire référence au théorème de collapsus simultané,
qui n'est qu'une reformulation de la coupure, inhérente à la structure de \trdi.

\noindent \hum{expliquer cela en détail dans l'annexe?
\vspace{-.5em}}
\end{proof}

Voici l'analogue du corollaire \ref{corNstfH2}.

\begin{corollary} 
\label{corNstfH2Tr} 
Soit $\cP=(J,U)$  un \proi dans un \trdi $T$ et
 $(I,V)$ le saturé de~$\cP$. Le \tcg implique que $I$ est l'intersection 
des idéaux premiers qui contiennent~$J$ et ne coupent pas $U$,
tandis que $V$ est l'intersection des complémentaires de ces mêmes 
idéaux premiers. En particulier tout idéal de $T$ est intersection des \ideps qui le contiennent.
\end{corollary}

Un corolaire de \ref{propTr2} est le théorème de représentation suivant
(théorème de Birkhoff).
\begin{theorem}
\label{thRep} {\em (Théorème de représentation)}
\Tcgi L'application
$\theta_T: T\rightarrow \scP(\Spec(T))$ définie par
$a\mapsto \left\{\varphi\in \Spec(T)\; ;\; \varphi(a)=1 \right\}$
est un \homo injectif de \trdis. Autrement dit, tout \trdi peut être
représenté comme un sous-treillis du treillis des parties d'un ensemble.
\end{theorem}

Un autre corolaire est le suivant.
\begin{proposition}
\label{propRep2} \emph{(Lying over, \trdis)}
\Tcgi Soit $\varphi:T\rightarrow T'$ un \homo de \trdis. Alors $\varphi$
est injectif \ssi le morphisme dual $\Spec(\varphi):\Spec(T')\rightarrow \Spec(T)$ est
surjectif.
\end{proposition}
\begin{proof}
On a les équivalences:
$\quad  a\not= b\quad \Longleftrightarrow\quad  a\vi b\not= a\vu b\quad
\Longleftrightarrow\quad
a\vu b\not\leq a\vi b
$.
\\
Supposons que $\Spec(\varphi)$ est surjectif. Si $a\not= b$ dans $T$, soit
$a'=\varphi(a)$, $b'=\varphi(b)$ et
  $\psi\in\Spec(T)$ tel que $\psi(a\vu b)=1$ et
 $\psi(a\vi b)=0$. Puisque  $\Spec(\varphi)$ est surjectif il existe
$\psi'\in\Spec(T')$ tel \hbox{que $\psi=\psi'\varphi$} donc  $\psi'(a'\vu b')=1$ 
et
 $\psi'(a'\vi b')=0$, donc $a'\vu b'\not\leq a'\vi b'$ et $ a'\not= b'$.\\
Supposons que  $\varphi$ est injectif. Identifions $T$ à un sous-\trdi 
de $T'$. \\
Si $\psi\in\Spec(T)$, \hbox{soient $I=\psi^{-1}(0)$} et $F=\psi^{-1}(1)$.
Alors $(I,F)$ ne peut pas s'effondrer dans~$T'$ car cela le ferait s'effondrer
dans $T$. Donc il existe $\psi'\in\Spec(T')$ tel que $\psi'(I)= 0$ 
\hbox{et $\psi'(F)= 1$}, ce qui signifie~\hbox{$\psi=\psi'\varphi$}.
\end{proof}

Naturellement, il est difficile d'attribuer un contenu calculatoire
(autre que le théorème~\ref{thColSimT1}) aux quatre résultats
précédents des \clama.
Une solution intuitive est de dire qu'on ne risque s\^urement pas de se
tromper {\em en faisant comme si} un \trdi donné $T$ était un sous-treillis du treillis des parties d'un ensemble.
Le but du programme de Hilbert est de donner un contenu précis à
cette affirmation intuitive.
Que faut-il entendre précisément par
\gui{ne risquer s\^urement pas de se tromper} et
par \gui{en faisant comme si}?

\subsection{Chaînes partiellement spécifiées d'idéaux
premiers}\label{subsecproctreillis}
\begin{definition}{\rm 
\label{defprochTreil} Dans un \trdi $T$
\begin{enumerate}
\item  Une {\em spécification partielle pour une chaîne
d'\ideps} (ce que nous abrégeons en {\em \proc}) est définie comme
suit. Une \prolo $\ell$ est un liste de $\ell+1$ \prois de $T$:
$\cC=(\cC_0,\ldots,\cC_\ell)$. La \proc est dite finie si
tous
les $\cC_i$ sont finis. Une \prolo $0$ n'est autre qu'un
{\em \proi}.
\item  Une \proc est dite {\em saturée} si les $\cC_i$ sont saturés
 et si $\cC_i\subseteq \cC_{i+1}$ pour chaque~$i$.
\item  On dit qu'une \proc $\cC'= (\cC'_0,\ldots,\cC'_\ell)$
est
{\em un raffinement de la \proc~$\cC$}
si on a $\cC_k\leq  \cC'_k$ pour chaque $k$.
\item  On dit qu'une \proc $\cC$ {\em collapse} ou encore {\em
s'effondre}
si la seule \proc saturée qui raffine $\cC$ est la \proc triviale
$\big((T,T),\ldots,(T,T)\big)$.
\end{enumerate}
}
\end{definition}
\begin{lemma}
\label{lemColTr}
Une \proc $\cC=(\cC_0,\dots,\cC_\ell)$  s'effondre dès que l'un des $\cC_j$ s'effondre.
\end{lemma}
\begin{proof}
Soit $\big((I_0,F_0),\ldots,(I_\ell,F_\ell)\big)$ une \proc saturée qui
raffine $\cC$. Si le \proi  $(I_h,F_h)$  s'effondre, puisque  $I_h$  et
$F_h$ sont conjugués, on a $1\in  I_h$ et  $0\in  F_h$.
Pour tout indice $ j>h$  on a donc  $1\in  I_j$ donc   $0\in  F_j$.
De même pour tout indice $ j<h$  on a   $0\in  F_j$  donc
$1\in  I_j$.
\end{proof}

Dans le théorème suivant les points \emph{3} et \emph{2} sont analogues aux points \emph{1} 
et \emph{2} dans les théorèmes~\ref{thColSimT1} et \ref{ThColSimKrA}.
Notez qu'ici l'analogue des lemmes de finitude évidents \ref{lemfinicol}, \ref{lemfinicol2} et~\ref{lemfinicol3}, n'est pas évident.
Il résulte du theorème qui suit.
\begin{theorem}
\label{thColSimT2} {\em (Effondrement d'une \proc, saturation et collapsus simultané dans les
\trdis)} ~\\
Soit une  \proc $\cC=\big((J_0,U_0),\ldots,(J_\ell,U_\ell)\big)$ dans un \trdi
$T.$
\begin{enumerate}
\item  
\begin{itemize}
\item La  \proc $\cC$ s'effondre si elle raffine une \proc finie qui s'effondre.
\item Lorsque $\cC$ est finie elle s'effondre \ssi 
 il existe $x_1,\ldots,x_\ell\in T$ tels que:
\begin{equation} \label {eqthColSimT2}
\begin{array}{rcl}
 x_1,\;   U_0& \vda  &  J_0  \\
 x_2,\;   U_1& \vda  &  J_1 ,\;  x_1  \\
\vdots\quad & \vdots  & \quad\vdots  \\
 x_\ell,\;   U_{\ell-1}& \vda  &  J_{\ell-1} ,\;  x_{\ell-1}  \\
  U_\ell& \vda  &  J_\ell ,\;  x_\ell
\end{array}
\end{equation}

\end{itemize}
\item  La \proc $\cC$ engendre une \proc saturée minimum:
on ajoute dans~$U_i$  tout  élement $x$ de $T$  tel que la \proc
$\big((J_0,U_0),\ldots,(J_i,x;U_i),\ldots,(J_\ell,U_\ell)\big)$ s'effondre
et on ajoute dans $J_i$ tout élement $x$  tel que la \proc $\big((J_0,U_0),\ldots,(J_i;U_i,x),\ldots,(J_\ell,U_\ell)\big)$
s'effondre.
\item  Soit $x\in T$. Supposons que les \procs
$\big((J_0,U_0),\ldots,(J_i,x;U_i),\ldots,(J_\ell,U_\ell)\big)$
et
$\big((J_0,U_0),\ldots,(J_i;U_i,x),\ldots,(J_\ell,U_\ell)\big)$
s'effondrent toutes les deux, alors $\cC$ s'effondre également.
\end{enumerate}
\end{theorem}
\begin{proof}
Voyons d'abord les deux premiers points, qui se démontrent simultanément. Il n'est pas vraiment restrictif 
de supposer que la \proc $\cC$ est finie,
car une fois le résultat établi dans ce cas, on passe au cas
général en regardant la \proc donnée comme limite inductive
des \procs finies dont elle est un raffinement.
\\
Nous notons $\cC^1=\big((I_0,F_0),\ldots,(I_\ell,F_\ell)\big)$ la \proc
construite en~\emph{2}.
On va démontrer
\begin{itemize}
\item [$(\alpha)$] Si $\cC$ vérifie les inégalités (\ref{eqthColSimT2}),
toute \proc saturée qui raffine $\cC$ est triviale (i.e., $\cC$
s'effondre).
\item [$(\beta)$]  La \proc $\cC^1$ est bien  saturée.
\item [$(\gamma)$] Toute \proc saturée qui raffine $\cC$ raffine
également $\cC^1$.
\item [$(\delta)$] Si $\cC^1$ est triviale, $\cC$ vérifie les
inégalités (\ref{eqthColSimT2}).
\end{itemize}
Cela suffira bien à montrer \emph{1} et \emph{2}.\\
$(\alpha)$
Soit $\cC'=\big((I'_0,F'_0),\ldots,(I'_\ell,F'_\ell)\big)$ une \proc saturée qui
raffine $\cC$. Considérons les inégalités (\ref{eqthColSimT2})
$$\begin{array}{rcl}
 x_1,\; U_0& \vda  & J_0  \\
 x_2,\; U_1& \vda  & J_1 ,\; x_1  \\
\vdots\qquad & \vdots  & \qquad\vdots  \\
 x_\ell,\; U_{\ell-1}& \vda  & J_{\ell-1} ,\; x_{\ell-1}  \\
 U_\ell& \vda  & J_\ell ,\; x_\ell
\end{array}$$
Puisque  $I'_0$ et $F'_0$ sont
conjugués, la première inégalité donne $x_1\in I'_0$.
Donc  $x_1\in I'_1$,
et la deuxième inégalité donne  $x_2\in I'_1$.
En poursuivant, on obtient à la dernière inégalité
$U_\ell \vda   J_\ell ,\; x_\ell$. avec $x_\ell\in I'_\ell$ et $U_\ell\subseteq F'_\ell$. Donc $(I'_\ell,F'_\ell)$ s'effondre (lemme \ref{lemfinicol4}). Donc $\cC'$ s'effondre (lemme \ref{lemColTr}).
\\
$(\beta)$
Nous faisons la démonstration avec $\ell=3$.
Montrons d'abord que les $I_j$  sont des idéaux.
Nous faisons la démonstration avec $j=1$.
Pour montrer $0\in I_1$ on prend
\hbox{$x_1=0,\; x_2=x_3=1$}. De même pour montrer $J_1\subseteq  I_1$ on prend 
un
$x\in J_1$ et $x_1=0$, \hbox{$x_2=x_3=1$}.
Le fait que $x\in I_1$  et $y\leq x$ impliquent $y\in I_1$ est
immédiat:
on garde les mêmes $x_i$. Supposons maintenant $x,y\in I_1$ et
montrons
$x\vu y\in I_1$. Par hypothèse nous avons des $x_i$  et $y_i$
vérifiant
les inégalités suivantes.
$$\begin{array}{rclcrcl}
 x_1,\; U_0& \vda  & J_0  &\qquad &
 y_1,\; U_0& \vda  & J_0  \\
 x_2,\; U_1,\; x& \vda  & J_1 ,\; x_1  &&
 y_2,\; U_1,\; y& \vda  & J_1 ,\; y_1  \\
 x_3,\; U_{2}& \vda  & J_{2} ,\; x_{2}  &&
 y_3,\; U_{2}& \vda  & J_{2} ,\; y_{2}  \\
 U_3& \vda  & J_3 ,\; x_3  &&
 U_3& \vda  & J_3 ,\; y_3
\end{array}$$
On en déduit, en utilisant la distributivité
$$\begin{array}{rclcrcl}
 (x_1\vu y_1),\; U_0& \vda  & J_0   \\
 (x_2\vi y_2),\; U_1,\; (x\vu y)& \vda &J_1,\;(x_1\vu y_1)\\
 (x_3\vi y_3),\; U_{2}& \vda  & J_{2} ,\; (x_2\vi y_2)  \\
 U_3& \vda  & J_3 ,\; (x_3\vi y_3)  \\
\end{array}$$
Montrons maintenant que idéaux et filtres sont conjugués. Par exemple
avec $I_1$  et $F_1$. Nous supposons  $x\vi y\in I_1,\; y\in F_1$ et
nous montrons  $x\in I_1$.
Par hypothèse nous avons des $x_i$  et $y_i$  vérifiant les
inégalités suivantes.
$$\begin{array}{rclcrcl}
 x_1,\; U_0& \vda  & J_0  &\qquad &
 y_1,\; U_0& \vda  & J_0  \\
 x_2,\; U_1,\; (x \vi y)& \vda  & J_1 ,\; x_1  &&
 y_2,\; U_1& \vda  & J_1 ,\; y_1 ,\; y \\
 x_3,\; U_{2}& \vda  & J_{2} ,\; x_{2}  &&
 y_3,\; U_{2}& \vda  & J_{2} ,\; y_{2}  \\
 U_3& \vda  & J_3 ,\; x_3  &&
 U_3& \vda  & J_3 ,\; y_3
\end{array}$$
On en déduit, en utilisant la distributivité
$$\begin{array}{rclcrcl}
 (x_1\vu y_1),\; U_0& \vda  & J_0   \\
 (x_2\vi y_2),\; U_1,\; x,\; y& \vda&J_1,\;(x_1\vu y_1)\\
 (x_2\vi y_2),\; U_1,\; x& \vda&J_1,\;(x_1\vu y_1),\; y\\
 (x_3\vi y_3),\; U_{2}& \vda  & J_{2} ,\; (x_2\vi y_2)  \\
 U_3& \vda  & J_3 ,\; (x_3\vi y_3)
\end{array}$$
Les inégalités \num 2 et 3 donnent par coupure
$$\begin{array}{rclcrcl}
 (x_2\vi y_2),\; U_1,\; x& \vda &J_1,\;(x_1\vu y_1)\\
\end{array}$$
La démonstration est complète.
\\
$(\gamma)$
Nous faisons la démonstration avec $\ell=3$.
Soit $\big((I'_0,F'_0),\ldots,(I'_3,F'_3)\big)$ une \proc
saturée qui raffine $\cC$. Montrons par exemple $I_1\subseteq I'_1$.
Soit $x\in I_1$, on a donc
$$\begin{array}{rclcrcl}
 x_1,\; U_0& \vda  & J_0   \\
 x,\; x_2,\; U_1& \vda  & J_1 ,\; x_1  \\
 x_3,\; U_{2}& \vda  & J_2 ,\; x_2  \\
 U_3& \vda  & J_3 ,\; x_3
\end{array}$$
On en déduit successivement $x_1\in I'_0\subseteq I'_1$,
 $x_3\in F'_3\subseteq F'_2$,
$x_2\in F'_2\subseteq F'_1$,  et enfin
$x\in I'_1$. Remarquez que la démonstration du point $(\alpha)$ peut être vue
comme un cas particulier de celle du point~$(\gamma)$.
\\
Le point $(\delta)$     est immédiat

\sni Prouvons enfin \emph{3}. Nous avons $x\in I_i$  et  $x\in F_i$,  donc
$\cC^1$ s'effondre (lemme \ref{lemColTr}). Donc  $\cC$ s'effondre.
\end{proof}
\hum{Ici on n'a pas eu besoin d'utiliser une \recu sur $\ell$ ? Cela contraste avec la démonstration du théorème \ref{thKrJoy}.}

\begin{corollary} \label{corthColSimT2}
Soit $\cC=(\cP_0,\dots,\cP_k)$ une \proc saturée et $\cQ$ un \proi saturé tel
que $\cP_k\subseteq \cQ$. La \proc $(\cP_0,\dots,\cP_k,\cQ)$ s'effondre \ssi
la \proc $(\cP_k,\cQ)$ s'effondre.  
\end{corollary}
Notez que  $(\cP_0,\dots,\cP_k,\cQ)$ est saturée par définition.

\begin{definition}{\rm 
\label{defproch3T} 
 Deux \procs qui engendrent la même \proc saturée sont dites
{\em
équivalentes}.
 Une \proc équivalente à une \proc finie est appelée une {\em
\proc de type fini}.
}
\end{definition}

Nous regardons une \proc $\cC=(\cP_0,\dots,\cP_k)$  comme une spécification partielle
pour une chaîne croisssante d'\ideps $P_0,\ldots,P_k$ où chaque $P_i$
raffine $\cP_i$.
Du théorème de collapsus simultané on déduit le résultat
suivant qui justifie cette idée de spécification partielle.
\begin{theorem}
\label{th.nstformelTreil} {\em (\nst formel pour les chaînes
d'\ideps dans un \trdi)} \Tcgi Soit $T$ un \trdi et
$\big((J_0,U_0),\ldots,(J_\ell,U_\ell)\big)$ une \proc dans $T$. \Propeq
\begin{itemize}
\item [$(a)$] Il existe $\ell+1$ idéaux premiers $P_0\subseteq
\cdots\subseteq
P_\ell$  tels que
 $J_i\subseteq P_i$, $U_i\cap P_i=\emptyset $, $(i=0,\ldots,\ell)$.
\item [$(b)$] La \proc ne s'effondre pas.
\end{itemize}
\end{theorem}
La démonstration est essentiellement la même que celle du théorème \ref{th.nstformel}. \hum{à vérifier quand même}

\subsection{Définition constructive de la dimension d'un \trdi}
\label{subsubsecDimTrdi}

On peut développer une théorie \cov de la dimension de Krull d'un \trdi selon les mêmes lignes que ce qui a été
développé pour les anneaux commutatifs dans la section~\ref{secKrA},
c'est ce que nous faisons dans la section présente. 

On peut aussi  utiliser une idée de Joyal:
construire un treillis $\Kr_\ell(T)$ universel
attaché à $T$  tel que les points de $\Spec(\Kr_\ell(T))$  soient
(de manière naturelle) les chaînes d'idéaux premiers de longueur
$\ell$.

La théorie de Joyal sera expliquée dans la section~\ref{secJoyal}.
On montrera que les deux points de vue sont
isomorphes.
\begin{definition}{\rm 
\label{defDiTr}~
\begin{enumerate}
\item Une {\em \proel}  dans un \trdi $T$ est une \proc de la forme suivante
$$ \big((0;x_0),(x_0;x_1),\ldots,(x_\ell;1)\big)
$$
(avec les $x_i$ dans $T$). On la note $\overline{(x_0,\ldots,x_\ell)}$.

\item Un \trdi $T$ est dit {\em de dimension $\leq \ell$}
si toute \proel $\overline{(x_0,\ldots,x_\ell)}$ de longueur $1+\ell$ s'effondre, i.e. s'il existe $ y_0,\dots,y_\ell\in T$ tels 
que
\begin{equation} \label {eqdefDiTr}
\begin{array}{rclll}
1 & \vda  &  y_\ell,\;x_\ell    \\
y_\ell,\;x_\ell& \vda  &y_{\ell-1},\;x_{\ell-1}      \\
\vdots\qquad & \vdots  & \qquad  \vdots    \\
y_1,\;x_1& \vda  &  y_0,\;x_0    \\
y_0,\;x_0& \vda  &  0    \\
\end{array}
\end{equation}
On notera $\dim\, T\leq \ell$. Par exemple pour la dimension $2$ cela correspond au dessin suivant dans $T$
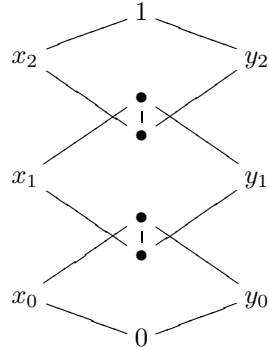
$$\SCO{x_0}{x_1}{x_2}{y_0}{y_1}{y_2}$$ 
\item Deux suites $(x_0,\dots,x_\ell)$ et $(y_0,\dots,y_\ell)$ qui vérifient les relations du point 2 sont dites \emph{complémentaires}.
\end{enumerate}
}
\end{definition}

 En particulier le \trdi $T$ est de dimension $\leq -1$ \ssi $1=0$
dans~$T$, et  de dimension
$\leq 0$ \ssi $T$ est une algèbre de Boole (tout élément de
$T$ a un complément).

Nous ne donnerons pas dans le cadre des \trdis de définition pour
$\dim(T)<\ell$ ni pour $\dim(T)\geq \ell$ et nous nous contenterons de 
dire
que $\dim(T)>\ell$ signifie l'impossibilité de $\dim(T)\leq \ell$.
On pourra introduire les mêmes raffinements
que dans la section \ref{subsecPsr}
lorsqu'on a une relation d'inégalité dans $T$ qui n'est pas la simple
impossibilité de l'égalité.

\smallskip Le corollaire suivant du \nst formel \ref{th.nstformelTreil} nous dit que la définition \cov est bien équivalente à la définition usuelle donnée en \clama.

\begin{theorem}
\label{th.defdimTr} Le \tcg implique que \propeq
\begin{itemize}
\item [$(a)$] Les chaînes strictement croissantes d'\ideps ont une longeur $< \ell$.
\item [$(b)$] Toute \proelo $\ell$ s'effondre.
\end{itemize}
\end{theorem}


\section{La Théorie de Joyal pour la dimension de Krull}
\label{secJoyal}

\subsection{Treillis de Krull associés à un \trdi}
\label{subsecJoyalKrull}

L'idée de Joyal est de construire un treillis $\Kr_\ell(T)$ universel
attaché à $T$  tel que les points de $\Spec(\Kr_\ell(T))$  soient
(de manière naturelle) les chaînes croissantes d'idéaux premiers 
de longueur~$\ell$.

Donner une telle chaîne revient à donner des \homos
$\mu_0\geq \mu_1\geq \cdots\geq  \mu_\ell$ de~$T$~vers~$\Deux$.

Si donc on a un \trdi $K_\ell$ et $\ell+1$ \homos
$\varphi_0\geq \varphi_1\geq \cdots\geq  \varphi_\ell$ de $T$ vers $K_\ell$ 
tels que, pour tout treillis $ T'$ et tous  
$\psi_0\geq \psi_1\geq \cdots\geq\psi_\ell\in\Hom(T,T')$
on ait une \homo unique  $\eta :K_\ell\rightarrow T'$  vérifiant 
$\eta\varphi_0=\psi_0$,
$\eta \varphi_1=\psi_1,$ $\ldots$,  $\eta \varphi_\ell=\psi_\ell$, 
alors on aura les éléments de $\Spec(K_\ell)$ qui s'identifient 
naturellement aux chaînes d'\ideps de longueur $\ell$ dans $T$.

L'avantage est que $K_\ell$ est toujours un objet qu'on peut construire
explicitement à partir de $T$, contrairement aux idéaux premiers et
aux spectres.
En effet le fait qu'un tel objet universel  
$$
(\Kr_\ell(T),\varphi_0, \dots,  \varphi_\ell)
$$ 
existe toujours de manière
\cov résulte de considérations générales d'algèbre
universelle\footnote{Ou encore de la description des problèmes universels qui admettent une solution dans la théorie des structures de Bourbaki.}.

Résumons la discussion précédente.

\begin{lemma} \label{lemKrJoyal} {\em (Un autre \nst formel pour les chaînes
d'\ideps dans un \trdi)}
Soit $(J_k)_{0\leq k\leq \ell}$ et $(U_k)_{0\leq k\leq \ell}$ des parties d'un \trdi $T$ non trivial. \\
Il revient au même de se donner 
\begin{itemize}
\item une chaîne d'\ideps détachables $P_k$ de $T$
qui vérifient $J_k\subseteq P_k$ et $U_k\cap P_k=\emptyset$,
\item ou un homomorphisme  $\psi:\Kr_\ell(T)\to \Deux$ (i.e. un point de $\Spec(\Kr_\ell(T))$) tel que \hbox{$\psi(\varphi_k(J_k))=0$} et $\psi(\varphi_k(U_k))=1$ pour $0\leq k\leq \ell$. 
\end{itemize}
\end{lemma}

Ceci motive la définition \cov alternative suivante.

\begin{defa} \label{defDiTr2} (alternative à la définition \ref{defDiTr})
\\
{\rm Pour $0\leq i\leq \ell+1$, notons $\varphi_i:T\to \Kr_{\ell+1}(T)$  les $\ell+2$ morphismes universels.
\\
On dit que le \trdi $T$ est de dimension $\leq \ell$ si pour tous éléments \hbox{$x_1,\dots,x_{\ell+1}\in T$}, on a
$
\varphi_0(x_1),\dots,\varphi_{\ell}(x_{\ell+1})
\vda
\varphi_1(x_1),\dots,\varphi_{\ell+1}(x_{\ell+1})
$
dans $\Kr_{\ell+1}(T)$. 
}
\end{defa}

Et l'on obtient comme corolaire du lemme \ref{lemKrJoyal} la proposition suivante.

\begin{proposition} \label{proplemKrJoyal} \Tcgi
 Un \trdi $T$ est de dimension $\leq \ell$ au sens usuel des \clama \ssi il satisfait la définition constructive \ref{defDiTr2}.
 \end{proposition}
%
\begin{proof}
La définition classique usuelle dit qu'il n'existe pas de chaîne d'\ideps
$P_0,\dots,P_{\ell+1}$ avec $x_1,\dots,x_{\ell+1}$ satisfaisant $x_1\in P_1\setminus P_0$, $x_2 \in P_2\setminus P_1$, \dots, $x_{\ell+1} \in P_{\ell+1}\setminus P_{\ell}$.
Autrement dit, d'après le \tcg et le lemme \ref{lemKrJoyal}, dans le treillis $\Kr_{\ell+1}(T)$, le \proi 
$$
(\varphi_1(x_1),\dots,\varphi_{\ell+1}(x_{\ell+1});\varphi_0(x_1),\dots,\varphi_{\ell}(x_{\ell+1}))
$$ 
doit s'effondrer. Cela signifie d'après le lemme  \ref{lemfinicol4} que
l'on a 
$$
\varphi_0(x_1),\dots,\varphi_{\ell}(x_{\ell+1})
\vdash_{\Kr_{\ell+1}(T)}
\varphi_1(x_1),\dots,\varphi_{\ell+1}(x_{\ell+1}).
$$
\end{proof}

La définition alternative est assez sympathique parce qu'elle permet
d'obtenir facilement la caractérisation plus simple suivante dans laquelle intervient seulement un système générateur $S$ du treillis.

\begin{lemma}
\label{lemDimGen}
Un \trdi $T$ engendré par une partie $S$ est  de dimension~$\leq \ell$
(définition~\ref{defDiTr2}) \ssi  pour tous $x_1,\dots,x_{\ell+1}\in S$, on a
$$
\varphi_0(x_1),\dots,\varphi_{\ell}(x_{\ell+1})
\vdash_{\Kr_{\ell+1}(T)}
\varphi_1(x_1),\dots,\varphi_{\ell+1}(x_{\ell+1}).
$$
\end{lemma}
En effet les règles de distributivité permettent par exemple de
déduire
$$ a\vu a',A\vda b\vu b',B$$
de $a,A\vda b,B$ et $a',A\vda b',B$. Par ailleurs tout élément de $T$
est un inf de sups d'éléments de $S$.

\subsection{Comparaison des deux points de vue}
\label{subsubsecCompar}

La description explicite de $\Kr_\ell(T)$ a été simplifiée par la
considération des \entrels (\cite[section 4.4]{cc}). Plus précisément on a
le théorème suivant.

\begin{theorem}
\label{thKrJoy}
Soit $T$ un \trdi. On considère le problème universel suivant, que
nous appelons \gui{problème de Krull}:
trouver  un \trdi $K_\ell$ et $\ell+1$ \homos
$\varphi_0\geq \varphi_1\geq \cdots\geq  \varphi_\ell$
de $T$ vers $K_\ell$ tels que,
pour tout treillis $ T'$ et tous
$\psi_0\geq \psi_1\geq \cdots\geq  \psi_\ell\in\Hom(T,T')$
on ait un \homo unique  $\eta :K_\ell\rightarrow T'$  vérifiant
$\eta \varphi_0=\psi_0$,
$\eta \varphi_1=\psi_1,$ $\ldots$,  $\eta \varphi_\ell=\psi_\ell$.
Ce problème universel admet une solution (unique à \iso unique 
près).
 On note  $\Kr_\ell(T)$ le \trdi correspondant.
Il peut être décrit comme le treillis engendré par la réunion
disjointe  $S_\ell$  de $\ell+1$  copies de $T$  (on note $\varphi_i$ la
bijection de $T$ sur la copie indexée par $i$) munie de la \entrel
$\vdash_{S_\ell}$ définie comme suit.
Si $A_i$ et $B_i$  $(i=0,\ldots,\ell)$ sont des parties finies de~$T$
on a
$$ \varphi_0(A_0),\ldots,\varphi_\ell(A_\ell) \,\vdash_{S_\ell}\,
\varphi_0(B_0),\ldots,\varphi_\ell(B_\ell) \eqno(*)
$$
\ssi il existe $u_1,\ldots,u_\ell\in T$ avec les relations suivantes dans
$T$ (où $\vda$ est  la \entrel de $T$):
$$\begin{array}{rcl}
u_1,\;  A_0& \vda  &  B_0  \\
u_2,\;  A_1& \vda  &  B_1 ,\; u_1  \\
    \vdots\quad & \vdots  & \quad\vdots  \\
u_\ell,\;  A_{\ell-1}& \vda& B_{\ell-1} ,\; u_{\ell-1}  \\
 A_\ell& \vda  &  B_\ell ,\; u_\ell  \\
\end{array} \eqno (**)$$
\end{theorem}
\begin{proof}
On montre d'abord que la relation $\vdash_{S_\ell}$ définie sur $\Pfe({S_\ell})$ dans 
le
théorème est bien une \entrel.
Le seul point délicat est la règle de coupure. Pour simplifier
les notations, on fait la démonstration avec $\ell=3.$
Il y a alors 3 cas possibles, et on analyse un des cas possibles, celui où
$X,\varphi_1(z)\vdash_{S_3} Y$ et $X\vdash_{S_3} Y,\varphi_1(z)$, les autres cas
étant semblables.
Par hypothèse on a $u_1,u_2,u_3,v_1,v_2,v_3$ tels que
$$\begin{array}{rclcrcl}
 u_1,\; A_0& \vda  & B_0  &\qquad &
 v_1,\; A_0& \vda  & B_0  \\
 u_2,\; A_1,\; z& \vda  & B_1 ,\; u_1  &&
 v_2,\; A_1& \vda  & B_1 , v_1 ,\; z \\
 u_3,\; A_{2}& \vda  & B_{2} ,\; u_{2}  &&
 v_3, A_{2}& \vda  & B_{2} , v_{2}  \\
 A_3& \vda  & B_3 ,\; u_3  &&
 A_3& \vda  & B_3 ,\; v_3
\end{array}$$
Les deux \entrels sur la deuxième ligne donnent
$$\begin{array}{rclcrcl}
 u_2,\;v_2,\; A_1,\; z& \vda  & B_1 ,\; u_1,\;v_1 &\qquad
 u_2,\;v_2,\; A_1& \vda  & B_1 ,\; u_1 ,\; v_1 ,\; z \end{array}$$
donc par coupure
$$\begin{array}{rclcrcl}
 u_2,\;v_2,\; A_1& \vda  & B_1 ,\; u_1,\;v_1
\end{array}$$
\cad
$$\begin{array}{rclcrcl}
 u_2\vi v_2,\; A_1& \vda  & B_1 ,\; u_1\vu v_1
\end{array}$$
Finalement,  en utilisant la distributivité
$$\begin{array}{rclcrcl}
 (u_1\vu v_1),\; A_0& \vda  & B_0   \\
 (u_2\vi v_2),\; A_1& \vda&B_1,\;(u_1\vu v_1)\\
 (u_3\vi v_3),\; A_{2}& \vda  & B_{2},\;(u_2\vi v_2)  \\
 A_3& \vda  & B_3, \;(u_3\vi v_3)
\end{array}$$
et donc $\varphi_0(A_0),\dots,\varphi_3(A_3)\,\vdash_{S_3}\,
\varphi_0(B_0),\dots,\varphi_3(B_3)$.\\
Il reste à voir que le treillis  $\Kr_\ell(T)$ défini à partir de
$({S_\ell},\vdash_{S_\ell})$ vérifie bien la propriété universelle voulue. La
solution du problème universel posé existe pour des raisons
générales d'algèbre universelle.
Et on sait que  $\Kr_\ell(T)$  est nécessairement engendré par ${S_\ell}$. Il
suffit donc de définir  sur $S_\ell$ la \entrel la moins contraignante possible
(avec la condition que les $\varphi_i$ forment une suite décroissante
d'\homos).

\noindent Montrons donc que pour une \entrel $\vdash_\ell$ adéquate sur  $S_\ell$,
la relation $(*)$ est toujours
satisfaite lorsque sont satisfaites les relations $(**)$ dans $T$.
Notons $X=\varphi_0(A_0),\ldots,\varphi_\ell(A_\ell)$ \hbox{et $Y=\varphi_0(B_0),\ldots,\varphi_\ell(B_\ell)$}. Par monotonie
et puisque les $\varphi_i$ sont des homomorphismes,	on déduit  par exemple des relations $A_0,u_1\vda B_0\,$ et  $\,A_1,u_2\vda B_1,u_1$ les relations
$$
X,\varphi_0(u_1)\vdash_\ell Y\;\hbox{  et  }\;X,\varphi_1(u_2)\vdash_\ell Y,\varphi_1(u_1)
.$$
Comme $\varphi_0(u_1)\geq \varphi_1(u_1) $, on obtient
par coupures $X,\varphi_1(u_2)\vdash_\ell Y,\varphi_0(u_1)$ puis $X,\varphi_1(u_2)\vdash_\ell Y$. Nous avons gagné un cran. Nous pouvons continuer jusqu'à obtenir $X\vdash_\ell Y$. 
\end{proof}

On note que les \homos $\varphi_i$ sont injectifs: on voit facilement que
pour $a,b \in T$ la relation $\varphi_i(a)\vdash_{S_\ell}\varphi_i(b)$ implique
$a\vda b$, donc $\varphi_i(a)=\varphi_i(b)$ implique $a=b$.

\smallskip 
On obtient le corollaire suivant des théorèmes
\ref{thColSimT2} et \ref{thKrJoy}.

\begin{corollary} \label{corthKrJoy}
Une \proc $\cC=\big((J_0,U_0),\ldots,(J_\ell,U_\ell)\big)$ s'effondre dans $T$ \ssi le \proi
 $\cP=\big(\varphi_0(J_0),\ldots,\varphi_\ell(J_\ell);
\varphi_0(U_0),\ldots,\varphi_\ell(U_\ell)\big)$
s'effondre dans $\Kr_\ell(T)$. En particulier les définitions
\covs pour la dimension de Krull, \ref{defDiTr} et \ref{defDiTr2}, sont équivalentes
\end{corollary}

En \clama le corolaire \ref{corthKrJoy} est assuré comme conséquence du théorème \ref{th.nstformelTreil}  et du lemme \ref{lemKrJoyal}. Cela permet donc de
déduire l'un de l'autre les deux théorèmes  \ref{thColSimT2}~\emph{1} et \ref{thKrJoy}.
On est d'ailleurs frappé par la ressemblance des démonstrations des théorèmes
\ref{thColSimT2} et \ref{thKrJoy}.

Pour établir directement le corollaire \ref{corthKrJoy}, i.e. l'équivalence des deux théorèmes \ref{thColSimT2}~\emph{1} et \ref{thKrJoy}  il suffit d'expliquer comment
la donnée d'une \proc saturée  $\big((I_0,F_0),\ldots,(I_\ell,F_\ell)\big)$   
de $T$ permet de produire une suite décroissante d'\homos 
$(\psi_0,\ldots,\psi_\ell)$ de $T$ dans un
\trdi  avec  $\psi_k^{-1}(0)=I_k$ et
$\psi_k^{-1}(1)=F_k$ ($k=0,\ldots,\ell$). 
On doit donc donner une démonstration directe du 
lemme
suivant.
\begin{lemma}
\label{lemProcEtKrl}
%
%
Soit $\;\cC=\big((I_0,F_0),\ldots,(I_\ell,F_\ell)\big)$ une \proc saturée dans un \trdi~$T$. 
Soit  $T_\cC$ le \trdi quotient de  $\Kr_\ell(T)$ obtenu en forçant
$$\varphi_0(I_0)=\cdots=\varphi_\ell(I_\ell)=0,\;
\varphi_0(F_0)=\cdots=\varphi_\ell(F_\ell)=1.$$ Soit $\pi$ la projection
canonique de $\Kr_\ell(T)$ sur $T_\cC$ et $\psi_k=\pi\circ\varphi_k$.
Alors $\psi_k^{-1}(0)=I_k$ et $\psi_k^{-1}(1)=F_k$ \hbox{pour $k=0,\ldots,\ell)$}.

%
%
%

\end{lemma}
\begin{proof}
Par exemple l'idéal $\psi_k^{-1}(0)=\left\{x\in T;\varphi_k(x)=_{T_\cC}0\right\}$
est, par application de la proposition \ref{propIdealFiltre} égal~à
$$ \left\{x\in 
T\;;\;\varphi_k(x),\varphi_0(F'_0),\ldots,\varphi_\ell(F'_\ell)
\;\vdash_{\Kr_\ell(T)}\; \varphi_0(I'_0),\ldots,\varphi_\ell(I'_\ell) 
\right\}
$$
pour des $I'_k\in\Pfe(I_k)$ et $F'_k\in\Pfe(F_k)$. 
Ainsi $\psi_k^{-1}(0)$ est l'ensemble des $x$ tels qu'il existe $x_1,\ldots,x_\ell$ tels que
$$\begin{array}{rcl}
x_1,\;  F'_0& \vdash_T  &  I'_0  \\
x_2,\;  F'_1& \vdash_T  &  I'_1 ,\; x_1  \\
    \vdots\quad & \vdots  & \quad\vdots  \\
x,x_{k+1},\;  F'_k& \vdash_T  &  I'_k ,\; x_k  \\
    \vdots\quad & \vdots  & \quad\vdots  \\
x_\ell,\;  F'_{\ell-1}& \vdash_T& I'_{\ell-1} ,\; x_{\ell-1}  \\
 F'_\ell& \vdash_T  &  I'_\ell ,\; x_\ell  \\
\end{array}$$
Comme la \proc $\cC$ est saturée on a de proche en proche
$x_1\in I_1\subseteq I_2$, $x_2\in I_2$, \dots,  $x_k\in I_k$, et  
$x_\ell\in
F_\ell$, \dots,
$x_{k+1}\in F_{k+1}\subseteq F_k$, d'où enfin $x\in I_k$.
\end{proof}

\subsection{Connections avec la définition de Joyal}
\label{subsubsecJoyal2}

 Si $T$ est un treillis distributif, Joyal \cite{esp} donne
la définition suivante de $\dim(T)\leq \ell$. \\
On note
$\varphi^\ell_i:T\rightarrow \Kr_\ell(T)$ les $\ell+1$ morphismes
universels ($0\leq i\leq \ell$). 
Par universalité de $\Kr_{\ell+1}(T)$, on
obtient $\ell+1$ morphismes $\sigma_i^\ell:\Kr_{\ell+1}(T)\rightarrow
\Kr_\ell(T)$
définis par les égalités 
$$
\sigma_i^\ell\circ \varphi^{\ell+1}_j = \varphi^\ell_j\;\hbox{  si  }\;j\leq i\;\;\;\hbox{ 
et }\;\;\;\sigma_i^\ell\circ \varphi^{\ell+1}_j = \varphi^\ell_{j-1}\;\hbox{  si  }\;j>i.
$$

Joyal définit alors $\dim(T)\leq \ell$ par le fait que
$\sigma^\ell:=(\sigma_0^\ell,\dots,\sigma_\ell^\ell):\Kr_{\ell+1}(T)\rightarrow
\Kr_\ell(T)^{\ell+1}$
est injectif. Cette définition peut se motiver à l'aide de
la proposition \ref{propRep2} qui dit que $\sigma^\ell$ est injectif \ssi 
le morphisme dual $\Spec(\sigma^\ell): \bigcup_i\Spec(K_\ell)\to K_{\ell+1}$ est surjectif. Les éléments dans l'image
de $\Spec(\sigma_i^\ell)$ sont les cha\^ines d'\ideps\footnote{On rappelle qu'un \idep $T$ est ici identifié à un morphisme $\alpha:T\to \Deux$.} 
$(\alpha_0,\dots,\alpha_{\ell+1})$
avec $\alpha_i=\alpha_{i+1}$, et $\Spec(\sigma_0^\ell,\dots,\sigma_\ell^\ell)$
est surjectif \ssi pour toute cha\^ine $(\alpha_0,\dots,\alpha_{\ell+1})$
il existe $i<\ell$ tel que $\alpha_i=\alpha_{i+1}$. Ceci dit donc
exactement qu'il n'y a pas de cha\^ines d'\ideps 
 de longueur~\hbox{$\ell+1$}. En utilisant le \tcg, on voit donc
l'équivalence de la définition classique avec la définition  \ref{defDiTr2} (qui fait le lien avec la définition \ref{defDiTr}
via le théorème~\ref{thKrJoy}). 

Nous donnons maintenant une démonstration directe de l'équivalence
de la définition de Joyal avec la
 la définition  \ref{defDiTr2}.

\begin{theorem}
On a $\dim(T)\leq \ell$, au sens de la définition  \ref{defDiTr2}
\ssi l'homomorphisme $(\sigma_0^\ell,\dots,\sigma_\ell^\ell):\Kr_{\ell+1}(T)\rightarrow
\Kr_\ell(T)^{\ell+1}$
est injectif.
\end{theorem}

\begin{proof}
 Pour simplifier les notations, nous nous limitons au cas
$\ell = 2$ et nous écrivons $\phi_i$ pour $\varphi^2_i$ et
$\psi_i$ pour $\varphi^3_i$. Nous notons aussi $K_\ell$ pour $\Kr_\ell(T)$
et $\sigma_i$ pour $\sigma_i^\ell$.

\noindent     Si $(\sigma_0,\sigma_1,\sigma_2)$ est injectif on montre que l'on a
$$
\psi_0(x_1),\psi_1(x_2),\psi_{2}(x_{3})\,\vdash_{K_3}\,
\psi_1(x_1),\psi_2(x_2),\psi_{3}(x_{3})
$$
pour toute suite $x_1,x_2,x_{3}$ en observant que l'on a
pour $i=0,1,2$
$$
\sigma_i\psi_0(x_1),\sigma_i\psi_1(x_2),\sigma_i\psi_{2}(x_{3})
\,\vdash_{K_2}\,
\sigma_i\psi_1(x_1),\sigma_i\psi_2(x_2),\sigma_i\psi_{3}(x_{3})
$$
car $\sigma_i\psi_{i+1} = \phi_i = \sigma_i\psi_{i}$, d'où $\sigma_i\psi_{i+1}(x_{i+1})  = \sigma_i\psi_{i}(x_{i+1})$: il y a donc un terme à droite égal à un terme à gauche.\\
Réciproquement, supposons que l'on a
$$
\psi_0(x_1),\psi_1(x_2),\psi_{2}(x_{3})
\,\vdash_{K_3}\,
\psi_1(x_1),\psi_2(x_2),\psi_{3}(x_{3})
$$
pour toute suite $x_1,x_2,x_{3}$ et montrons que
$(\sigma_0,\sigma_1,\sigma_{2})$ est injectif.
On doit montrer que l'on a $X\,\vdash_{K_3}\, Y$ si
$\sigma_i(X)\,\vdash_{K_2}\,\sigma_i(Y)$ pour chaque $i$.
Comme $K_3$ est engendré par $\bigcup_{0\leq i\leq 2} \psi_i(T)$
il suffit en fait de montrer que si l'on a
$$
\sigma_i\psi_0(a_0),\dots,\sigma_i\psi_{3}(a_{3})
\,\vdash_{K_2}\,
\sigma_i\psi_0(b_0),\dots,\sigma_i\psi_{3}(b_{3})
$$
pour chaque $i$ alors on a
$$
\psi_0(a_0),\dots,\psi_{3}(a_{3})
\,\vdash_{K_3}\,
\psi_0(b_0),\dots,\psi_{3}(b_{3})
$$

\noi Pour ceci, on remarque que l'on a par hypothèse
$$
\psi_0(a_1),\dots,\psi_{2}(a_{3})
\,\vdash_{K_3}\,
\psi_1(a_1),\dots,\psi_{3}(a_{3})
$$
et donc, il suffit de montrer
$$
(\psi_j(a_j))_{0\leq j<i},\,\psi_{i+1}(a_i\vi a_{i+1}),\,
(\psi_j(a_j))_{i+1< j\leq 3}
\,\vdash_{K_3}\,
\psi_0(b_0),\,\dots,\,\psi_{3}(b_{3})
$$
pour chaque $i=1,2,3$. Montrons par exemple le cas $i=1$. Soit
$$
\psi_0(a_0),\psi_2(a_1\vi a_2),\psi_3(a_3)\,\vdash_{K_3}\, 
\psi_0(b_0),\dots,\psi_{3}(b_{3})
$$
On a par hypothèse
$$
\sigma_1\psi_0(a_0),\dots,
\sigma_1\psi_{3}(a_{3})
\,\vdash_{K_2}\,
\sigma_1\psi_0(b_0),\dots,
\sigma_1\psi_{3}(b_{3})
$$
ce qui s'écrit
$$
\phi_0(a_0),\phi_1(a_1),\phi_1(a_2),\phi_2(a_3)
\,\vdash_{K_2}\,
\phi_0(b_0),\phi_1(b_1),\phi_1(b_2),\phi_2(b_3)
$$
soit
$$
\phi_0(a_0),\phi_1(a_1\vi a_2),\phi_2(a_3)
\,\vdash_{K_2}\,
\phi_0(b_0),\phi_1(b_1\vu b_2),\phi_2(b_3)
$$
Par universalité de $\phi_0\geq \phi_1\geq \phi_2$ appliquée
à $\psi_0\geq\psi_2\geq\psi_3$ ceci entraine
$$
\psi_0(a_0),\psi_2(a_1\vi a_2),\psi_3(a_3)
\,\vdash_{K_3}\,
\psi_0(b_0),\psi_2(b_1\vu b_2),\psi_3(b_3)
$$
et comme $\psi_2(b_1)\,\vdash_{K_3}\, \psi_1(b_1)$
ceci entraine bien
$$
\psi_0(a_0),\psi_2(a_1\vi a_2),\psi_3(a_3)
\,\vdash_{K_3}\,
\psi_0(b_0),\psi_1(b_1),\psi_2(b_2),\psi_3(b_3)
$$
\end{proof}

\subsection{Connections avec le travail d'Espa\~nol}
\label{subsubsecEspanol}

 Soit $T$ un treillis distributif,
 Espa\~nol \cite{esp} donne une caractérisation élégante de
$\dim(T)\leq {\ell}-1$ en terme de l'algèbre de Boole
engendrée par $T$. Le but de cette section est de
présenter cette caractérisation et de montrer
son équivalence avec la définition \ref{defDiTr}.

\begin{lemma}\label{lattice3}
Si $x_1\leq x_2\leq \dots\leq x_\ell$ et $a_1,\dots,a_\ell$ sont deux suites complémentaires
et si on pose $b_\ell = a_\ell$, $b_{\ell-1}=a_\ell\vi a_{\ell-1}$, \dots, alors $b_1\leq b_2\leq \dots\leq b_\ell$ et les deux suites $(x_1,\dots,x_\ell)$ et $(b_1,\dots,b_\ell)$ sont complémentaires. 
\end{lemma}
\begin{proof} Les conditions  $b_1\vi x_1= 0$ et $b_\ell\vu x_\ell=1$ sont immédiates. Soit maintenant $1\leq i<\ell$, on doit vérifier que $b_{i+1}\vi x_{i+1}\leq b_i\vu x_i$. D'une part $b_{i+1}\vi x_{i+1}\leq a_{i+1}\vi x_{i+1}\leq a_i\vu x_i$ car $b_{i+1}\leq a_{i+1}$. D'autre part 
$b_{i+1}\vi x_{i+1}\leq b_{i+1} \leq b_{i+1}\vu x_i$. Donc par distributivité  $b_{i+1}\vi x_{i+1}\leq (a_i\vi b_{i+1})\vu x_i=b_i\vu x_i$.
\end{proof}

La définition \ref{defDiTr} dit
que l'on a $\dim(T)\leq \ell-1$ \ssi
toute \proel $\overline{(x_1,\dots,x_{\ell})}$ s'effondre.
Le lemme suivant dit que l'on peut se limiter aux
 \proels $\overline{(x_1,\dots,x_{\ell})}$ telles
que  $x_1\leq x_2\leq\cdots\leq x_{\ell}$.

\begin{lemma}\label{lemEspanol}
Si pour toute suite $x_1\leq x_2\leq\cdots\leq x_{\ell}$
la  \proel $\overline{(x_1,\dots,x_{\ell})}$ s'effondre
alors $\dim(T)\leq \ell-1$.
\end{lemma}

\begin{proof}
Soit $y_1,\dots,y_{\ell}$ une suite arbitraire. On considère alors
$x_1=y_1,$ $x_2=y_1\vu y_2$, \dots, $x_\ell=y_1\vu\dots\vu y_\ell$.
Soit $(a_1,\dots,a_{\ell})$ une suite complémentaire de $(x_1,\dots,x_{\ell})$.
On prend $b_1=a_1$  et $b_{i+1} = a_{i+1}\vu  x_i$  \hbox{pour $1\leq i<\ell$} (de sorte que $x_{i}\vu  a_i=y_i\vu b_i$ pour $1\leq i\leq \ell$). \\
On a tout d'abord $0 = x_1\vi a_1=y_1\vi b_1$ et $1=x_{\ell}\vu a_{\ell}=y_\ell\vu b_\ell$. Voyons maintenant une inégalité intermédiaire pour $1\leq i<\ell$.
On a $x_{i+1}\vi a_{i+1}\leq  x_i\vu a_i$, et donc
\[ 
y_{i+1} \vi a_{i+1} \leq  x_{i+1} \vi a_{i+1}  \leq  x_i\vu a_i= y_i\vu b_i
\]
On a alors
\[ y_{i+1}\vi b_{i+1}=  y_{i+1}\vi (a_{i+1}\vu  x_i)=
(y_{i+1}\vi a_{i+1})\vu (y_{i+1}\vi x_i)\leq (y_{i+1}\vi a_{i+1})\vu x_i
\]
Comme les deux derniers termes après $\leq $ sont majorés par $x_i\vu a_i= y_i\vu b_i$, on obtient bien l'inégalité $y_{i+1}\vi b_{i+1}\leq y_i\vu b_i$.
\end{proof}

Nous décrivons maintenant la caractérisation présentée
dans \cite{esp}. 

Soit $B=\mathrm{Bo}(T)$ l'algèbre de Boole engendrée
par $T$.
Rappelons
\cite{Macneille1} que tout élément de $B$ peut être
décrit comme une union finie
$\Vu (a_i-b_i)$ de différences formelles,
avec $a_{i+1}\leq b_i\leq a_i\in T.$
En général, on ne
peut rien dire sur la longueur minimale d'une telle suite.

 Espa\~nol \cite{esp} donne la caractérisation suivante de
$\dim(T)\leq \ell$  (au sens de Joyal):

\smallskip \noindent 
\emph{Soit $k\geq 0$:
\begin{itemize}
\item on a $\dim(T)\leq 2k+1$ \ssi tout élément de $B$ peut s'écrire
$\Vu _{1\leq i\leq k}(a_i-b_i)$,
\item on a  $\dim(T)\leq 2k$
\ssi tout  élément de $B$ peut s'écrire
$a\vu \Vu _{1\leq i\leq k-1}(a_i-b_i)$.
\end{itemize}
 }

\smallskip Pour simplifier les notations, si $a_1,\dots,a_{\ell}\in T$
on écrira
$$(a_1-a_2)\vu (a_3-a_4)\vu\cdots$$
pour
$$ (a_1-a_2)\vu (a_3-a_4)\vu\cdots\vu (a_{2k-1}-a_{2k})$$
si $\ell=2k$ et pour
$$ (a_1-a_2)\vu (a_3-a_4)\vu\cdots\vu (a_{2k-1}-a_{2k})\vu a_{2k+1}$$
si $\ell=2k+1.$

Avec cette notation, la condition d'Espa\~nol devient:
on a $\dim(T)\leq \ell$ \ssi tout élément de $B$ peut s'écrire
$$(a_1-a_2)\vu \cdots$$
pour une suite $a_1,\dots,a_{\ell}$ dans $T$.
Pour montrer l'équivalence entre ces deux caractérisations, on
utilisera le fait suivant.

\begin{lemma}\label{lattice2}
Si $x_1\geq \dots\geq x_{\ell}$ et $a_1\geq \dots\geq a_{\ell}$
vérifient
$$(1-x_1)\vu (x_2-x_3)\vu\cdots = (a_1-a_2)\vu (a_3-a_4)\vu\cdots$$
on a
$$1 = x_1\vu a_1,\;0 = x_{\ell}\vi a_{\ell},~~x_i\vi a_i\leq x_{i+1}\vu
a_{i+1} \; (1\leq i<\ell)$$
\end{lemma}

\begin{proof}
Simple vérification.
\end{proof}

\begin{theorem}
On a $\dim(T)\leq \ell-1$, au sens de la définition \ref{defDiTr},
\ssi tout élément de $B=\mathrm{Bo}(T)$ peut s'écrire
$$(a_1-a_2)\vu \cdots$$
pour une suite $a_1\geq \cdots\geq a_{\ell}$ dans $T$.
\end{theorem}

\begin{proof}
On suppose $\dim(T)\leq \ell-1$ au sens de la définition \ref{defDiTr}.
En général si $x_1\geq x_2\geq\dots$ on a
$$ (x_1-x_2)\vu (x_3-x_4)\vu\dots = x_1 - ((x_2-x_3)\vu \dots)$$
et il suffit de montrer que pour toute suite $x_0\geq x_1,\dots\geq
x_{\ell}$
on peut trouver $a_1\geq \dots\geq a_{\ell}$ telle que
$$ (x_0-x_1)\vu (x_2-x_3)\vu\dots = (a_1-a_2)\vu\dots$$
Comme $(x_0-x_1)\vu (x_2-x_3)\vu\dots = x_0\vi ((1-x_1)\vu (x_2-
x_3)\vu\dots)$
on se ramène au cas $x_0=1$, qui résulte du lemme \ref{lattice3}.

 La réciproque résulte directement du lemme \ref{lattice2} qui donne l'hypothèse du lemme \ref{lemEspanol}.
\end{proof}

Signalons enfin l'article plus récent \cite[Español]{Espa08} qui présente une autre variante pour la définition \cov de la dimension de Krull. Cette variante
est étudiée en exercice\footnote{Une erreur est rectifiée dans la version sur arXiv.} XIII-17 de \cite{ACMC}.

\section{Treillis de Zariski et de Krull dans un anneau commutatif}
\label{secZariKrull}
\subsection{Treillis de Zariski}
\label{subsubsecZar}
Dans un anneau commutatif $A$, le {\em treillis de Zariski} $\Zar(A)$ a 
pour
éléments les radicaux d'\itfs (la relation d'ordre est l'inclusion).
Il est bien défini en tant que treillis.
Autrement dit $\sqrt{I_1}=\sqrt{J_1}$ et $\sqrt{I_2}=\sqrt{J_2}$
impliquent
$\sqrt{I_1I_2}= \sqrt{J_1J_2} $
(ceci définit $\sqrt{I_1}\vi\sqrt{I_2}$) et
$\sqrt{I_1+I_2}=\sqrt{J_1+J_2}$
(ceci définit $\sqrt{I_1}\vu\sqrt{I_2}$). Le treillis de Zariski de
$A$
est toujours un \trdi, mais en général l'égalité n'est pas
testable. Néanmoins une inclusion $\sqrt{I_1}\subseteq \sqrt{I_2}$ peut
être certifiée de manière finie si l'anneau $A$ est discret.
Ce treillis contient toutes les informations nécessaires au
développement du point de vue \cof concernant la théorie abstraite
et non \cov du spectre de Zariski.

Tout morphisme $\varphi:A\to B$ d'anneaux commutatifs induit un morphisme naturel $\Zar(\varphi)=\Zar(A)\to\Zar(B)$. 

Nous notons $\wi{a}$ pour $\sqrt{\gen{a}}$. Pour une partie $S$ de $A$
nous notons  $\wi{S}$ la partie de $\Zar(A)$ formée des  $\wi{s}$ pour
$s\in S$. On a $\wi{a_1}\vu\cdots\vu\wi{a_m}=\sqrt{\gen{a_1,\ldots,a_m}}$ 
et
$\wi{a_1}\vi\cdots\vi\wi{a_m}=\wi{a_1\cdots a_m}$ \\
Soient $U$ et  $J$  deux familles finies dans $A$, on a
$$\Vi\wi{U}\vdash_{\Zar(A)} \Vu\wi{J}
\quad\Longleftrightarrow \quad
\prod\nolimits_{u\in U} u  \in \sqrt{\gen{J}}
\quad\Longleftrightarrow \quad
\cM(U)\cap \gen{J}\neq \emptyset
$$
où $\cM(U)$ est le \mo multiplicatif engendré par $U$, \cad encore
$$ \hbox{ le \proi}(J,U){\rm \;s'effondre \;dans  \;  } A
\,\Longleftrightarrow \,
(\wi{J},\wi{U})\hbox{ s'effondre dans } \Zar(A)
$$
Cela suffit à décrire ce treillis. Plus précisément on a:
\begin{proposition}
\label{propZar} Le treillis $\Zar(A)$ d'un anneau commutatif $A$  est
(à \iso près) le treillis engendré par  $(A,\vda)$  où $\vda$
est la plus petite \entrel vérifiant
$$\begin{array}{rclcrclcrcl}
  0_A   & \vda  &      &\qquad  & x,\; y & \vda  &  x y   \\
        & \vda  &  1_A &\qquad &  xy    & \vda  &   x   &\qquad &
  x+y   & \vda  & x,\; y  \\
\end{array}$$
\end{proposition}
\begin{proof}
Il est clair que la relation $U\vda J$ définie par \gui{$\cM(U)$
intersecte
$\gen{J}$} vérifie ces axiomes. Il est clair aussi que la relation
implicative engendrée par ces axiomes contient cette relation.
Montrons donc que cette relation est une \entrel.
Seule la règle de coupure n'est pas directe.
Supposons que $\cM(U,a)$ intersecte
$\gen{J}$ et que $\cM(U)$ intersecte $\gen{J,a}$. On peut alors trouver
$m_1,m_2\in \cM(U)$ et $k,x$ tels que
$a^k m_1 \in \gen{J},~m_2+ax\in \gen{J}$.
En éliminant $a$ on obtient que $\cM(U)$ intersecte $\gen{J}.$
\end{proof}
Notez que l'application canonique de $A$ dans $\Zar(A)$ est
$a\mapsto \wi{a}$ et que l'on a $\wi{a}=\wi{b}$ \ssi $a$ divise une
puissance de $b$ et $b$ divise une puissance de $a$.

\Subsubsection{Idéaux et filtres de $A$ et de $\Zar\,A$}
\begin{fact}
\label{factQuotients} \emph{(quotients par les idéaux)}
Soit $I$ un idéal de $A$. 
Le morphisme naturel $A\to A/I$ induit un morphisme surjectif $\Zar(A)\to\Zar(A/I)$. Ce dernier donne par factorisation un isomorphisme $\Zar(A/I)\simeq
\Zar(A)/(\cI =0)$, où ${\cI} = \sotq{ J \in \Zar(A)}{J \subseteq \sqrt I}.$
\end{fact}

\begin{proposition}
\label{propZar2}
On a une bijection naturelle entre les
idéaux du treillis $\Zar(A)$ et les idéaux radicaux de l'anneau $A$. La bijection fonctionne comme suit.
\begin{itemize}
\item Si $I$ est un idéal radical de $A$ on lui  associe l'idéal
${\cI} = \{ J \in \Zar(A)~|~ J \subseteq I\}$.
\item  Si $\cal{I}$ est un idéal
de $\Zar(A)$ on lui associe l'idéal
$I=\bigcup\nolimits _{J\in\cal{I}}J=\{ x\in A~|~\wi{x}\in \cal{I}\}.$
\item Dans cette bijection les idéaux premiers correspondent aux idéaux
premiers. 
\end{itemize}
\end{proposition}
\begin{proof}
Nous expliquons seulement la dernière affirmation.\\
Si $I$ est un idéal premier de $A$, si $J,J'\in \Zar(A)$ et
$J\vi J'\in \cal{I}$, soient $a_1,\dots,a_n\in A$ des \gui{générateurs}
de $J$ (\cad $J=\sqrt{\gen{a_1,\dots,a_n}}$) et $b_1,\dots,b_m\in A$ des
générateurs de~$J'.$ On a alors
$a_ib_j\in I$ d'où $a_i\in I$ ou $b_j\in I$ pour tout $i,j.$ Il en
résulte que l'on a $a_i\in I$ pour tout $i$ ou
$b_j\in I$ pour tout $j$. Donc $J\in \cal{I}$ ou $J'\in \cal{I}$ et
$\cal{I}$ est un idéal premier de $\Zar(A).$
\\
Réciproquement si $\cal{I}$ est un idéal premier de $\Zar(A)$
et si on a $\wi{xy}\in \cal{I}$ alors on a
$\wi{x}\vi \wi{y}\in \cal{I}$ et donc $\wi{x}\in \cal{I}$
ou $\wi{y}\in \cal{I}$ ce qui montre que $\{ x\in A~|~\wi{x}\in \cal{I}\}$
est un idéal premier de $A.$
\end{proof}

Un \emph{filtre} dans un anneau commutatif est un \mo $S$ qui  vérifie
$xy\in S\Rightarrow x\in S$ (on dit aussi que le \mo $S$ est saturé). 

Un \emph{filtre premier} est un filtre 
qui vérifie $x+y\in S\Rightarrow x\in S\;\mathrm{ou}\;y\in S$ (c'est 
le complémentaire d'un \idep). 

Un $x\in A$ engendre un \emph{filtre principal} $\wh x=\sotq{y}{\exists n\in\N\;x^n\in\gen{y}}=\sotq{y\in A}{x\in \wi y}$.

Pour $x\in A$ le filtre $\uar{\wi x}$ de $\Zar(A)$ satisfait l'équivalence: 
\fbox{$\wi y\in  \uar{\wi x} \Leftrightarrow y\in \wh x$}.

\begin{fact}
\label{factLocalises} \emph{(localisés)}
Soit $S$ un filtre de $A$. 
Alors le morphisme naturel $A\to A_S$ induit un morphisme surjectif $\Zar(A)\to\Zar(A_S)$. Ce dernier donne par factorisation un isomorphisme $\Zar(A_S)\simeq
\Zar(A)/(\uar S =1)$.
\end{fact}

\begin{proposition}
\label{factFZ} Pour un filtre $S$ de $A$, notons $\uar S=\bigcup_{s\in S} \uar \wi s$. C'est un filtre de $\Zar(A)$.
\begin{itemize}
\item 
L'application $S\mapsto  \uar S$ établit une correspondance
\emph{injective} croissante des filtres de $A$ vers les filtres de 
$\Zar(A)$. 
\item  Cette application préserve les sups finis
(le sup de $S_1$ et $S_2$  est engendré par les 
$s_1s_2$ où
$s_i\in S_i$).
\item Elle  se restreint en une bijection des filtres
premiers de $A$ sur ceux de $\Zar(A)$. 
\end{itemize}
\end{proposition}

Notez cependant que le filtre principal de $\Zar(A)$ engendré par
$\wi{a_1}\vu\cdots \vu \wi{a_n}$ (\cad l'intersection des filtres 
$\uar\wi
{a_i}$), ne correspond en général à aucun filtre de $A$.

\subsection{Treillis de Krull d'ordre $\ell$}
\label{subsubsecKrulac}

\begin{definition}{\rm 
\label{defKruA}
On définit $\Kru_\ell(A):=\Kr_\ell(\Zar(A))$. On l'appelle
le {\em treillis de Krull d'ordre $\ell$ de l'anneau $A$}
}
\end{definition}

\begin{theorem}
\label{corCollaps} Soit $\cC=\big((J_0,U_0),\ldots,(J_\ell,U_\ell)\big)$  une
\proc dans un anneau commutatif $A$.
Elle s'effondre \ssi la \proc
$\big((\wi{J_0},\wi{U_0}),\ldots,
(\wi{J_\ell},\wi{U_\ell})\big)$  s'effondre dans
$\Zar(A)$. Par exemple si $\cC$ est finie, \propeq
\begin{enumerate}
\item  il existe $j_i\in \gen{J_i}$, $u_i\in\cM(U_i)$,
$(i=0,\ldots,\ell)$,
vérifiant l'égalité
$$u_0\cdot(u_1\cdot(\cdots(u_\ell+j_\ell)+\cdots)+j_1)+j_0=0,
$$
\item  il existe $x_1,\ldots,x_\ell\in \Zar(A)$ avec les relations
suivantes
dans $\Zar(A)$:
$$\begin{array}{rcl}
 x_1,\; \widetilde{U_0}& \vda  &  \widetilde{J_0}
\\
 x_2,\; \widetilde{U_1}& \vda  &  \widetilde{J_1} ,\;  x_1
\\
\vdots\qquad & \vdots  & \qquad\vdots
\\
x_\ell,\; \widetilde{U_{\ell-1}}&\vda& \widetilde{J_{\ell-1}},\;x_{\ell-1}
\\
\widetilde{U_\ell}& \vda & \widetilde{J_\ell} ,\; x_\ell
\end{array}$$
\item  même chose mais avec $x_1,\ldots,x_\ell\in \wi{A}$.
\end{enumerate}
\end{theorem}
\begin{proof}
Il est clair que \emph{1} entraine \emph{3}: on prend simplement
$$v_\ell = u_\ell+j_\ell,~v_{\ell-1} = v_\ell u_{\ell-1} +
j_{\ell-1},\dots,~v_0 = v_1u_0 + j_0
~~~~{\rm  et}~~~ x_i= \wi{v_i}
$$
et que \emph{3} entraine \emph{2}.
Le fait que \emph{2} entraine \emph{1} peut se voir en reformulant
\emph{2} de la manière suivante.
On considère la \proc $\cC^1=\big((K_0,V_0),\ldots,(K_\ell,V_\ell)\big)$ obtenue
en saturant la \proc $\cC$.
On définit  $\ell+1$ idéaux radicaux
$I_0,\dots,I_{\ell}$ de $A$
\begin{itemize}
\item $I_0 = \{x\in A~|~\cM(x,U_0)\cap \gen{J_0} \neq \emptyset \}$
\item $I_1 = \{x\in A~|~\cM(x,U_1)\cap (\gen{J_1} + I_0) \neq \emptyset\}$
\item ~~$\vdots$
\item $I_{\ell-1} = \{x\in A~|~\cM(x,U_{\ell-1})\cap (\gen{J_{\ell-1}} +
I_{\ell-2}) \neq \emptyset \}$
\item $I_\ell=\gen{J_\ell}+I_{\ell-1}$
\end{itemize}
Il est clair que $I_i\subseteq K_i$ ($i=0,\ldots,\ell$).
Dans la corrrespondance donnée en \ref{propZar2} ces idéaux
correspondent aux idéaux suivants de $\Zar(A)$
  \begin{itemize}
\item ${\cI}_0 =
  \{u\in \Zar(A)~|~ u,\;  \widetilde{U_0} \vda  \widetilde{J_0}\}$
\item ${\cI}_1 =
  \{u\in \Zar(A)~|~ (\exists v\in {\cI}_0)~
                    u,\;  \widetilde{U_1} \vda \widetilde{J_1},\;
v\}$
\item ~~$\vdots$
\item ${\cI}_{\ell-1} =
  \{u\in \Zar(A)~|~ (\exists v\in {\cI}_0)~
                     u,\;  \widetilde{U_{\ell-1}} \vda
\widetilde{J_{\ell-1}},\; v\}$
\end{itemize}
La condition  \emph{2}  signifie alors que l'on a
$\;\widetilde{U_{\ell}}\vda  \widetilde{J_{\ell}},\; v\;$ pour un
$\;v\in {\cI}_{\ell-1}$.
Autrement dit, $\cM(U_{\ell})$ intersecte
$I_{\ell}$, ou $I_{\ell}\subseteq K_{\ell}$. 
Donc $\cC^1$ s'effondre, donc $\cC$ s'effondre.

Donnons une autre démonstration, directe, que  \emph{2} entraine \emph{3}.
On réécrit les \entrels de \emph{2} comme suit.
Chaque $\wi{U_i}$ peut être remplacé par un $\wi{u_i}$ avec $u_i\in 
A$,
chaque~$\wi{J_i}$ peut être remplacé par un radical d'\itf $I_i$  de
$A$, et on note $L_i$ à la place de~$x_i$ pour se rappeler qu'il s'agit 
du
radical d'un \itf. On obtient:
$$\begin{array}{rcl}
 L_1,\; \widetilde{u_0}& \vda  &  I_0  \\
 L_2,\; \widetilde{u_1}& \vda  &  I_1 ,\;  L_1  \\
 L_3,\; \widetilde{u_2}& \vda  &  I_2 ,\; L_2  \\
\widetilde{u_3}& \vda & I_3 ,\; L_3
\end{array}$$
La dernière ligne signifie que $\cM(u_3)$ coupe $I_3+ L_3$
ou encore  $I_3+ \gen{y_3}$ pour un élément $y_3$ de
$L_3$, et donc aussi que l'on a
$ \wi{u_3} \vda  I_3 ,\; \wi{y_3}$.
Comme $\wi{y_3}\leq L_3$ dans $\Zar(A)$ l'avant-dernière ligne
implique  $ \;\wi{y_3},\; \wi{u_2} \vda  I_2 ,\; L_2$.
On a donc remplacé dans ses deux occurences  $L_3$ par~$\wi{y_3}$
On raisonne comme précédemment et on voit qu'on peut remplacer
les deux occurences de~$L_2$ par un $\wi{y_2}$ convenable, puis
les deux occurences de $L_1$ par un $\wi{y_1}$ convenable.
On a bien obtenu~\emph{3}.
\end{proof}

\begin{corollary}
\label{corCompDim}
La dimension de Krull d'un anneau commutatif $A$ est $\leq \ell$ (définition \ref{def.dimKrull}) \ssi
la dimension de Krull du treillis de Zariski $\Zar(A)$ est $\leq \ell$
(définition~\ref{defDiTr2}).
\end{corollary}

\begin{proof}
Appliquer le théorème précédent et le lemme \ref{lemDimGen}.
\end{proof}

Notons que nous avons obtenu ce résultat de manière purement \cov, sans avoir à utiliser la définition classique usuelle en termes de chaînes d'\ideps.

\section{Dimension de Krull relative} 
\label{secKrullRel}
\subsection{Généralités sur la dimension de Krull relative } 

Nous développons ici un analogue \cof pour les chaînes croissantes 
formées par des idéaux premiers qui coupent tous un sous-anneau donné selon un même idéal premier. Ce paragraphe est en fait suffisamment général 
pour fonctionner dans le cadre de la section~\ref{secKrA} (anneau commutatif) ou dans celui d'un \trdi arbitraire. Il n'y a pas de vrais calculs, seulement un peu de combinatoire.

Dans la suite \gui{Soit $A\subseteq B$} est une abréviation pour \gui{Soit $A\subseteq B$ des anneaux commutatifs ou des \trdis}.
\begin{definition}{\rm  
\label{defColRel} 
Soit $A\subseteq B$ et 
$\cC=\big((J_0,U_0),\ldots,(J_\ell,U_\ell)\big)$ une \proc dans $B$. 
\begin{enumerate}
\item On dit que {\em la \proc $\cC$ s'effondre au dessus de $A$} s'il 
existe un entier $k\geq 0$ et des éléments $a_1,\ldots,a_k$ de $A$ tels que  pour tout couple de parties 
complémentaires $(H,H')$ de $\{1,\ldots,k\}$, 
on ait le collapsus dans $B$ de la \proc
\begin{equation} \label {eqdefColRel}
(\{(a_h)_{h\in H}\}\cup J_0,U_0),(J_1,U_1)\ldots,(J_\ell,U_\ell\cup 
\{(a_h)_{h\in H'}\})
\end{equation}
\item On dit que {\em la dimension de Krull (relative) de l'extension 
$B/A$ est $\leq  \ell-1$} si toute \proel  
$\big((0;x_1),(x_1;x_2),\ldots,(x_\ell;1)\big)$ s'effondre au dessus de $A$.
\end{enumerate}
}
\end{definition}

Il est clair qu'une \proc de $B$ qui raffine un \proc qui s'effondre au dessus de $A$ s'effondre elle-même au dessus de $A$.
Le lemme qui suit justifie l'emploi du lemme de Zorn en \clama.
\begin{lemma} \label{lemdefColRel}
Toute \proc de $B$ qui s'effondre au dessus de $A$ raffine une \proc finie qui s'effondre au dessus de $A$.
\end{lemma}

Pour ce qui concerne la dimension de Krull relative, on obtient
le critère suivant, qui reprend le point \emph{2} de la \dfn \ref{defColRel}. 
\begin{lemma} \label{lemKrullRel}
La dimension de Krull (relative) de l'extension 
$B/A$ est $\leq  \ell$ \ssi pour toute liste $(x_0,\dots,x_\ell)$ dans $B$
il existe un entier $k\geq 0$ et des éléments $a_1,\ldots,a_k\in A$ tels que  pour tout couple de parties 
complémentaires $(H,H')$ de $\{1,\ldots,k\}$, il existe $ y_0,\dots,y_\ell\in B$ tels 
que
\begin{equation} \label {eqdefDiTrRel}
\begin{array}{rclll}
\Vi_{j\in H'} a_j & \vda  &  y_\ell,\;x_\ell    \\
y_\ell,\;x_\ell& \vda  &y_{\ell-1},\;x_{\ell-1}      \\
\vdots\qquad & \vdots  & \qquad  \vdots    \\
y_1,\;x_1& \vda  &  y_0,\;x_0    \\
y_0,\;x_0& \vda  &  \Vu_{j\in H} a_j    \\
\end{array}
\end{equation}
\end{lemma}
Par exemple pour la dimension relative $\leq 2$ cela correspond au dessin suivant dans $B$ avec $u=\Vi_{j\in H'} a_j$ et $i=\Vu_{j\in H} a_j$.
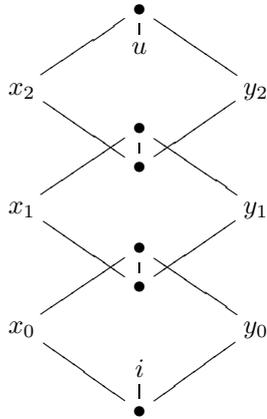
$$\SCOR{x_0}{x_1}{x_2}{y_0}{y_1}{y_2}{u}{i}$$
\begin{remark} \label{remdefColRel} 
{\rm   On peut considérer le cas d'une extension d'anneaux ou de \trdis plus générale:  
on a un \homo $A\rightarrow B$  non nécessairement injectif. 
On peut alors adapter la définition précédente en remplaçant~$A$  par son image dans $B$. 
}\end{remark}

Voici un théorème de collapsus simultané relatif.
\begin{theorem} 
\label{thColSimRel} {\em (Collapsus simultané relatif pour les \procs)} 
Soit $A\subseteq B$  et $\cC$ une \prolo $\ell$ 
dans $B$. 
\begin{itemize}
\item [$(1)$] Soit $x\in B$ et $i\in\left\{0,\ldots ,\ell\right\}$.
Supposons que les \procs 
$\cC\,\& \left\{x\in \cC(i) \right\} $
\hbox{et $\cC\,\& \left\{x\notin \cC(i) \right\} $}
s'effondrent toutes les deux au dessus de $A$, alors $\cC$ s'effondre 
également au dessus de $A$.
\item [$(2)$] Soit $x\in A$. Supposons que les \procs 
$\cC\,\& \left\{x\in \cC(0) \right\} $
et 
$\cC\,\& \left\{x\notin \cC(\ell) \right\} $
s'effondrent toutes les deux au dessus de $A$, alors $\cC$ s'effondre 
également au dessus de~$A$.
\end{itemize}
\end{theorem}
C'est une conséquence facile des théorèmes (non relatifs) 
\ref{ThColSimKrA} pour les anneaux commutatifs et~\ref{thColSimT2} pour les \trdis, que nous laissons au lecteur.
De là, on déduit (en \clama) une caractérisation des \procs qui 
s'effondrent relativement.
 
\begin{theorem} 
\label{th.nstformelRel} {\em (\nst formel pour les chaînes d'idéaux 
premiers dans une extension)} 
 \Tcgi \\
Soit $A\subseteq 
B$  et soit $\cC=\big((J_0,U_0),\ldots,(J_\ell,U_\ell)\big)$ une 
\proc dans $B$. \Propeq
\begin{itemize}
\item [$(a)$] Il existe un idéal premier détachable $P$ de $A$ et 
$\ell+1$ idéaux premiers détachables 
$P_0\subseteq \cdots\subseteq P_\ell$ de $B$ tels que
$J_i\subseteq P_i$, $U_i\cap P_i=\emptyset $ et $P_i\cap A=P$  
$(i=0,\ldots,\ell)$. 
\item [$(b)$] La \proc $\cC$ ne s'effondre pas au dessus de $A$.
\end{itemize}
\end{theorem}
\begin{proof}
L'implication $(a)\Rightarrow (b)$ est presque immédiate. On dispose d'une suite croissantes $P_0,\dots,P_\ell$ d'\ideps de $B$ qui raffine $\cC$ (chaque $\wh{P_i}$ raffine $(J_i,U_i)$), avec un \idep $P$ de~$A$ tel que l'on a $P_i\cap A=P$ pour $i=0,\dots,\ell$.
Le fait que $\cC$ ne collapse pas au dessus de~$A$ signifie que pour tous $a_1,\dots,a_k\in A$ il existe une partie $H$ de $\so{1,\dots,k}$ telle que
la suite~(\ref{eqdefColRel}) ne collapse pas dans $B$. Il suffit en effet de prendre pour $H$ les indices $i$ tels \hbox{que $a_i\in P$}, car alors $(\wh{P_0},\dots,\wh{P_\ell})$
raffine la suite~(\ref{eqdefColRel}). 

\noindent Pour $(b)\Rightarrow (a)$ nous 
faisons la démonstration qui s'appuie sur le \pte et le lemme de 
Zorn (plus immédiate que celle qui utilise le \tcg). On considère une \proc 
$\cC^1=\big((P_0,S_0),\ldots,(P_\ell,S_\ell)\big)$ 
maximale (pour la relation de raffinement) parmi les \procs qui raffinent 
$\cC$ et qui ne s'effondrent pas au dessus de~$A$. Vu le collapsus 
simultané relatif, la même démonstration que pour les théorèmes 
\ref{th.nstformel} \hbox{et \ref{th.nstformelTreil}} montre que c'est une chaîne croissante d'idéaux 
premiers (avec leurs compléments). 
Il reste à voir que les $P_i\cap A$  sont tous égaux, ce qui revient 
à dire $S_0\cap P_\ell\cap A=\emptyset$. 
Si ce n'était pas le cas soit $x\in S_0\cap P_\ell\cap A$.   
Alors $\big((P_0,x;S_0),\ldots,(P_\ell,S_\ell)\big)$ et 
$\big((P_0,S_0),\ldots,(P_\ell;S_\ell,x)\big)$ s'effondrent dans $B$, 
donc $\cC^1$ s'effondre au dessus de $A$  (avec la suite à un seul élément $x$). 
C'est absurde.
\end{proof}

\subsection{L'inégalité de base} 

On a le résultat \cof suivant.
\begin{theorem} 
\label{thDim1} 
Soit $A\subseteq B$ des anneaux commutatifs ou des \trdis. 
\begin{itemize}
\item [$(1)$] Supposons que la dimension de Krull de $A$ est $\leq m$ 
et que la dimension de Krull de l'extension $B/A$ est $\leq n$, 
alors  la dimension de Krull de $B$ est $\leq (m+1)(n+1)-1$.
\item [$(2)$] Supposons que $A$ et $B$ sont munis du prédicat $\neq0$
défini comme la négation de $=0$.
 Supposons que la dimension de Krull de l'extension $B/A$ 
est $\leq n$ et qu'on a un test pour le collapsus des \proels dans $A$. Si 
on donne une suite pseudo régulière de longueur $(m+1)(n+1)$  dans 
$B$, on peut construire une suite pseudo régulière de longueur $m+1$  
dans~$A$.
\end{itemize}
\end{theorem}
\begin{proof}
En \clama la démonstration est immédiate. On considère une chaîne 
strictement croissante dans $B$ avec $(m+1)(n+1)+1$ termes. Comme $n+2$ 
termes consécutifs ne peuvent avoir la même intersection avec $A$ cela 
fournit une chaîne strictement croissante de $m+2$ idéaux de~$A$, 
ce qui est absurde. En \coma on peut mimer cette démonstration et on obtient le 
résultat sous forme \cov, ce qui nous donne une véritable information 
de nature algorithmique. Voici ce que cela donne.\\
Nous prouvons d'abord le point (1).\\
Tout d'abord nous traitons le cas où $A$ est de dimension zéro avec la 
dimension relative égale~à~$n$. Nous voulons montrer que toute 
\proelo$n$ dans $B$ s'effondre. Soit~\hbox{$\cC=\big((0,x_1),\ldots,(x_n,1)\big)$} une telle 
\proel. Nous savons qu'elle s'effondre au dessus de $A$. 
Soit $F=\{a_1,\ldots,a_k\}$ la partie finie de $A$ correspondante. Si $F$  
est vide c'est bon. Si~$F$ n'est pas vide, on va montrer qu'on peut 
enlever un élément à $F$ et on aura gagné par induction. 
Soit donc $F=H\cup \{a_k\}$. Soit maintenant $G$  et $G'$  deux parties  
complémentaires de $H$. On a les collapsus (tout court) des \procs 
suivantes, la première par collapsus au dessus de $A$, la deuxième 
parce qu'elle contient $\big((0;a_k),(a_k;1)\big)$ et que $A$  est de dimension 
$0$:
$$\begin{array}{rcl} 
\big((G;x_1),\ldots,(x_n;G',a_k)\big) &\; \; {\rm donc\; aussi}\; \; &  
\big((G;x_1,a_k),\ldots,(x_n;G',a_k)\big)\\[.3em] 
\big((0;x_1,a_k),\ldots,(a_k,x_n;1)\big) &\; \; {\rm donc\; aussi}\; \; &  
\big((G;x_1,a_k),\ldots,(a_k,x_n;G')\big) 
\end{array}$$
Par collapsus simultané, les deux colonnes de droite donnent le 
collapsus de
$$\big((G;x_1,a_k),\ldots,(x_n;G')\big).$$
Comme on a aussi, 
puisque $\left\{a_k \right\}\cup G\cup G'=F $,
le collapsus (tout court) de la \proc $$\big((a_k,G;x_1),\ldots,(x_n;G')\big),$$ un
collapsus simultané donne celui de la \proc $\big((G;x_1),\ldots,(x_n;G')\big)$. 
Et nous avons gagné.
\\
Nous passons au cas général. Nous traitons l'exemple où $m=2,\; 
n=3$, qui, nous l'espérons, est suffisamment éclairant. Nous devons 
donc considérer une suite de $3\times 4=12$ éléments $x_i$ de $B$ et 
montrer que  la \proel suivante s'effondre:
$$\big((0;x_1),(x_1;x_2),\ldots, (x_{11};x_{12}),(x_{12};1)\big) 
$$
Par hypothèse les \procs
$$\begin{array}{cl} 
\big((0;x_1),\ldots,(x_4;x_5)\big),  \; 
\big((x_4;x_5), \ldots,(x_8;x_9)\big) \; \hbox{ et }\;
\big((x_8;x_9),\ldots, (x_{12};1)\big)
\end{array}$$
s'effondrent au dessus de $A$, ce qui nous fournit trois listes finies 
$F_1$, $F_2$, $F_3$  d'éléments de $A$ et les collapsus correspondants 
dans $B$.
Nous pouvons montrer que la \proc $\cC$ s'effondre en utilisant un processus 
de recollement analogue à celui utilisé lorsque $m=0$.
Nous commençons par remarquer que toute \proc du type
$$ \big((0;x_1,b_1),\ldots, (b_1,x_4;x_5,b_2),
 \ldots, (b_2,x_8;x_9,b_3),\ldots, (b_3,x_{12};1)\big)
$$
avec les $b_i\in A$ s'effondre parce qu'elle \gui{contient} la \proc
$$ \big((0;b_1),(b_1;b_2), (b_2;b_3), (b_3;1)\big)
$$
On va établir que toute \proc du type 
$$ \big((0;x_1),\ldots, (x_4;x_5,b_2),
 \ldots, (b_2,x_8;x_9,b_3),\ldots, (b_3,x_{12};1)\big)
$$
avec $b_2,\,b_3\in A$ s'effondre elle aussi. Après on pourra recommencer 
(on supprimera $b_2$ puis $b_3$). On voit qu'on peut reproduire à 
l'identique le cas déjà traité ($A$ zéro dimensionnel): en 
ajoutant partout la queue $;x_5,b_2),(x_5,x_6) \ldots, 
(b_2,x_8;x_9,b_3),\ldots, (b_3,x_{12};1))$.\\
Nous prouvons ensuite le point (2).\\
Il faut reprendre la démonstration du point (1) (en faisant comme un semblant 
raisonnement par l'absurde) et regarder à quel endroit elle ne fonctionne  
plus. Elle est basée sur des collapsus de \proels dans $A$ et des 
collapsus de \proels de $B$ au dessus de $A$.
Les seuls endroits où cela ne fonctionne plus à tout coup, c'est avec les 
collapsus de \proels dans $A$. Or justement, on suppose qu'on a un test 
pour ces collapsus-là. Donc l'un d'entre eux au moins ne fonctionne pas, 
explicitement, et cela fournit la suite pseudo régulière que l'on 
cherche.
\end{proof}

\subsection{Comparaison avec une autre approche \cov dans le cas des anneaux commutatifs} 

Dans \cite{ACMC}, pour les anneaux commutatifs,  la dimension de Krull relative de l'extension $B/A$ (ou du morphisme $A\to B$) est définie comme la dimension de Krull de $A\bul\otimes_A B$, où $A\bul$ est l'anneau zéro-dimensionnel réduit engendré par $A$. La justification intuitive est que les \ideps $P$ de $A$ se retrouvent sous la forme des \idemas de $A\bul$ (deux à deux incomparables) et que lorsqu'on a un \idep $P_0$ de $B$ au dessus de $P$, \hbox{si $P_0\subseteq \dots \subseteq P_\ell$} dans $A\bul\otimes_A B$, alors les $P_i$ sont tous au dessus de $P$ dans $B$. La définition dans \cite{ACMC} peut être considérée comme plus élégante dans la mesure où elle réclame comme seul nouvel outil l'anneau $A\bul$ qui résout un problème universel naturel. 

Néanmoins, l'inégalité de base (théorème \ref{thDim1}) est plus laborieuse à établir dans \cite{ACMC}.

\medskip Nous démontrons maintenant que les deux définitions \covs proposées pour la dimension de Krull relative coïncident.
Il suffit de démontrer qu'une \proc $\cC$ de~$B$ s'effondre au dessus de $A$ selon la définition  \ref{defColRel} \ssi $\cC$ s'effondre dans l'anneau~\hbox{$A\bul\otimes _A B$}.

Pour cela on considère un étage fini de la construction de l'anneau~$A\bul$.
Il est obtenu à partir d'une suite $a_1,\dots,a_k$ dans $A$: on  associe à cette suite l'anneau
$$
A[a_1\bul,\dots,a_k\bul]=\prod\nolimits_{H\subseteq \so{1,\dots,k}} A[1/a_{H'}]/a^{H}
$$
avec $a_{H'}=\prod_{i\in H}a_i$ et $a^{H}=\gen{(a_i)_{i\in H}}$; où $H'$
est le complément de $H$ dans $\so{1,\dots,k}$.

L'anneau~$A\bul$ est la limite inductive de ces anneaux $A[a_1\bul,\dots,a_k\bul]$
et l'on voit que le collapsus d'une \proc $\cC$ de $B$ au dessus de $A$
défini en \ref{defColRel}
revient à dire que $\cC$ s'effondre dans un des anneaux $A[a_1\bul,\dots,a_k\bul]\otimes_A B$. Autrement dit $\cC$ s'effondre dans 
$A\bul\otimes _A B$.

Dans les sous-sections \ref{secextent} et \ref{secextpol} nous traitons des résultats qui sont démontrés dans \cite{ACMC} en utilisant la définition de \cite{ACMC} pour la dimension de Krull relative.

\subsection{Cas des extensions entières d'anneaux commutatifs}\label{secextent} 

Dans la proposition suivante le point (1) est la version \cov du 
\gui{théorème d'incom\-pa\-rabilité} (voir le théorème
13.33 du livre de Sharp \cite{Sha}.)

\begin{proposition} 
\label{propDim2}  Soit $A\subseteq B$ des anneaux commutatifs. 
\begin{itemize}
\item [$(1)$] Si $B$  est entier sur $A$ la dimension de Krull relative de 
l'extension $B/A$  est nulle.
\item [$(2)$] Plus généralement on a le même résultat si tout 
élément de $B$ est un zéro d'un polynôme de $A[X]$ ayant un 
coefficient égal à $1$. Par exemple si $A$ est un 
anneau de Prüfer intègre, cela s'applique à n'importe quel 
sur-anneau de $A$ dans son corps de fractions.
\item [$(3)$] En particulier, en appliquant le théorème 
\ref{thDim1}
si  $\dim(A)\leq n$ alors $\dim(B)\leq n$, ce que l'on abrège en $\dim \,B\leq \dim\,A$.   
\end{itemize}
\end{proposition}
\begin{proof}
On montre $(2)$. On veut montrer que pour tout $x\in B$ la \proc
$\big((0;x),(x;1)\big)$ s'effondre au dessus de $A$. 
La liste finie dans $A$ est celle donnée par les coefficients du 
polynôme qui annule~$x$. Supposons que
$x^k=\sum_{i\neq k,i\leq r}a_ix^i$. 
Soit  $\,G,\; G'\,$ deux parties complémentaires de
$\left\{a_i;i\neq k\right\}$. 
Le collapsus de $\big((G,x),(x,G')\big)$ est donné par une égalité 
$\,x^m(g'+bx)=g\,$ avec $g\in \gen{G}_B$, \hbox{$g'\in \cM(G')$} et $b\in B$. 
En fait nous prenons $g\in G[x]$ et $b\in A[x]$. 
Si $G'$ est vide on prend $m=k,\; g'=1$. 
Sinon soit $h$  le plus petit indice $\ell$ tel que $a_\ell\in G'$. 
Tous les $a_j$ avec $j<h$ sont dans $G$. Si $h<k$ on prend $m=h,\; 
g'=a_h$.  Si $h>k$ on prend $m=k,\; g'=1$. \\
NB: on remarque que la disjonction porte sur  $r$ cas et non pas sur 
$2^r$ cas:
\begin{itemize}
\item $ a_0\in G'$, ou
\item $  a_0\in G,a_1\in G'$, ou
\item $  a_0,a_1\in G,a_2\in G'$, ou
\item $~~~~~~~\cdots$
\end{itemize}
\end{proof}


\subsection{Dimension de Krull relative des anneaux de polynômes} \label{secextpol} 
Nous donnons la version constructive du théorème classique concernant  
la dimension de Krull relative de l'extension
$A[x_1,\ldots, x_n]/A$.

Nous aurons besoin d'un lemme d'algèbre linéaire élémentaire.

\begin{lemma} 
\label{lemALE} Soit $V_1,\ldots,V_{n+1}$ des vecteurs de 
$A^{n}$. 
\begin{itemize}
\item  Si $A$ est un corps discret, il existe un indice 
$k\in \left\{1,\ldots,n+1\right\}$ tel que
$V_k$ est une \coli des vecteurs qui suivent (si $k=n+1$ cela signifie 
que $V_{n+1}=0$).
\item  Si $A$ est un anneau commutatif, notons $V$ la matrice dont les 
vecteurs colonnes sont les $V_i$. Soit $\mu_1,\ldots,\mu_\ell$
(avec $\ell=2^n-1$) la liste des mineurs de $V$ extraits sur les 
$n$ ou $n-1$ ou $\ldots$ ou $1$  dernières colonnes, rangée par ordre 
de taille décroissante. On pose $\mu_{\ell+1}=1$ 
(le mineur correspondant à la matrice extraite vide).
Pour chaque $k\in\left\{1,\ldots,\ell+1\right\}$ on pose
$I_k=\gen{(\mu_i)_{i< k}}$ et $S_k=\cS(I_k;\mu_k)$.
Si le mineur $\mu_k$ est d'ordre $j$, le vecteur $V_{n+1-j}$ 
est, dans l'anneau $(A/I_k)_{S_k}$,
égal à une \coli des vecteurs qui suivent.
\end{itemize}
\end{lemma}
\begin{proof}
Pour le deuxième point, on applique les formules de Cramer.
\end{proof}

\begin{proposition} 
\label{propKrRelPol} Soit $B=A[X_1,\ldots, X_n]$  un anneau de polynômes ($n\geq 1$).
La dimension de Krull relative de l'extension
$B/A$ est égale à $n$. 
Donc si  $\dim\, A \le r$, 
alors $\dim\, B\le r+n+rn$. 
Par ailleurs si  $\dim\, B\le r+n$, alors de $\dim\, A \le r$. On abrège cela en disant que $\dim\, B=n+\dim\,A$. 
\end{proposition}

\begin{proof}
La troisième affirmation résulte du fait que si la suite 
 $(a_1,\ldots ,a_r,X_1,\ldots,X_n)$ est singulière dans $B$, 
alors la suite $(a_1,\ldots ,a_r)$ est singulière dans $A$: on a 
en effet, en prenant pour simplifier $m=n=2$ une égalité dans $B$ de
la forme
$$ a_1^{m_1}a_2^{m_2}X_1^{p_1}X_2^{p_2}+ 
a_1^{m_1}a_2^{m_2}X_1^{p_1}X_2^{p_2+1}Q_4+
a_1^{m_1}a_2^{m_2}X_1^{p_1+1}Q_3+
a_1^{m_1}a_2^{m_2+1}Q_2+
a_1^{m_1+1}Q_1 = 0
$$
En regardant dans le polynôme du premier membre le coefficient de 
$X_1^{p_1}X_2^{p_2}$ on trouve
$$ a_1^{m_1}a_2^{m_2}+ 
a_1^{m_1}a_2^{m_2+1}q_2+
a_1^{m_1+1}q_1 = 0
$$
qui donne le collapsus de $\ov{(a_1,a_2)}$ dans $A$.\\
La deuxième affirmation résulte de la première (cf. théorème 
\ref{thKDP}(1)).\\
La démonstration de la première affirmation en \clama s'appuie de manière 
directe sur le cas de corps. Nous donnons une démonstration \cov qui s'appuie 
également sur la cas de corps (discrets). 
Nous reprenons la démonstration de la proposition \ref{propKrDimetDegTr} et nous 
lui faisons subir une relecture
(le corps $K$ est remplacé par un anneau $A$) qui nous permet de faire 
fonctionner le définition du collapsus au dessus de $A$. 
Soit $(y_1,\ldots,y_{n+1})$ dans $A[X_1,\ldots, X_n]$. 
Considérons une démonstration simple (\cad directement écrite comme une 
démonstration d'algèbre linéaire) du fait que les
$y_i$ sont algébriquement dépendants sur $A$ 
lorsque $A$ est un corps discret.
Par exemple si les $y_i$ sont des polynômes de degré $\le d$, les 
polynômes
$y_1^{m_1}\cdots y_{n+1}^{m_{n+1}}$ avec  $\sum_i m_i\le m$ sont dans 
l'espace vectoriel des polynômes de degre $\le dm$, qui est de dimension 
$\le {dm+n\choose n}$, et ils sont au nombre de ${m+n+1 \choose n+1}$.
Pour une valeur explicite de $m$ on a 
${m+n+1 \choose n+1}>{dm+n\choose n}$ (car un polynôme de degré $n+1$ 
l'emporte sur un polynôme de degré $n$).
On considère désormais que $m$ est fixé à cette valeur.
Rangeons les \gui{vecteurs} $y_1^{m_1}\cdots y_{n+1}^{m_{n+1}}$ 
correspondants
(\cad tels que   $\sum_i m_i\le m$) dans l'ordre lexicographique pour 
$(m_1,\ldots,m_{n+1})$. Nous pouvons nous limiter à ${dm+n\choose n}+1$
vecteurs.
En appliquant le lemme \ref{lemALE}, on obtient que dans chacun des 
anneaux
$(A/I_k)_{S_k}$ un vecteur 
$y_1^{m_1}\cdots y_{n+1}^{m_{n+1}}$ est égal à une \coli des vecteurs 
qui suivent.
Ceci donne, comme dans la démonstration de
la proposition \ref{propKrDimetDegTr} un collapsus, mais cette fois-ci
nous devons ajouter au début et à la fin de la \proel 
$\overline{(y_1,\ldots,y_{n+1})}$ les \gui{hypothèses supplémentaires}:
c'est donc la \proc 
$$
\big((\mu_i)_{i<k},y_1;y_2),(y_2;y_3)\ldots,
      (y_{n-1};y_{n}),(y_{n};y_{n+1},\mu_k)\big)
$$ 
qui s'effondre (pour chaque $k$). En effet, pour le collapsus
d'une \proc on peut toujours passer au quotient par le premier 
des idéaux (ou par un idéal plus petit) et localiser en 
le dernier des \mos (ou en un \mo plus petit).\\
Finalement, tous ces collapsus fournissent le collapsus de 
$\overline{(y_1,\ldots,y_{n+1})}$ au dessus de
$A$ en utilisant la famille finie $(\mu_i)$.
\end{proof}

\section{Going Up et Going Down} 
\label{secGUGD}

Cette section est simplifiée par rapport au rapport technique de 2001.
Ceci grâce à la \dfn simple du going up pour les \trdis donnée dans \cite{ACMC}. Ici nous démontrons l'équivalence en \clama avec les \dfns usuelles à base d'\ideps.
\subsection{Généralités}\label{secGUpTrdi} 

\Subsubsection{Relèvement des \ideps (lying over)}

En \clama on dit qu'un \homo $\alpha: T\to V$ de \trdis \gui{possède la propriété de relèvement des \ideps} lorsque l'\homo dual \hbox{$\Spec\,\alpha:\Spec \,V\to\Spec\, T$} est surjectif, autrement dit lorsque tout \idep de $\Spec\, T$
est l'image réciproque d'un \idep de $\Spec\, V$. Pour abréger on dit
aussi que le morphisme est \gui{lying over}.

Nous donnons maintenant une \dfn \cot pertinente sans utiliser l'\homo dual. L'équivalence en \clama avec la \dfn via les spectres
est établie dans la proposition \ref{propRep2}.

\begin{definition}\label{defiLYO}~
\begin{enumerate}
\item Un \homo $\alpha: T\to V$ de \trdis est dit \emph{lying over} lorsqu'il est injectif. Autrement dit $\alpha$ réfléchit les inégalités: 
$$\hbox{pour } a ,\, b\in T\hbox{ on a : }
 \alpha(a)\leq\alpha(b) \ \Longrightarrow\ a\leq b.$$
\item Un \homo $\varphi: A\to B$ d'anneaux
commutatif est dit \emph{lying over} lorsque
l'\homo  $\Zar\,\varphi:\Zar\, A\to\Zar\, B$ est injectif.
\end{enumerate}
\end{definition}
\index{lying over!morphisme ---}

\begin{remark} \label{remLY} 
{\rm    On a aussi les formulations équivalentes suivantes pour les morphismes lying over.
\begin{itemize}
\item Pour les \trdis:
\begin{itemize}
\item Pour tout $b\in T$, $\alpha^{-1}(\dar \alpha(b)\big)=\dar b$.
\item Pour tout idéal $J$ de $ T$,  
$\alpha^{-1} (\gen{(\alpha(J)}_V)=J$.
\end{itemize}
\item Pour les anneaux commutatifs:
\begin{itemize}
\item Pour tout $x\in A$ et tout \itf $I$ de $ A$ on a
l'implication

\centerline{$\varphi(x)\in\gen{\varphi(I)} \
\Longrightarrow\ x\in\sqrt{I}.$}
\item  Pour tout \itf $I$ de $ A$ on a
$\varphi^{-1}(\gen {\varphi(I)}) \subseteq \sqrt{I}$.
\item   Pour tout idéal $I$ de $ A$ on a
$\varphi^{-1}\big(\sqrt{\gen{\varphi(I)}}\big)=\sqrt{I}$. 
\end{itemize}
\end{itemize}
}\end{remark}

\begin{lemma} \label{lemLYO}
Soit $\varphi:T\to V$ un morphisme de \trdis et $S$ un système générateur de $T$. Pour que $\varphi$ soit lying over, il faut et suffit que pour tous
$a_1,\dots,a_n,b_1,\dots,b_m\in S$ soit satisfaite l'implication
$$
\varphi(a_1),\dots,\varphi(a_n)\vdi V \varphi(b_1),\dots,\varphi(b_m)
\;\;\Rightarrow\;\;
a_1,\dots,a_n\vdi T b_1,\dots,b_m
$$ 
\end{lemma}
%
\begin{proof}
Soient $a,b\in T$, on écrit $a=\Vu_i \Vi A_i$, $b=\Vi_j \Vu B_j$
pour des parties finies $A_i$ et $B_j$ de $S$.
Puisque $\varphi$ est un morphisme de treillis, on a $\varphi(a)=\Vu_i \Vi \varphi(A_i)$ et $\varphi(b)=\Vi_j \Vu \varphi(B_j)$.
Supposons $\varphi(a)\leq_V \varphi(b)$. On a, pour chaque $i,j$, $\Vi \varphi(A_i)\leq_V \Vu\varphi(B_j)$, \cad 
$\varphi(A_i)\vdi V \varphi(B_j)$. Par hypothèse on obtient $A_i\vdi T B_j$, \cad $\Vi A_i\leq_T \Vu B_j$
pour chaque $i,j$. Donc $a\leq _T b$. 
\end{proof}
%

\Subsubsection{Montée (going up)} 
\label{subsecGu}
En \clama on dit qu'un \homo $\alpha: T\to V$ de \trdis \emph{poss\`ede la \prt de montée pour les chaînes d'\ideps}, ou plus simplement qu'il est \emph{going up} lorsque  la \prt suivante est satisfaite.

\emph{Si $ Q _1\in\Spec\, V$ et $\alpha^{-1}( Q _1)= P _1$, toute chaîne $P_1\subseteq\cdots\subseteq P_n$ d'\ideps de $T$ 
est l'image réciproque d'une  chaîne $Q_1\subseteq\cdots\subseteq Q_n$ d'\ideps de $V$.
}

\smallskip Naturellement on peut se limiter au cas $n=2$.
On voit alors que la \dfn peut se relire comme suit: \emph{Si $ Q \in\Spec V$ et 
$(\Spec\,\alpha)(Q)=P$, le morphisme $\Spec\,\alpha':\Spec( V/ Q )\to\Spec( T/ P )$ est surjectif.} Ceci nous ramène au cas d'un morphisme lying over.

Les mêmes \dfns sont utilisées pour les morphismes d'anneaux commutatifs.

Nous montrons maintenant en \clama une proposition qui permet de donner une
\dfn \cov satisfaisante du going up.

\begin{proposition} \label{propGu} \emph{(Going up versus lying over)}
\\
Pour un morphisme $\alpha:T\to V$ de \trdis \propeq
\begin{enumerate}
\item  Pour tout \idep $Q$ de $V$, en notant $P=\alpha^{-1}(Q)$ le morphisme 
   $$\alpha':T/(P=0)\to V/(Q=0)$$  
est injectif (i.e., lying over).
\item  Pour tout idéal $I$ de $V$, avec $J:=\alpha^{-1}(I)$, le morphisme 
$\alpha':T/(J=0)\to V/(I=0)$ est injectif. 
\item  Pour tout $i\in V$, avec $J=\alpha^{-1}(\dar i)$ le morphisme 
$\alpha':T/(J=0)\to V/(i=0)$ est injectif. 
\item  Pour tous $a,b\in T$ et $i\in V$ on a
$$
\alpha(a)\vdi V \alpha(b),\, i \quad\Longrightarrow\quad\exists j\in T \;\;\; a \vdi T b ,\, j \;\;\; \hbox{et}\;\;\; \alpha(j)\leq_V i.
$$
\end{enumerate} 
\end{proposition}
%
\begin{proof} Notons que dans le point \emph{3}, $j\in \alpha^{-1}(\dar i)$ signifie  $\alpha(j)\leq_V i$.
Le point \emph{4} est donc simplement la traduction du point \emph{3} en tenant compte de la description de l'inégalité~\hbox{$a\leq b$} dans un quotient par un idéal: $a\vdi{S/(j=0)}b$ équivaut à $a\vdi{S}b,\,j$ (point \emph{3} de la proposition~\ref{propIdealFiltre}).
\\
Le point \emph{3} est un cas particulier du point \emph{2}. Pour l'implication réciproque si $\alpha(a)\vdi{V/(I=0)}\alpha(b)$, il y a un $i\in I$ tel que
$\alpha(a)\vdi{V/(i=0)}\alpha(b)$. 
\\
Le point \emph{1} est un cas particulier du point \emph{2}.
Il reste à prouver \emph{1} $\Rightarrow$ \emph{4}, ce qui ne peut se faire qu'avec des arguments non \cofs. 

\smallskip \noindent \emph{Premier argument, avec Zorn et tiers exclu.}
Les \elts $a,b,i$ satisfaisant $\alpha(a)\vdi V \alpha(b),\,i$ sont fixés et on cherche un $j$ convenable. Par hypothèse (point \emph{1}) pour tout \idep~$Q$ de $V$ contenant $i$, on a un $j_Q\in T$ vérifiant: $a\vdi T b,\,j_Q$ et $\alpha(j_Q)\in Q$. Soit $F$ le filtre de~$V$ engendré par les~$\alpha(j_Q)$. Si $F$ contient $i$, on obtient un nombre fini de $j_{Q_k}$ dont la borne inférieure  $j$ satisfait~\hbox{$\alpha(j)\leq i$} et $a\vdi T b,\,j$ (par distributivité): on a gagné.
Si ce n'est pas le cas, une zornette nous donne un filtre $G$ maximal parmi les filtres contenant~$F$ et ne contenant pas $i$, et ce filtre est premier par construction. D'où un contradiction en considérant l'\idep $V\setminus G$
qui contient~$i$ mais aucun des $\alpha(j_Q)$.  

\smallskip \noindent \emph{Deuxième argument, avec le théorème de complétude.} On considère une théorie dynamique avec deux sortes. Une pour  $T$ et une pour $V$. On a un symbole de fonction $\alpha$ pour un morphisme de $T$ dans $V$ et un prédicat unaire $Q(x)$ sur $V$ pour $x\in Q$. Les constantes sont les éléments de $T$ et $V$. Pour $a,b,i$ fixés satisfaisant 
$\alpha(a)\vdi V \alpha(b),\,i$, on considère la formule $F$: $\exists j\;(\alpha(j)\leq i,\, a\vdi T b,\,j)$. Puisqu'on suppose que \emph{1} est satisfait
la formule $F$ est valide dans tous les modèles de la théorie, car il suffit qu'elle soit valide pour les modèles donnés par $T,\,V,\,\alpha$ et un \idep arbitraire de~$V$). Par le théorème de complétude, la formule $F$ est démontrable.
Mais comme la théorie est sans axiomes existentiels, une formule existentielle $F= \exists j \,G(j)$ n'est démontrable que si est démontrable une disjonction $G(j_1) \hbox{ ou } \dots \hbox{ ou }G(j_\ell)$ pour des termes clos $j_k$ (nécessairement prouvablement égaux à des \elts de $T$). La borne inférieure $j$ de ces $j_k$ nous donne la conclusion voulue.   
\end{proof}
%

\begin{definition}\label{defiGoingup}~
\begin{enumerate}
\item Un \homo $\alpha:T\to V$ de \trdis est dit \emph{going up}
lorsque  pour
\hbox{tous $a,b\in T$} et $y\in V$ on a
$$
\alpha(a)\vdi V \alpha(b), y \quad\Longrightarrow\quad\exists x\in T\;\;\; (a \vdi T b , x\hbox{ et }\alpha(x)\leq y).
$$

\item Un \homo $\varphi: A\to B$ d'anneaux
commutatifs est dit \emph{going up} lorsque
l'\homo  $\Zar\,\varphi:\Zar\, A\to\Zar\, B$ est going~up.
\end{enumerate}
\end{definition} 

Notons que dans le point \emph{1} l'implication réciproque est toujours satisfaite.

Nous donnons maintenant la version \gui{\proc} de la notion usuelle de going up formulée en termes de relèvement d'une chaîne d'\ideps en \clama. La longueur 1 est suffisante.

\begin{proposition} 
\label{propGUclass} {\em (Going Up, version \procs)} 
Soit $\alpha:T\to V$ un morphisme de \trdis.  
Soit $\cP_0$ un  \proi saturé de~$V$, $\cQ_0=\alpha^{-1}(\cP_0)$, $\cQ_1$ 
un \proi  de $T$ tel que $\cQ_0\subseteq \cQ_1$, et $\cP_1=\alpha(\cQ_1)$.   \Propeq
\begin{enumerate}
\item  Dans toute situation décrite précédemment, la \proc 
$(\cP_0,\cP_1)$ 
s'effondre dans~$V$ \ssi la \proc $(\cQ_0,\cQ_1)$ s'effondre dans $T$.
\item  Le morphisme $\alpha$ est going up.
\end{enumerate}
Dans ce cas le saturé de $\cP_1$ est réfléchi dans le saturé de $\cQ_1$.
\end{proposition}
Si $\cP_0=\wh{P_0}$ et $\cQ_1=\wh{Q_1}$ avec $P_0$ et $Q_1$ des \ideps, le \nst formel \ref{th.nstformelTreil} avec le point~\emph{1} ci-dessus montrent qu'on retrouve bien la \dfn classique du going up.
\begin{proof}
Nous noterons
$\cP_0=(J_0,W_0),\; \cQ_0= (I_0,U_0),\; \cQ_1= (I_1,U_1) \hbox{ et }
\cP_1= (J_{1},W_{1})$. 

\smallskip \noindent 
\emph{1} $\Rightarrow$ \emph{2}. On considère $a,b\in T$ et $y\in V$. On prend pour $\cP_0=(\dar y,\emptyset)$ et $\cQ_1=(I_0,b;a)$. On a $I_0=\so{x\in T \,; \,\alpha(x)\leq y}$.
Si $(\cP_0,\cP_1)$ s'effondre, on a un $z\in V$ tel que $z\vdi V y$ \hbox{et $\alpha(a)\vdi V \alpha(b), z$}. Cela revient à dire que $\alpha(a)\vdi V \alpha(b), y$. Sous cette hypothèse on sait que $(\cP_0,\cP_1)$ s'effondre dans $T$, ce qui signifie qu'il existe $x\in T$ tel que $\alpha(x)\leq y$ et $a\vdi T b, x$. 

\smallskip \noindent 
\emph{2} $\Rightarrow$ \emph{1}. On peut supposer \spdg que $I_1$ est un idéal et $U_1$ un filtre.
On suppose que $(\cP_0,\cP_1)$ s'effondre, on doit montrer que $(\cQ_0,\cQ_1)$ s'effondre.  \\
On a $z\in V$, $w\in W_0$, $j\in J_0$,  $a\in U_1$ et $b\in I_1$ tels que $z,w\vdi V j$ et $\alpha(a) \vdi V \alpha(b), z$. La première relation donne $z\in J_0$ parce que $\cP_0$ est saturé. Puique $\alpha$ est going up, on a un $x\in T$ tel que $\alpha(x)\leq z$ (donc $x\in I_0$) et $a \vdi T b,x$. Or $b$ et $x\in I_1$, donc $\cQ_1$ s'effondre. 
 \end{proof}

On vient de traiter une \prolo 1 dans $T$. Le corollaire \ref{corthColSimT2} permet d'obtenir le résultat analogue au précédent pour  une \prolo $\ell$ arbitraire dans $T$.

\Subsubsection{Going down} 
\label{subsecGd}

Il s'agit de la notion opposée au going up, obtenue en renversant l'ordre.
On prend donc la définition \ref{defiGoingdown}, justifiée par la proposition
\ref{propGd} obtenue en renversant l'ordre dans la proposition
\ref{propGu}.

\begin{proposition} \label{propGd} \emph{(Going down versus lying over)}
\\
Pour un morphisme $\alpha:T\to V$ de \trdis \propeq
\begin{enumerate}
\item Pour tout filtre premier $Q$ de $V$, en notant $P=\alpha^{-1}(Q)$ le morphisme 
$$\alpha':T/(P=1)\to V/(Q=1)$$ est injectif (i.e., lying over).
\item Pour tout filtre $I$ de $V$, avec $J:=\alpha^{-1}(I)$, le morphisme $\alpha':T/(J=1)\to V/(I=1)$ est injectif. 
\item Pour tout $i\in V$, avec $J=\alpha^{-1}(\dar i)$ le morphisme $\alpha':T/(J=1)\to V/(i=1)$ est injectif. 
\item Pour
tous $a,b\in T$ et $i\in V$ on a
$$
\alpha(a), y\vdi V \alpha(b) \quad\Longrightarrow\quad\exists x\in T\;\;\; (a , x \vdi T b\hbox{ et }\alpha(x)\geq y).
$$
\end{enumerate} 
\end{proposition}

\begin{definition}\label{defiGoingdown}~
\begin{enumerate}
\item Un \homo $\alpha:T\to V$ de \trdis est dit \emph{going down}
lorsque  pour
\hbox{tous $a,b\in T$} et $y\in V$ on a
$$
\alpha(a), y\vdi V \alpha(b) \quad\Longrightarrow\quad\exists x\in T\;\;\; (a , x \vdi T b\hbox{ et }\alpha(x)\geq y).
$$

\item Un \homo $\varphi: A\to B$ d'anneaux
commutatifs est dit \emph{going down} lorsque
l'\homo  $\Zar\,\varphi:\Zar\, A\to\Zar\, B$ est going~down.
\end{enumerate}
\end{definition} 

Naturellement on obtient la version \proc du going down en renversant l'ordre dans la proposition \ref{propGUclass}.

\subsubsection*{Conséquences pour la dimension de Krull} 
\label{subsecGuGdKrull}

\begin{theorem} \label{thLYGUKdim}
Si un morphisme $\alpha:T\to V$ de \trdis est lying over et going up
(ou bien lying over et going down) on a 
$\dim(T)\leq \dim(V)$. 
\end{theorem}
%
\begin{proof}
La démonstration est simple et naturelle avec les définitions \covs adoptées ici. Voir \cite[section XIII-9]{ACMC}.
\end{proof}
%

\Subsubsection{Incomparabilité} 
\label{subsecIncomp}
En \clama on dit qu'un \homo $\alpha: T\to V$ de \trdis (ou d'anneaux commutatifs) \emph{poss\`ede la \prt d'incomparabilité} lorsque les fibres de l'\homo dual
$\Spec\,\alpha:\Spec\, V\to\Spec\, T$ sont constituées d'\elts
deux \`a deux incomparables. 
Autrement dit, pour $ Q_1$
et $ Q_2$ dans~$\Spec\, V$, si $\alpha^{-1}( Q_1)=\alpha^{-1}( Q_2)$ \hbox{et $ Q_1\subseteq Q_2$}, alors $Q_1=Q_2$.

La \dfn \cov correspondante est que le morphisme $ T\to V$ est zéro-dimensionnel.

 La principale conséquence de la situation d'incomparabilité
pour un \homo $\alpha: T\to V$ est le
fait que $\dim (V)\leq \dim (T)$. Ceci est un cas particulier du théorème \ref{thDim1}.

\subsection{Cas des anneaux commutatifs}\label{secGUpAC} 

Nous donnons deux ou trois énoncés caractéristiques. Les démonstrations et un traitement beaucoup plus complet se trouve dans \cite[section~\hbox{XIII-9}]{ACMC}. 
\Subsubsection{Lying over} 

En application de la remarque \ref{remLY} on a la caractérisation suivante.

\begin{lemma} 
\label{lemGU1} Soit $\varphi:A\to B$ un morphisme d'anneaux commutatifs. 
\Propeq
\begin{enumerate}
\item Le morphisme $\varphi$ est lying over.
\item Pour tout idéal  $I$ de~$A$  et tout $x\in A$, on a:
$\, \varphi(x)\in \varphi(I)B\, \Rightarrow \, x\in\sqrt[A]{I}.
$  
\end{enumerate}

\end{lemma}

\begin{lemma} 
\label{corLyO} 
{\rm (Le lying over usuel)} Soit $A\subseteq B$ des anneaux commutatifs avec $B$ 
entier sur~$A$. Le morphisme $A\to B$ est lying over.
\end{lemma}
\begin{proof}
Supposons  $ x\in {IB}$, i.e., $x=\sum j_ib_i$, avec
$j_i\in I,\;b_i\in B$. Les $b_i$ et 1 engendrent une sous-$A$-algèbre $C$ de $B$ qui est un $A$-module 
fidèle et \tf.   Le \polcar de la matrice de la multiplication 
par~$x$ (exprimée sur un système générateur fini du $A$-module~$C$) a donc tous ses 
coefficients (sauf le coefficient dominant) dans $I$.     
\end{proof}

\begin{lemma} \label{lemGuAnFp}
Soit $A\subseteq B$ une $A$-\alg 
fidèlement plate sur $A$. Le morphisme $A\to B$ est lying over. 
\end{lemma}

\Subsubsection{Going Up} 
\label{subsecGU}

\begin{lemma} 
\label{lemGU2} Soit $\varphi:A\to B$ un morphisme d'anneaux commutatifs. 
\Propeq
\begin{enumerate}
\item Le morphisme $\varphi$ est going up (\dfn \ref{defiGoingup}).
\item Pour tout idéal $I$ de $B$, en notant $J=\varphi^{-1}(I)$, le morphisme
$\varphi_I:A/J\to B/I$ obtenu par factorisation est lying over. 
\item Même chose avec les $I$  de type fini.
\end{enumerate}
En \clama, on peut se limiter aux \ideps dans le point 2.
\end{lemma}
%
\begin{proof}
Cela résulte du résultat  analogue pour les \trdis (proposition \ref{propGu}) et de la comparaison des idéaux de $A$ et de $\Zar\,A$ (fait \ref{factQuotients} et proposition \ref{propZar2}). 
\end{proof}

En application du point \emph{2} et du lying over pour les extensions entières on obtient immédiatement le going up pour les extensions entières.
\begin{corollary} 
\label{corLy1} 
{\rm (Un going up classique)} Soit $A\subseteq B$ des anneaux commutatifs avec $B$ 
entier sur $A$. Alors le morphisme $A\to B$ est going up. 
Donc $\dim\, A\leq \dim\, B$, et vu la proposition~\ref{propDim2}, $\dim\, A= \dim\, B$.
\end{corollary}

Un corolaire du résultat précédent et du théorème \ref{thKDP} 
est le théorème suivant, qui nous dit que la dimension de Krull d'une 
algèbre \pf sur un corps discret est bien celle que nous donne la mise 
en position de Noether.
\begin{theorem} 
\label{thKDP2} Soit $K$ un corps discret,  $I$ un idéal \tf de l'anneau 
$K[X_1,\ldots,X_\ell]$ et $A$ l'algèbre quotient $K[\underline{X}]/I$. 
Si $1\notin I$ une mise en position de Noether de l'idéal $I$ 
fournit un entier $r$  et des éléments $y_1$, \ldots, $y_r$ de $A$ qui 
sont algébriquement indépendants sur $K$ et tels que $A$ est
un module \tf sur $K[y_1,\ldots,y_r]$. Alors~$\dim\,A=r$.
\end{theorem}

Notons que le test $1\in I\,?$ et la mise en position de Noether sont 
réalisées par des algorithmes explicites.

\Subsubsection{Going Down} 
\label{subsecGDAC}

\begin{lemma}\label{lem1Gdown}
\label{lemGD2} Soit $\varphi:A\to B$ un morphisme d'anneaux commutatifs. 
\Propeq
\begin{enumerate}
\item Le morphisme $\varphi$ est going down (\dfn \ref{defiGoingdown}).
\item Pour tous $b$, $a_1$, \ldots, $a_q\in A$ et $y\in B$ tels que~$\varphi(b)y\in\sqrt[B]{\gen{\varphi(a_1,\dots,a_q)}}$,
il existe 
$x_1$, \dots, $x_p\in A$ tels que:
$$
\gen{bx_1,\dots,bx_p}\subseteq  \sqrt[A]{\gen{a_1,\dots,a_q}} \;\hbox{ et }\;y\in\sqrt[B]{\gen{\varphi(x_1), \dots, \varphi(x_p)}}.
$$
\item (en \clama) Pour tout \idep $P$ de $B$, \hbox{avec $Q=\varphi^{-1}(P)$}, le morphisme $A_Q\to B_P$ obtenu par passage au quotient est lying over.
\end{enumerate}
\end{lemma}
\begin{proof}
L'équivalence des points \emph{1} et \emph{3} en \clama résulte du résultat  analogue pour les \trdis (proposition \ref{propGd}), de la comparaison des filtres de $B$ et de $\Zar\,B$ (fait \ref{factLocalises} et proposition \ref{factFZ}) en notant que tout filtre premier de $\Zar\,B$ provient d'un filtre premier de $B$.
\\
Pour l'équivalence des points \emph{1} et \emph{2}, dans la définition
nous avons remplacé  un \elt arbitraire $\sqrt[A]{\gen{\underline b}}$ de~$\Zar A$ et un
\elt arbitraire $\sqrt[B]{\gen{\underline y}}$ de $\Zar B$ par des générateurs~$\sqrt[A]{\gen{b}}$ et~$\sqrt[B]{\gen{y}}$.
Comme les générateurs   engendrent $\Zar A$ et $\Zar B$
par sups finis, les r\`egles de distributivité impliquent que la restriction \`a ces générateurs
est suffisante.
\end{proof}
\hum{On aimerait bien avoir un énoncé \cof similaire au point \emph{3} en remplaçant le filtre premier complémentaire de $P$ par un filtre arbitraire de l'anneau $B$, mais je n'ai pas réussi à faire la chose.}

\begin{theorem} 
\label{thGDplat} {\em (Going Down)} 
Soit $A\subseteq B$ des anneaux. Le morphisme $A\to B$ est going down dans les deux cas suivants. 
\begin{enumerate}
\item $B$ est plat sur $A$.
\item $B$ est intègre et entier sur $A$, et $A$ est intégralement clos.
\end{enumerate}
\end{theorem}


\newpage
\markboth{Annexe}{Annexe}

\setcounter{section}{1}\setcounter{subsection}{0}
\setcounter{theorem}{0}\def\thesection{\Alph{section}}  
\addcontentsline{toc}{section}{Annexe: Complétude, compacité, LLPO et théories géométriques}
\section*{Annexe A: Complétude, compacité, LLPO et théories géométriques}\label{Annexe}
\subsection{Théories propositionnelles et modèles}

Considérons un ensemble $V$ fixé de {\em propositions atomiques}, ou
{\em variables de propositions}, ou \emph{atomes}. 

Une \emph{proposition} $\phi,\psi,\dots$
est un objet syntactique construit sur les atomes  $p,q,r\in V$ avec les connecteurs logiques usuels
$$ 0,\;\;1,\;\;\phi\wedge \psi,\;\;\phi\vee\psi,\;\;\phi\rightarrow\psi,\;\;
\neg\phi
$$

Notons $P_V$ l'ensemble des propositions. Soit $F_2$
l'algèbre de Boole à deux éléments.
Une {\em valuation}
est une fonction $v\in F_2^V$ qui assigne une valeur de vérité à  
chaque proposition atomique. Une telle valuation s'étend en une application
$P_V\rightarrow \{0,1\},\;\phi\longmapsto v(\phi)$
de manière naturelle\footnote{En utilisant les tables de vérité des connecteurs.}. 

Une {\em théorie propositionnelle}~$\cT$ est un sous-ensemble de $P_V$.
Un {\em modèle} de $\cT$ est une  valuation $v$ qui vérifie $v(\phi)=1$
pour tous \hbox{les $\phi\in \cT$}.

Plus généralement, étant donnée une algèbre de Boole $B$ on peut définir
une $B$-valuation comme une fonction $v\in B^V$. Une telle valuation s'étend 
également en une application
$P_V\rightarrow B$, $\phi\longmapsto v(\phi)$ de manière naturelle.  
Un {\em $B$-modèle} de $\cT$ est une  valuation $v$ qui vérifie $v(\phi)=1$
pour tous \hbox{les $\phi\in \cT$}. La notion usuelle de modèle est un cas particulier en prenant pour $B$ l'\alg $F_2$.

Pour une théorie $\cT$ donnée il existe  une \agB 
universelle pour laquelle~$\cT$ est un modèle. On l'appelle l'\emph{algèbre de Lindenbaum de $\cT$}. C'est une algèbre de Boole définie par générateurs et relations. Les générateurs sont les éléments de $V$ et chaque \elt $\phi$ de~$\cT$ fournit la relation $v(\phi)=1$. La théorie $\cT$ est formellement consistante \ssi son algèbre de Lindenbaum est non triviale. 

\subsection{Le théorème de complétude}

\begin{theorem}
(théorème de complétude) Une théorie $\cT$ est formellement consistante \ssi elle admet un modèle.
\end{theorem}

Ceci est le théorème de complétude pour la logique des propositions, démontrable en \clama.

Un tel théorème est fortement relié au programme de Hilbert, qui peut être vu comme une tentative de remplacer la question de  l'existence d'un modèle 
pour une théorie par le fait purement formel que la théorie n'est pas 
contradictoire.

Soit $B$ l'algèbre de Lindenbaum de $\cT$. Pour prouver la complétude
il suffit de trouver un morphisme $B\rightarrow F_2$ quand $B$ n'est pas
triviale, ce qui est la même chose que trouver un \idep dans une algèbre de Boole non triviale.

Notez que l'existence est claire quand $B$ est finie, car elle est isomorphe à une algèbre~$F_2^k$. Dans ce cas le théorème de complétude est direct.

\subsection{Le \tcg}

Le théorème de complétude peut être vu comme une conséquence du résultat fondamental suivant. 

\begin{theorem}
(\tcg) Soit $\cT$ une théorie.
Si tous les sous-ensembles finis de~$\cT$ admettent un modèle, il en va de même pour $\cT$.
\end{theorem}

En fait cet énoncé est équivalent au théorème de complétude. En effet une théorie est formellement consistante \ssi toutes ses sous-théories finies sont formellement consistantes. Et toute théorie finie consistante admet un modèle.

Une preuve simple générale du \tcg consiste à munir $W:=\{0,1\}^V$ de la topologie produit et à remarquer que l'ensemble des modèles de $\cT$ est un sous espace fermé. Le théorème est alors un corolaire de la compacité de l'espace $W$ si l'on exprime la compacité (en \clama) sous la forme: si une famille de sous-ensembles fermés de $W$ a des intersections finies non vides, alors son intersection est non vide. 

\subsection{LPO et LLPO}

Lorsque $V$ est dénombrable (\cad discret et énumérable) on a l'argument alternatif suivant. On écrit $V = \{p_0,p_1,\dots \}$ et l'on construit par \recu une valuation partielle~$v_n$ sur $\{p_i~|~i<n\}$ de sorte que: 
\begin{itemize}
\item $v_n$ étend $v_{n-1}$,
\item chaque sous-ensemble fini de $\cT$ admet un modèle qui étend un $v_n$ 
\end{itemize}
Pour définir $v_{n+1}$ on essaie d'abord $v_{n+1}(p_n) = 0$.
Si cela échoue il y a une partie finie de $\cT$ telle que chacun de ses modèles qui étend $v_n$ satisfait $v(p_n) = 1$ et l'on prend $v_{n+1}(p_n) = 1$.

Le caractère non effectif de l'argument est contenu dans le choix de $v_{n+1}(p_n)$, qui demande de donner une réponse globale à un ensemble infini de questions élémentaires.

Maintenant supposons aussi que nous pouvons énumérer l'ensemble infini $\cT$.
Nous pouvons alors construire une suite croissante de parties finies de $\cT$,
$K_0\subseteq K_1\subseteq \dots$ qui recouvre~$\cT$.
Si nous avons construit une valuation $v_n$ qui s'étende en un modèle de chaque~$K_j$, pour définir~$v_{n+1}(p_n)$, nous devons donner une réponse globale à la question: est-ce que tous les $K_j$ ont un modèle $w$ qui prolonge~$v_{n}$ et qui vérifie $w(p_n)=0$~?  
Pour chaque $j$, il s'agit d'une question élémentaire qui a une réponse claire. 

Plus précisément définissons $g_{n}:\N\rightarrow \so{0,1}$ de la manière suivante: $g_{n}(j)=0$ s'il y a un modèle
$v_{n,j}$ de $K_{j}$ qui étend  $v_{n}$ avec $v_{n,j}(p_{n})=0$,
sinon~\hbox{$g_{n}(j)=1$}. Si $g_{n}(j)=1$ et \hbox{si $\ell\geq j$}
tous les modèles~$v_{n,\ell}$
de~$K_{\ell}$ qui étendent~$v_{n}$ satisfont~\hbox{$v_{n,\ell}(p_{n})=1$}.
Ainsi la fonction~$g_n$ est croissante.
Nous devons poser $v_{n+1}(p_n)=1$ si $g_{n}(j)=1$ pour un $j$.
Si $g_{n}(j)=0$ pour tout~$j$, nous pouvons poser $v_{n+1}(p_n)=0$.
Dans les deux cas nous savons que pour tout~$j$, il y a un modèle de $K_j$
qui étend $v_{n+1}$. La \recu fonctionne! 
 Ainsi nous \gui{construisons} par \recu une suite infinie de modèles partiels $v_n$ en utilisant à chaque étape le principe non constructif~LPO de Bishop (Least Principle of Omniscience):
\begin{itemize}
\item  [$\bullet$] \emph{pour toute fonction $f:\N\rightarrow \{0,1\}$,
\hbox{on a $f=0$} \hbox{ou $\exists j\in\N\;f(j)= 1$}.}
\end{itemize}
Ce principe est appliqué
à chaque étape $n$ à la fonction $g_n$. 

En fait nous pouvons modifier 
légèrement l'argument et utiliser seulement une combinaison du Choix Dépendant et du principe de Bishop LLPO (Lesser Limited
Principle of Omniscience), qui est connue pour être strictement plus faible que LPO: 
\begin{itemize}
\item  [$\bullet$] \emph{étant données deux fonctions croissantes
$g,h:\N\rightarrow \{0,1\}$ telles que, pour chaque $j$ on ait
$$g(j) = 0 \hbox{ ou } h(j) = 0$$
alors nous avons $g=0$ ou $h=0$.} 
\end{itemize}
En effet définissons
$h_{n}:\N\rightarrow \{0,1\}$ de manière symétrique:
 $h_{n}(j)=0$ s'il y a un modèle~$v_{n,j}$ de $K_{j}$ qui étend  $v_{n}$ avec $v_{n,j}(p_{n})=1$,
sinon $h_{n}(j)=1$. Les fonctions $g_n$ et $h_n$ sont croissantes. Par \hdr nous avons
$g_n(j) = 0 \vee h_n(j) = 0$ pour tout $j$. Donc par LLPO, nous pouvons définir $v_{n+1}(p_n)=0$
si $g_n=0$ et $v_{n+1}(p_n)=1$ si $h_n=0$.
Néanmoins nous devons utiliser le Choix Dépendant pour pouvoir faire
ce choix une infinité de fois car la réponse
\gui{$g=1$ ou $h=1$} donnée par l'oracle LLPO peut être ambigüe.

Réciproquement, on voit facilement que le \tcg restreint au cas dénombrable implique LLPO.

\subsection{Théories géométriques}

{\em What would have happened if topologies without points
had been discovered before topologies with points, or if
Grothendieck had known the theory of distributive lattices?}  (Gian Carlo Rota \cite{Rota}).


\smallskip 
Une formule est \emph{géométrique} \ssi elle est construite uniquement avec le connecteurs
 $0,1,\phi\wedge\psi,\phi\vee\psi$ à partir des symboles de proposition dans $V$. Une théorie (propositionnelle) $\cT$ est \emph{géométrique} 
 \ssi toutes les formules dans $\cT$ sont du type $\phi\rightarrow \psi$ 
où $\phi$
et $\psi$ sont géométriques.

Les formules d'une théorie géométrique peuvent être vues comme fournissant une \entrel sur $V$ pour un \trdi $L_\cT$ (engendré par $V$)
et l'\agL de $\cT$ n'est rien d'autre que l'\agB
 librement engendrée par le \trdi~$L_\cT$. 
 
Il résulte de la proposition
\ref{propTrBoo} que $\cT$ est formellement consistante \ssi le treillis $L_\cT$ est non trivial.
En outre un modèle de $\cT$ n'est rien d'autre qu'un \elt de~$\Spec(L_\cT)$.

\begin{theorem}
(théorème de complétude pour les théories géométriques) Si  le \trdi $L_\cT$ défini par une théorie géométrique $\cT$ est non trivial, alors
$\cT$ a un modèle.
\end{theorem}

\subsection{Théories dynamiques}

Voir  \cite{clr}. Une théorie propositionnelle dynamique est presque la même chose qu'une théorie géométrique. Elle est un moyen très élémentaire de produire une théorie géométrique $\cT$ à partir d'\emph{axiomes}. 

Elle est définie comme un système de preuves pour 
établir des \emph{règles dynamiques}.
Une règle dynamique est une règle de déduction du type
$$
P \,\vdash\, Q_1 \hbox{ ou } \dots \hbox{ ou } Q_n
$$
où $P$ et chaque $Q_i$ sont des listes dans $V$ (des listes de propositions atomiques).
Une règle sans second membre est lue comme $P \,\vdash\, \bot$. La conclusion $\bot$ signifie l'effondrement de la théorie.

Les axiomes de la théorie sont  des règles dynamiques et le fonctionnement des
preuves consiste à valider d'autres règles dynamiques à partir des axiomes en appliquant les règles de réflexivité, monotonie et transitivité évidentes. 

Si par exemple la règle suivante 
$$
(p_1,p_2,p_3)\; \vdash\; (q_{1,1},q_{1,2}) \hbox{ ou } (q_{2,1},q_{2,2},q_{2,3})
$$
est validée à partir des axiomes par le système de preuves, cela valide 
dans la théorie géométrique correspondante $\cT$ la présence de la formule  
$$
(p_1 \vi p_2 \vi p_3) \im (q_{1,1} \vi q_{1,2}) \vu (q_{2,1} \vi q_{2,2} \vi q_{2,3}).
$$

Une règle est dite algébrique s'il n'y a pas de \gui{ou} dans le second membre.
Pour une théorie algébrique, \cad une théorie où tous les axiomes sont algébriques, les conditions d'effondrement sont souvent faciles à établir.
Par exemple si l'effondrement d'un anneau commutatif\footnote{L'anneau est défini à partir d'un système générateur $G$ par une famille d'hypothèses $p_i=0$ pour des $p_i\in\Z[G]$.} est défini comme le fait d'avoir établi la règle $\,\vdash\,1=0$, l'effondrement signifiera que l'idéal des hypothèses $a=0$ contient $1$. 

\subsubsection*{Collapsus simultané}
Il arrive que l'ajout d'un axiome dans une théorie dynamique ne change pas
le collapsus (l'effondrement de la théorie). On dit que la nouvelle théorie
s'effondre simultanément avec la première.

Les théorèmes de collapsus simultané sont les versions \covs des théorèmes en \clama qui disent que tout modèle de la première théorie s'étend en un modèle de la seconde. Par exemple un anneau commutatif s'effondre dès que le fait de le voir comme un corps algébriquement clos dynamique le fait s'effondrer.
Cela constitue la version \cov de plusieurs théorèmes classiques importants: 1) tout anneau commutatif non trivial possède un idéal premier; 2) tout anneau intègre est un sous-anneau d'un corps; 3) tout corps possède une clôture algébrique. Mais seul le deuxième admet une preuve \cov.

Dans le cas où on ajoute un axiome $P\,\vdash\,Q_1 \hbox{ ou } Q_2$, démontrer le collapsus simultané de la nouvelle théorie revient à démontrer que
si l'ajout séparément de $\vdash\,Q_1$ ou de $\vdash\,Q_2$ produit le collapsus, alors l'ajout de $\vdash\, P$ produit également le collapsus.

Un autre exemple important est la démonstration du caractère inoffensif de l'ajout du tiers exclu dans les théories géométriques. Plus précisément, si pour un $p\in V$
on ajoute un symbole $p'$ (non dans $V$) et si sur le nouveau $V'$ on ajoute les axiomes $\vdash\,p \hbox{ ou }p'$ \hbox{et $p,p'\,\vdash \bot$}, la nouvelle théorie s'effondre simultanément avec la précédente (en fait elle prouve exactement les mêmes règles dynamiques formulées sans l'utilisation de $p'$).

\subsubsection*{Théories dynamiques du premier ordre \cite{clr}}

La généralisation usuelle des théories propositionnelles est donnée par les théories formelles du premier ordre.

Il correspond à cela des théories géométriques du premier ordre et
des théories dynamiques du premier ordre.

Ce qui relève d'une théorie dynamique du premier ordre admet une interprétation \cov évidente.
En outre, comme pour les théories propositionnelles, un théorème  explique le caractère inoffensif de
la logique classique du premier ordre lorsque l'on se situe dans une théorie formelle dont tous les axiomes sont des formules géométriques (\cite{clr}): \emph{lorsqu'on passe d'une théorie dynamique à la théorie formelle du premier ordre avec logique classique correspondante, les règles dynamiques prouvables (formulées dans la théorie dynamique de départ) restent les mêmes}. 

Ce résultat fondamental  est rendu
particulièrement simple par la considération des théories dynamiques du premier ordre. Cela tient à ce que le système de preuve dans une théorie dynamique
relève de calculs purs, sans utilisation de formules plus compliquées que les formules atomiques. Il s'agit en quelque sorte de preuves directes sans utilisation de la logique. 

Le \tcg pour les théories propositionnelles (valable en \clama) implique également le théorème de complétude pour les théories formelles du premier ordre, mais la démonstration (\cov) de cette implication est plus délicate.



\end{document}